\documentclass[11pt,a4paper]{article}

\usepackage[T1]{fontenc}

\usepackage[french,english]{babel}

\usepackage[dvips]{graphicx}
\usepackage{amsmath,amsfonts,amssymb,amsthm,bbm,latexsym,graphicx,color}
\usepackage{amsmath,amssymb,exscale,epsfig,psfrag}
\usepackage{float}

\numberwithin{equation}{section}

\newtheorem{lemme}{Lemma}[section]
\newtheorem{proposition}[lemme]{Proposition}

\newtheorem{remarque}[lemme]{Remark}
\newtheorem{definition}[lemme]{Definition}
\title{{\bf Estimation of the intensity parameter of the  germ-grain Quermass-interaction model when the number of germs is not observed}}

\author{
David {\sc Dereudre} \\
{\footnotesize  
Paul Painlevé Mathematics Institute, University of Lille 1}\\
{\footnotesize  UMR CNRS 8524}\\
{\footnotesize  59655 Villeneuve d'Ascq, France}\\
{\footnotesize e-mail~:~david.dereudre@univ-lille1.fr}\\[7mm]
Frédéric {\sc Lavancier}\\
{\footnotesize  
Jean Leray Mathematics Institute, University of Nantes}\\
{\footnotesize  UMR CNRS 6629}\\
{\footnotesize 2 rue de la Houssinière}\\
{\footnotesize 44322 Nantes Cedex 3,  France.}\\
{\footnotesize e-mail~:~frederic.lavancier@univ-nantes.fr}\\[7mm]
Kate\v{r}ina {\sc Sta\v{n}kov\'{a} Helisov\'{a}}\\
{\footnotesize Czech Technical University in Prague,} 
{\footnotesize Faculty of Electrical Engineering}\\ 
{\footnotesize Department of Mathematics, Technick\'{a}~2}\\
{\footnotesize 16627~Prague~6~-~Dejvice, Czech Republic.}\\
{\footnotesize e-mail~:~helisova@math.feld.cvut.cz}
}

\newcommand{\E}{\ensuremath{{\mathcal{E}}}}

\newcommand{\A}{\ensuremath{{\mathcal{A}}}}

\newcommand{\LL}{\ensuremath{{\mathcal{L}}}}

\newcommand{\R}{\ensuremath{{\mathbb{R}}}}

\newcommand{\Rd}{\ensuremath{{\mathbb{R}^2}}}

\newcommand{\U}{\ensuremath{\mathcal U}}
\renewcommand{\t}{\ensuremath{\theta}}
\newcommand{\x}{\ensuremath{\chi}}

\newcommand{\1}{\ensuremath{\mbox{\rm 1\kern-0.23em I}}}
\renewcommand{\L}{\ensuremath{\Lambda}}
\renewcommand{\l}{\ensuremath{\lambda}}
\newcommand{\oo}{\ensuremath{\omega}}
\renewcommand{\O}{\ensuremath{\Omega}}

\begin{document}
\maketitle

\vspace{-3mm}

\centerline{\textbf{Abstract}}

The Quermass-interaction model allows to  generalise the classical germ-grain Boolean model in adding a morphological interaction between the grains. It enables to model random structures with specific morphologies which are unlikely to be generated from a Boolean model.  The Quermass-interaction model depends in particular on an intensity parameter, which is impossible to estimate from  classical likelihood or pseudo-likelihood approaches because  the number of points is not observable  from a germ-grain set. In this paper, we present a procedure based on the Takacs-Fiksel method which is able to estimate all parameters of the Quermass-interaction model, including the intensity. An intensive simulation study is conducted to assess the efficiency of the procedure and to provide practical recommendations. It also illustrates that the estimation of the intensity parameter is crucial in order to identify the model.  The Quermass-interaction model is finally fitted by our method to P. Diggle's heather dataset.
\noindent

\vspace{5mm}
\noindent
KEY-WORDS: Gibbs Point Process, germ-grain model, Quermass-Interaction Process, Area-interaction Process, Perimeter-interaction Process, Takacs-Fiksel estimator. 
\newpage

\section{Introduction}

Physics, biology or agronomy are often confronted with problems involving complex random sets like liquid-vapour interface structures, micro-emulsions, porous media or propagation domains of plants. Admissible models for these random structures are the germ-grain models. These random set models are constructed from a point process (the germs), where each point is associated with a random set (the grain). The union of grains forms the final random structure. The most popular germ-grain model is certainly the Boolean model of balls. In this model, the grains correspond to balls, whose centres are distributed as a Poisson point process and the radii are independently and identically distributed. 
The probabilistic and statistical properties of the Boolean model  are well known (see for example \cite{Molchanov97, SKM}).  However, considering the Boolean model as a random set (see \cite{Lieshout}), the variety of morphological structures that it generates is limited. 
To reach more realistic morphologies, some Gibbs  modifications of the Poisson point process defining the germs of the Boolean model have  been developed.  The idea is to introduce an interaction (or  Hamiltonian) between the points, where the interaction depends on geometrical features of the associated union of grains. The density of this new process is maximal when the Hamiltonian is minimal. So, the random structures generated by the model tend to minimise the Hamiltonian and this produces  geometrical features that are unlikely to occur for Boolean models.


In the following  the grains are assumed to be balls. A first model, based on the area of the union of grains, has been introduced in \cite{WR}, known as the penetrable sphere model, or Widom-Rowlinson model. This model was not designed to model random sets, but  was introduced  as a marked point process with attraction, to model liquid-vapour phase transitions. The so-called area-interaction process  \cite{BV} is an extension of the Widom-Rowlison model to  both  the attractive and repulsive case. The area-interaction process can be viewed either as a marked point process,  in which case it aims at modelling interactions between points, or as a random set model, which we consider in the present paper. Motivated by  the former point of view, multi-scales area-interaction processes have been considered in \cite{Ambler,Goulard,Picard}. Other extensions, more adapted to the random geometry setting, can be found in  \cite{Kendall, KVB,LMW, Me}. These works consider a Hamiltonian that relies not only on the area but also on other functionals of the union of grains. Hadwiger's theorem \cite{Hadwiger} ensures, under mild conditions, that any function acting on an union of compact convex sets can be decomposed into a linear combination of the Minkowski (or Quermass) functionals.
These functionals correspond in $\R^2$ to the area, the perimeter and the Euler-Poincaré characteristic (number of connected components minus number of holes). Accordingly, any function of the union of balls  in $\R^2$ can be written as a linear combination of these three functionals. A Quermass-interaction model  \cite{Kendall, KVB} corresponds to the choice of a Hamiltonian equal to  some linear combination of the area, the perimeter and the  Euler-Poincaré characteristic of the union of grains. From Hadwiger's theorem, it appears as a very rich model to represent random structures. 

The Quermass-interaction model in $\R^2$ (also called Quermass model for short in the following) is entirely specified by the law of the radii and by four parameters: the three coefficients in the linear combination defining the Hamiltonian, and the intensity parameter of the underlying Poisson point process. The present paper deals with the estimation of these four parameters, while the law of the radii is assumed to be known and is included in the reference measure. Note that it is also possible to consider a parametric law for the radii, e.g. a continuous parametric law or a multi-scale discrete law as in \cite{Ambler,Goulard,Picard}, provided some modifications in the definition of the model, see Remark~\ref{rem_reference_measure}. The estimation of these extra parameters can be conducted with the same procedure as described in this paper, but we do not consider this generalisation in this work.

The main difficulty for inference comes from the nature of the observable data. Since  we use the  Quermass germ-grain model to model  random sets, we assume that we observe the  union of grains only. In particular the germs and the number of balls are not observed and cannot be used in the statistical procedure, which is untypical in estimation problems for Gibbs point processes. Note that this specific issue already occurs for the estimation of the intensity parameter of the classical Boolean model (viewed as a random set). In this case, some explicit estimating equations have been found, that express the intensity parameter in terms of the specific volume of the set, see \cite{Molchanov97}. In the presence of Gibbs interactions as for the Quermass-interaction model, it is well-known that the  computation of macroscopic quantities is intractable, so a similar explicit estimation procedure is not possible.
  
Assuming the intensity parameter of the underlying Poisson point process is known, a maximum likelihood approach has been  investigated in \cite{MH2} for the estimation of the three other parameters of the Quermass interaction. Unfortunately the intensity parameter cannot be estimated by this method due to the unobservability of the number of points. Section~\ref{MLE} gives more details about this procedure and explains the serious consequences of a misspecification of the intensity parameter in practice. In this paper, we estimate all parameters of the Quermass process, including the intensity, via a Takacs-Fiksel procedure (\cite{Takacs, Fiksel, CDDL}). The Takacs-Fiksel contrast function is based on an empirical counterpart of the Georgii-Nguyen-Zessin  equation (\cite{Georgii, NZ}), and depends on the choice of  test functions. In the context of point processes, the performance of the procedure depends on the latter choice. But in our setting of random sets,  the contrast function is  not computable in general due to the unobservability of germs. However, some specific choices of test functions lead to a contrast function that does not depend on the number of points (see Section~\ref{functions}). This particular case of the Takacs-Fiksel procedure allows the estimation of the Quermass model. The purpose of the present paper is to provide these relevant test functions,  to present a simulation study and to give some practical recommendations.
 
As an application, we finally  fit a Quermass model to heather data. The heather dataset was initially analysed by P. Diggle in \cite{Diggle81}, followed by many studies (\cite{Hall85, Hall88, Cressie, Mrkvicka09, MH2}). We show that our model seems to be a better approximation of this heather dataset, both from a visual impression and from a statistical diagnostic inspection.

In Section~\ref{QM}, we introduce the Quermass-interaction model and we recall the fundamental Georgii-Nguyen-Zessin equation. Section~\ref{sec:estimation} presents the estimation procedures. The limitation in using the maximum likelihood approach is explained in Section~\ref{MLE}. Then in Section~\ref{TFsection} we present the general Takacs-Fiksel procedure and its application to the Quermass-interaction model. Practical aspects for the implementation of the procedure are given in Section~\ref{PA}. A simulation study assessing the efficiency of the procedure is presented in Section~\ref{sec:simu}. A fit to the heather dataset is conducted in Section~\ref{sec:heather}. Finally, note that an  appealing alternative estimation procedure would consist of combining the likelihood  and the Takacs-Fiksel approaches, to take advantages of both methods. We have implemented such a mixed procedure (not shown in this paper),  but it turns out  to be very time consuming without being more efficient than the Takacs-Fiksel procedure.

\section{Quermass-interaction model}\label{QM}

\subsection{Notations}

We denote by $\E$ the space $\Rd \times [0,R_0]$ (where $R_0>0$ is a fixed positive real number)  endowed with its natural Euclidean Borel $\sigma$-algebra. It is the space of marked points $(x,R)$ where $x\in\Rd$ is the centre of a ball  and $R$ its radius. We assume for simplicity that the radii are  bounded by $R_0$. This assumption allows to define properly the Quermass model on the whole plane $\Rd$, see Section~\ref{sec:gibbs}, though this restriction is not mandatory (see \cite{Dereudre09}). For any bounded set $\L\subset\Rd$, we denote by $\E_\L:=\L\times [0,R_0]$ the restriction of $\E$ to $\L$.

By definition, a {\it configuration} of points  $\oo$ is a locally finite subset of $\E$, which means that the set $\oo_\L:=\oo\cap\E_\L$ is finite for any bounded set $\L\subset\Rd$. 
The space of all configurations  of points in $\E$ is denoted by  $\O$, while  for any bounded set $\L\subset\Rd$, $\O_\L$ denotes the subspace of configurations in $\E_\L$. 

For $x\in\Rd$, we write for short $x\in\oo$ if there exists $R\in[0,R_0]$ such that $(x,R)\in\oo$. For $(x,R)\in\E$ we write $\oo\cup (x,R)$ instead of $\oo\cup \{(x,R)\}$ and $\oo \backslash (x,R)$ instead of $\oo \backslash \{(x,R)\}$. 

For any configuration $\oo$ we denote by $\U_{\oo}$ its germ-grain representation defined by the following set 
$$ \U_{\oo}:= \bigcup_{(x,R)\in\oo}  B(x,R),$$
where $B(x,R)$ is the closed ball centred at $x$ with radius $R$. 

Let $\mu$ be a reference probability measure on $[0,R_0]$. We denote by $\l$ the Lebesgue measure on $\Rd$ and by $\pi^{\mu}$ the marked Poisson process on $\E$ with intensity measure $\l \otimes \mu$. For every bounded set $\L$, the probability measure $\pi^{\mu}_\L$ denotes the marked Poisson process on $\E_\L$ with intensity measure $\l_\L\otimes \mu$. Recall that the law of the random set $\U_{\oo}$ under the probability measure $\pi^{\mu}$ is nothing else than the standard homogeneous Boolean model with intensity one and  distribution of radii $\mu$. 

\subsection{Quermass-interaction model on a bounded window}

Following Kendall et al. \cite{KVB}, for any configuration  $\oo_\L$ in a bounded window $\L\subset\Rd$, the Quermass interaction  (or Quermass Hamiltonian) is defined by 
\begin{equation}\label{energy}
H^{\theta}(\oo_\L)= \t_1\;  \A(\U_{\oo_\L}) + \t_2\; \LL(\U_{\oo_\L}) +\t_3\;  \x(\U_{\oo_\L}),
\end{equation}   
where $\theta:=(\t_1, \t_2, \t_3)$  is a vector of  real parameters. The functionals $\A$, $\LL$ and $\x$ are the three fundamental Minkowski (or Quermass) functionals: area, perimeter and Euler-Poincaré characteristic (the number of connected components minus the number of holes). 

From Hadwiger's Theorem \cite{Hadwiger},  any additive functional  defined on the space of finite unions of convex compact sets and satisfying some continuity assumption (see \cite{Hadwiger})  can be decomposed as in (\ref{energy}). This universal representation explains the interest of the Quermass interaction for morphological modelling at mesoscopic scale  \cite{LMW,Me}.   

\begin{definition}\label{Quermassfini}
The Quermass point process on a bounded set $\L\subset\Rd$ with parameter $\theta\in\mathbb{R}^3$,  intensity $z>0$ and  distribution $\mu$ of the radii, is the probability measure $P^{z,\theta}_\L$ on $\O_\L$ which is absolutely continuous with respect to the marked Poisson Process $\pi_\L^{\mu}$ with  density 

\begin{equation}\label{dens_quermass} g_\L(\oo_\L; z, \theta)= \frac{1}{Z_\L(z,\theta)}\; z^{n(\oo_\L)} \mathrm{e}^{-H^{\theta}(\oo_\L)},\end{equation}
where $n(\oo_\L)$ denotes the cardinality of $\oo_\L$ and $Z_\L(z,\theta):=\int z^{n(\oo_\L)} \mathrm{e}^{-H^{\theta}(\oo_\L)} \pi_\L^{\mu}(\mathrm{d}\oo_\L)$ is a normalising constant called the partition function.
\end{definition}

\begin{figure}[htbp]
    \setlength{\tabcolsep}{0.1cm} \centerline{
  \begin{tabular}[]{ccc}
    \includegraphics[angle=0,scale=.2]{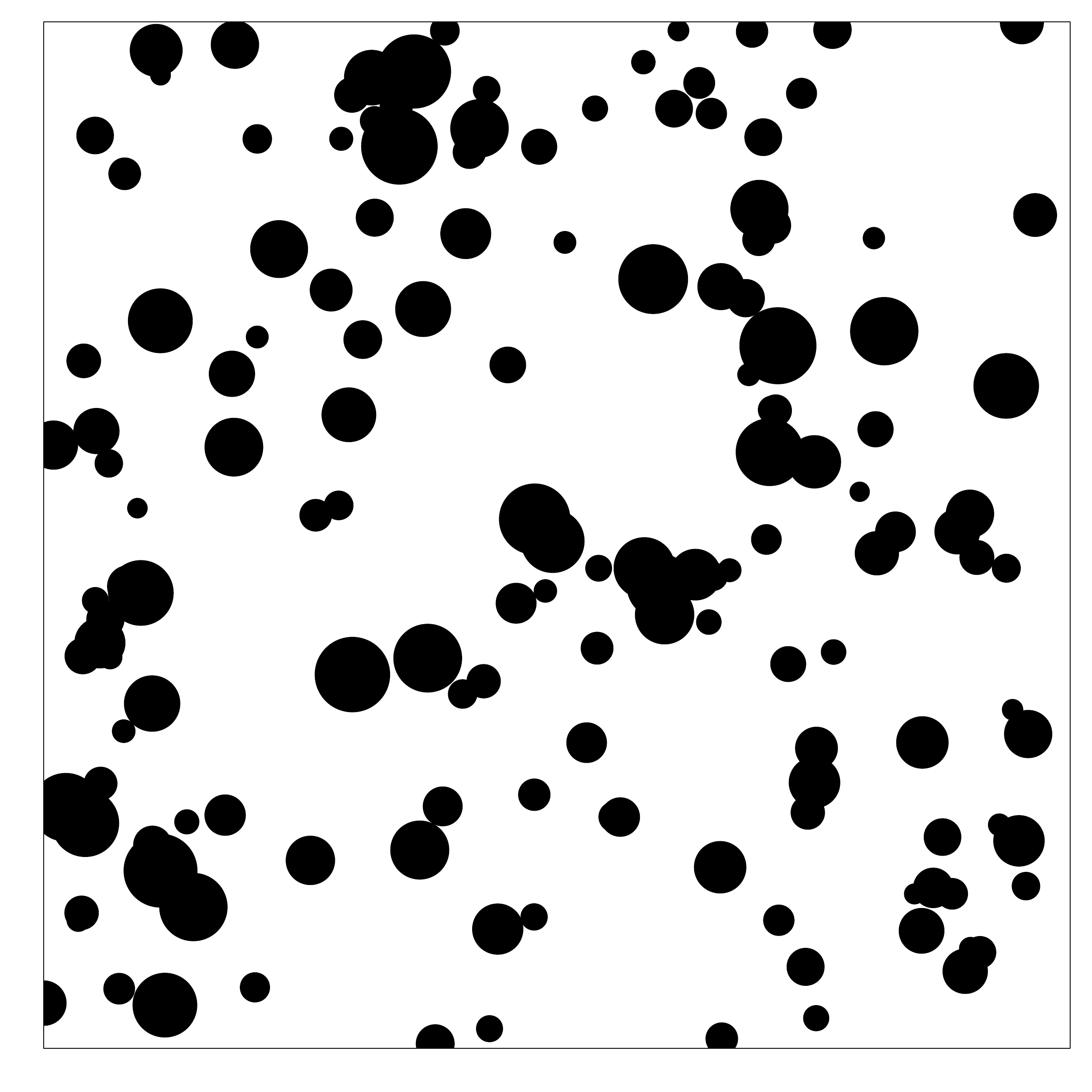}&
  \includegraphics[angle=0,scale=.2]{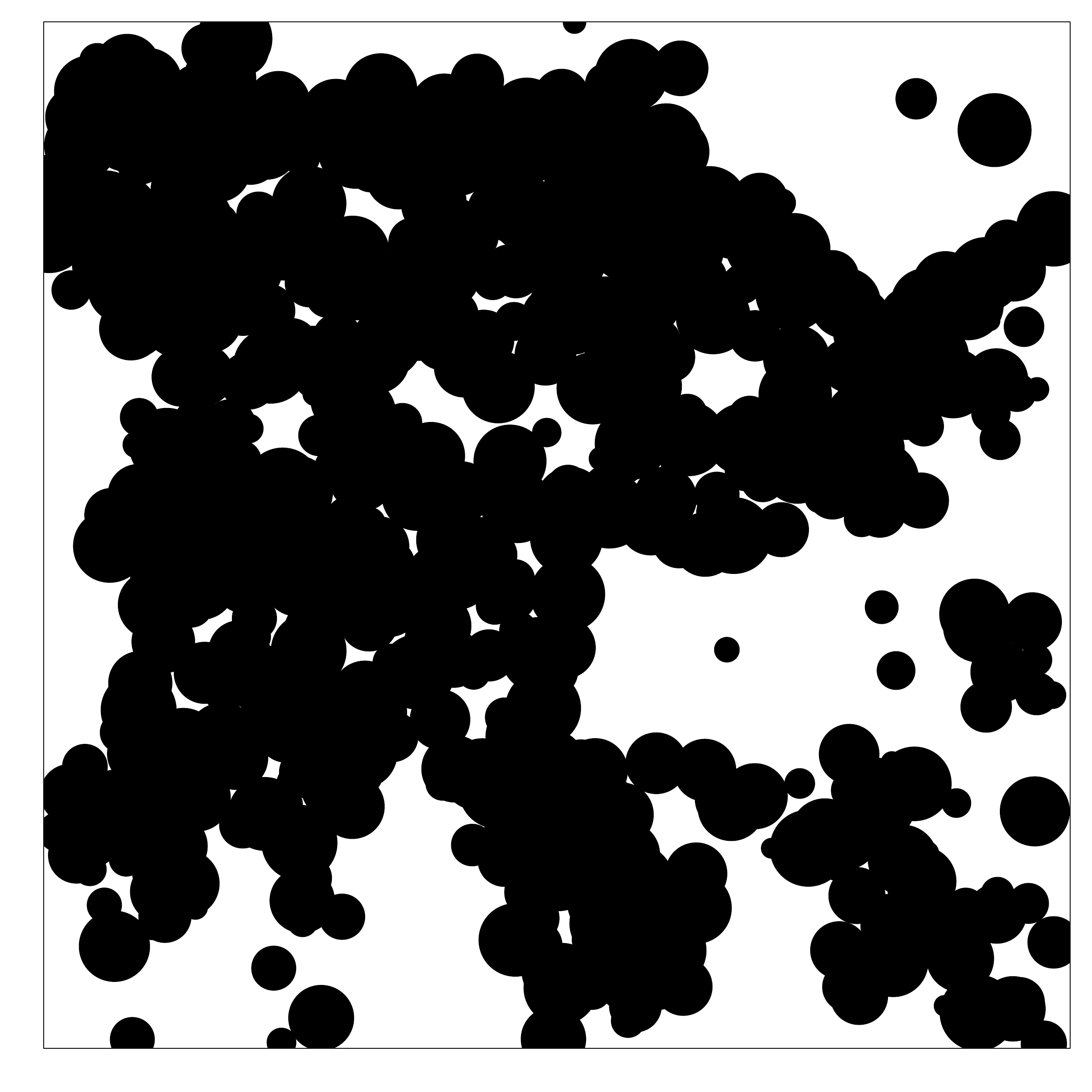}&
   \includegraphics[angle=0,scale=.2]{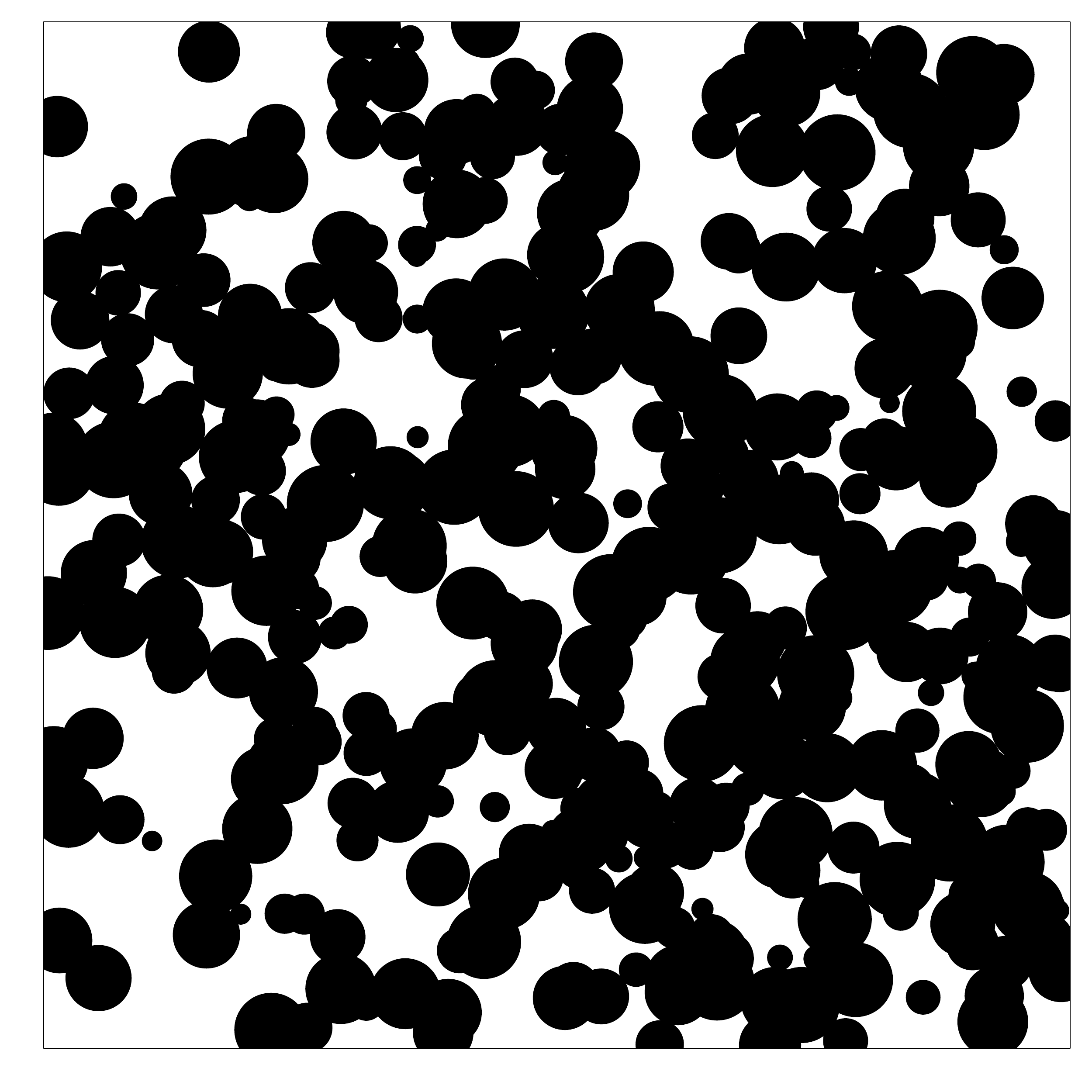}
  \end{tabular}
  }
      \caption{{\small Samples of: the area-interaction process with  $z=0.1$, $\t_1=0.2$ (left); the perimeter-interaction process with $z=0.2$, $\t_2=0.4$ (middle); the Euler-Poincar\'e-interaction process with $z=0.1$, $\t_3=1$ (right). The window is $[0,50]^2$ and $\mu$ is the uniform law on  $[0.5,2]$.}}\label{samples}
 \end{figure}

Some simulations are shown in Figure~\ref{samples}. They correspond to Quermass-interaction models involving only one non-null interaction parameter: the area-interaction process ($\t_2=\t_3=0$), the perimeter-interaction process ($\t_1=\t_3=0$), and the Euler-Poincar\'e-interaction process ($\t_1=\t_2=0$). These particular situations show the rich variety of random sets that the Quermass-interaction model can provide. Note that in the three situations  displayed in Figure~\ref{samples}, the interaction parameter is positive, so that the resulting random set is more likely to induce a lower value of Minkowski functional (resp. $\A$, $\LL$ or $\x$) than in the Boolean case. These simulations have been done by a   Metropolis-Hasting algorithm as presented in \cite{MH1}.

\begin{remarque}\label{rem_reference_measure}
Let us mention some possible generalisations of Definition~\ref{Quermassfini}. 

Consider first the case of a continuous distribution for the radii. To avoid an arbitrary choice of $\mu$, it is possible to simply choose $\mu$ as the uniform distribution on $[0,R_0]$ and then consider an extra term in the density \eqref{dens_quermass} for the law of the radii. For example, assume we would like to fit a  beta distribution to the radii, then $g_\L$ would become proportional to $$ z^{n} \mathrm{e}^{-H^{\theta}(\oo_\L)}\frac{R_0^n}{B(\nu_1,\nu_2)^n}\prod_{i=1}^{n}\left(1-\frac{R_i}{R_0}\right)^{\nu_1-1}\left(\frac{R_i}{R_0}\right)^{\nu_2-1}$$
where $\oo_\L = \{(x_1,R_1),\dots,(x_n,R_n)\}$, $\nu_1>0$, $\nu_2>0$ and $B(\nu_1,\nu_2)$ denotes the beta function. Then the parameters $\nu_1,\nu_2$ can be  estimated along with $\theta$ thanks to the Takacs-Fiksel  procedure described in Section~\ref{sec:estimation}, even if the germs are not observed. For simplicity of presentation, we do not include this generalisation in the following.

 A second generalisation concerns the multi-scale setting (see \cite{Ambler,Goulard,Picard}), where $\mu$ is a discrete measure, leading to different types of balls.  The interaction may then differ according to the type of balls.  For example, in presence of two types of balls with respective radius $R_1$ and $R_2$,  we choose $\mu=\delta_{R_1}+\delta_{R_2}$.  Then the density $g_\L$ for an area interaction becomes proportional to
$$z_1^{n_1}z_2^{n_2} \mathrm{e}^{-\theta_1\A(\U_{\oo_1})-\theta_2\A(\U_{\oo_2})+\theta_{12}\A(\U_{\oo_1}\cap\,\U_{\oo_2})}$$
where $\oo_1$ (resp. $\oo_2$) denotes the set of grains with radius $R_1$ (resp. $R_2$) in $\L$, $n_1$ (resp. $n_2$) denotes their number, and $\U_{\oo_1}$ (resp. $\U_{\oo_2}$) denotes their union. Here again, the parameters $z_1$, $z_2$, $\t_1$, $\t_2$ and $\t_{12}$ can be estimated with the procedure described in Section~\ref{sec:estimation}, provided the two sets $\U_{\oo_1}$ and $\U_{\oo_2}$ are observed. 
\end{remarque}

\subsection{Quermass model on the whole plane: The Markov property}\label{sec:gibbs}

In this section the Markov property of the Quermass model is displayed via the Georgii-Nguyen-Zessin (GNZ) equation (\cite{Georgii, NZ}). 
An alternative presentation, from a statistical physics point of view, can be found in  \cite{Dereudre09}.

First, let us define the local energy of a marked point $(x,R)$ with respect to a configuration $\oo$ by the following expression
\begin{equation}\label{localenergy}
h^{\theta}((x,R),\oo):=H^{\theta}(\oo_{B_x} \cup(x,R))-H^{\theta}(\oo_{B_x}),
\end{equation}
where ${B_x}$ denotes the ball $B(x,2R_0)$. The local energy is related to the Papangelou intensity $$\lambda^*((x,R),\oo):= \frac{g_{B_x}(\oo_{B_x} \cup(x,R); z, \theta)}{ g_{B_x}(\oo_{B_x}; z, \theta)},$$ 
by $\lambda^*((x,R),\oo)=z\,\mathrm{exp}(-h^{\theta}((x,R),\oo))$. Note that the above relation makes sense from Definition~\ref{Quermassfini} since the set ${B_x}$ is bounded.

We have the following characterisation of the Quermass process via the GNZ equation.

\begin{proposition}[Georgii \cite{Georgii}, Nguyen-Zessin \cite{NZ}] \label{NZprop}
For any bounded set $\L\subset\Rd$, a probability measure $P$ on $\O_\L$ is the Quermass point process on $\L$ with parameter $\theta\in\mathbb{R}^3$,  intensity $z>0$ and  distribution $\mu$ of the radii (i.e. $P=P^{z,\theta}_\L$) if and only if for any non-negative function $f$ from $\E\times\O_\L$ to $\R$ 

\begin{equation}\label{NZeq}
\mathbb E\left( \sum_{x\in\boldsymbol{\omega}_\L} f\left((x,R),\boldsymbol{\omega}_\L\backslash (x,R)\right) \right) = \mathbb E\left( \int_0^{R_0} \int_\L  z\;\mathrm{e}^{-h^{\theta}((x,R),\boldsymbol{\omega}_\L)}f\left((x,R),\boldsymbol{\omega}_\L\right) \mathrm d x\;\mu( \mathrm dR)\right),
\end{equation}
where $\mathbb E$ denotes the expectation with respect to $P$ and $\boldsymbol{\omega}$ follows the distribution $P$.
\end{proposition}

The GNZ equation involves the expectation under the Quermass process  of two completely different types of expressions. This  equation is the starting point of the Takacs-Fiksel estimation procedure to be presented in Section~\ref{TFsection}.

In the present paper, we mainly consider  Quermass processes on bounded windows as presented in Definition~\ref{Quermassfini}. However,  it is necessary to consider Quermass processes on the full plane $\R^2$ for questions involving asymptotic properties of estimators (such as consistency or asymptotic normality). This extension is not trivial because for an infinite  configuration $\oo$ in $\Omega$,  $\A(\U_{\oo})$ and  $\LL(\U_{\oo})$ are  infinite while $\x(\U_{\oo})$ can be infinite, minus infinite or an indeterminate form. So the energy $H^{\theta}(\oo)$ in \eqref{energy} is in general an indeterminate form and the definition  \ref{Quermassfini} makes no sense in this case. Nevertheless, a definition of the Quermass process on $\Rd$ is possible due to the GNZ equation \eqref{NZeq}. 

\begin{definition}\label{deffullplane}
A probability measure $P^{z,\theta}$ on $\O$ is a  Quermass point process on $\Rd$ with parameter $\theta\in\mathbb{R}^3$,  intensity $z>0$ and  distribution $\mu$ of the radii if for any non-negative function $f$ from $\E\times\O$ to $\R$, for any bounded set $\L\subset\Rd$, 
\begin{equation*}
\mathbb E\left( \sum_{x\in\boldsymbol{\omega}_\L} f((x,R),\boldsymbol{\omega}\backslash (x,R))\right)= \mathbb E\left( \int_0^{R_0} \int_\L  ze^{-h^{\theta}((x,R),\boldsymbol{\omega})}f((x,R),\boldsymbol{\omega}) \mathrm d x\;\mu( \mathrm dR)\right),
\end{equation*}
where $\mathbb E$ denotes the expectation with respect to $P^{z,\theta}$ and $\boldsymbol{\omega}$ follows the distribution $P^{z,\theta}$.
\end{definition} 

The existence of a measure $P^{z,\theta}$ (for any $z>0$ and any $\theta\in\R^3$) satisfying Definition~\ref{deffullplane} was proved recently in \cite{Dereudre09}.


\section{The estimation procedures}\label{sec:estimation}
Let us consider  a realisation $\oo$  of a Quermass point process defined on a set $\Lambda'$, where $\Lambda'\subseteq\R^2$. We assume that we observe the random set $\U_{\oo}\cap\Lambda$, where $\Lambda$ is the observation window, i.e. $\Lambda$ is a bounded set and  $\Lambda\subseteq\L'$. In practice two typical cases occur: $\Lambda'=\Lambda$, i.e. the observation window coincides with the domain of definition of $\oo$  (which is therefore bounded in this case),  or  $\Lambda'=\R^2$ and $\Lambda\subset\L'$. We denote by $z^*$ and $\theta^*$ the true parameters defining the law of the observed Quermass process in Definitions \ref{Quermassfini} and \ref{deffullplane}. Since $\oo_\L$ is not observed, but only  $\U_{\oo}\cap\Lambda$,  the positions of the marked points $(x,R)$ in $\oo_\L$, i.e.\ the germs,  are unknown. A challenging task is  to estimate the parameters $z$ and  $\theta=(\theta_1,\theta_2,\theta_3)$ while the germs are not observed, which is not a typical problem of inference for  Gibbs point processes. The first subsection explains the limitation in using the maximum likelihood procedure (MLE). The second subsection presents the main approach of this paper which is the Takacs-Fiksel (TF) method. 


\subsection{The maximum likelihood approach}\label{MLE}

When $z^*$ is known, the  classical maximum likelihood approach to  estimation of the parameter $\theta$ is possible. This  has been investigated in \cite{MH2}, where  the authors introduce an original procedure based only on the connected components that are completely included in the observation window. The edge effects are then reduced. However, the maximum likelihood procedure does not allow us to estimate the intensity parameter $z$ since the number of points is not observed. If $z^*$ is unknown, a two step procedure is proposed in   \cite{MH2}:  first  estimating $z$ as if the observable data come from a Boolean model (various methods are available in this setting),  then,  in a second step,  applying the MLE procedure to estimate $\theta$.

Unfortunately,  it seems that this two step procedure induces a strong restriction on  the possible values for the  interaction parameter $\theta$. Indeed, in the first step, the estimation $\hat z$ is chosen  such that the Boolean model with intensity $\hat z$ has the same specific volume as the observed real data. Similarly, in the second step,  the MLE estimation of $\theta$ ensures that the corresponding Quermass process has the same specific volume as the real data.
So, this two step procedure concerns Quermass processes with intensity parameter $\hat z$, that exhibit the same specific volume as the Boolean model with intensity $\hat z$. This turns out to be a strong restriction.  For example, it is well-known that the area-interaction process (i.e. $\theta_2=\theta_3=0$) with parameters $z$ and $\theta_1>0$ is strictly stochastically dominated by the Boolean model with intensity $z$ \cite{GK}, which   implies a strict comparison of their specific volumes. For this example, the lonely possible value of $\t_1$ following the previous restriction is thus $\t_1=0$. Therefore the two step procedure of \cite{MH2}  may imply a strong bias of the estimates, as confirmed in Figure~\ref{fig-comparison}.

\subsection{The Takacs-Fiksel approach}\label{TFsection}

The possibility to use the Takacs-Fiksel procedure for estimating all parameters of  the Quermass process, including $z$, was   recently emphasised in \cite{CDDL}. In this section we recall the procedure, which depends on the choice of some test functions. Next we propose various test functions, which allow to circumvent the unobservability issue of the Quermass random set.

\subsubsection{The general procedure}

Consider, for any non-negative function $f$ from $\E\times\O$ to $\R$, and any $z>0$, $\theta\in\mathbb R^3$, the random variable 
\begin{equation}\label{Delta}
 C_{\L}^{z,\theta}(\boldsymbol{\omega}; f)= \sum_{x\in\boldsymbol{\omega}_{\L}} f((x,R),\boldsymbol{\omega}\backslash (x,R))-  \int_0^{R_0} \int_{\L}  z\;\mathrm e^{-h^{\theta}((x,R),\boldsymbol{\omega})}f((x,R),\boldsymbol{\omega}) \mathrm d x\;\mu( \mathrm dR),
\end{equation}
where $\boldsymbol{\omega}$ follows the distribution $P^{z^*,\theta^*}$.
From the ergodic theorem, if $\Lambda$ is sufficiently large, $C_{\L}^{z,\theta}(\boldsymbol{\omega}; f)/|\L|$ is approximatively equal to   $\mathbb E_{z^*,\t^*}(C_{[0,1]^2}^{z,\theta}\left(\boldsymbol{\omega}; f)\right)$,  where $\mathbb E_{z^*,\t^*}$ is the expectation with respect to $P^{z^*,\theta^*}$. 
Moreover, from the GNZ equation (\ref{NZeq}) it is easy to show that $\mathbb E_{z^*,\t^*}(C_{[0,1]^2}^{z^*,\theta^*}\left(\boldsymbol{\omega}; f)\right)=0$. Therefore, for any function $f$, the random variable $C_{\L}^{z,\theta}(\boldsymbol{\omega}; f)$  should be close to zero when $z=z^*$ and $\t=\t^*$. This remark is the basis of the Takacs-Fiksel approach. 

Given $K$ functions $(f_k)_{1\le k \le K}$, the Takacs-Fiksel estimator is simply defined by

\begin{equation}\label{Takacs}
(\hat z,\hat \theta):= \mathop{\arg\min}_{z,\theta}  \sum_{k=1}^K \left(C_{\L}^{z,\theta}(\oo; f_k)\right)^2,
\end{equation}
where $\oo$ is a realisation of $\boldsymbol{\omega}$.
The strong consistency and  asymptotic  normality of $(\hat z,\hat \theta)$  when the observation window $\Lambda$ grows to $\Rd$ are discussed in \cite{CDDL}. The main ingredient for consistency is to ensure that $(z^*,\t^*)$ is the unique solution to the optimisation problem in \eqref{Takacs} when $\L$ is large enough. This identifiability is not easy to check. Some criteria are given in \cite{CDDL}. If $p$ is the number of parameters to estimate, it turns out that choosing $K=p$ test functions may be not sufficient to  ensure identifiability.  If $K>p$, then identifiability is in general achieved.

In the setting of Quermass model, it is in general not possible to compute $C_{\L}^{z,\theta}(\oo; f)$ given the observation $\U_{\oo}\cap\L$. Consider for instance the pseudo-likelihood estimator. It is a particular case of the Takacs-Fiksel procedure where  $f_k=\partial h^{\theta}((x,R),\oo)/\partial \theta_k$, $k=1,2,3$ and $f_4=1/z$ (see section 3.2 in \cite{CDDL}). In this case the computation of \eqref{Delta} requires the observation of all germs in $\oo_\L$ and so the pseudo-likelihood procedure is not feasible. However, it is possible to find some test functions $f$  such that  
\begin{enumerate}
\item for all $x\in\L$ and $R\in[0,R_0]$, $f((x,R),\oo)$ is computable, making  the integral term in \eqref{Delta}  calculable, 
\item the sum term in \eqref{Delta} is computable from the observation of $\U_{\oo}\cap\L$ only.
\end{enumerate}
Some examples are provided below.

Note finally that for any $x\in\L$ and any $R\in[0,R_0]$ the computation of the local energy $h^{\theta}((x,R),\oo)$ is always possible (up to some possibly edge effects that will be discussed in Section~\ref{edge}). Indeed, denoting $\tilde\U(x,R)=\U_{\oo}\cap B(x,R)$,  the additivity of $\mathcal A$, $\mathcal L$ and $\chi$ implies that
\begin{equation}\label{hlocalized} 
h^{\theta}((x,R),\oo)=\t_1\left[\pi R^2-\A(\tilde\U(x,R))\right] + \t_2\left[2\pi R-\LL(\tilde\U(x,R))\right] + \t_3\left[1-\x(\tilde\U(x,R))\right].
\end{equation}

\subsubsection{Some well-adapted test functions}\label{functions}

In this section, we introduce some test functions $f$ that make  $C_{\L}^{z,\theta}(\oo; f)$ in \eqref{Delta} computable from the observation of $\U_{\oo}\cap\L$ only.  This provides a  solution for estimation in  the Quermass model when the germs are not observed. We restrict our presentation here to the case $\oo=\oo_\L$, or equivalently $\L'=\L$,  meaning that the observation window coincides with the domain of definition of $P^{z^*,\t^*}$. The general case involves some edge effects which are discussed in Section~\ref{edge}.

Let us first consider the test function $f_0$ defined by 

\begin{equation}\label{f0}
f_0((x,R), \oo)=\text{Length} \Big(\mathcal S(x,R)\cap (\U_{\oo})^c\Big),
\end{equation}
where $\mathcal S(x,R)$ is the sphere with centre $x$ and radius $R$. The quantity $f_0((x,R), \oo)$ is actually the length of arcs from the sphere $\mathcal S(x,R)$ which are outside $\U_{\oo}$. Although the quantity  $f_0((x,R), \oo\backslash (x,R))$ is not observable for any $(x,R)$ in $\oo_\L$ if $R\not=0$, its sum over $(x,R)\in\oo_\L$ is nothing else than the total perimeter of $\U_{\oo_\L}$:
\begin{equation}\label{f0sum}\sum_{x\in\oo_{\L}} f_0((x,R),\oo\backslash (x,R))  =\LL(\U_{\oo_\L}).\end{equation} 
Moreover, it is possible to compute $f_0((x,R), \oo)$ for any $x\in\L$ and $R\in[0,R_0]$, making the integral in \eqref{Delta} calculable. It follows that $C_{\L}^{z,\theta}(\oo; f_0)$ is computable.

Similarly for any $\alpha>0$, we consider the test function 
\begin{equation}\label{falpha}
f_\alpha((x,R), \oo)=\text{Length} \Big(\mathcal S(x,R+\alpha)\cap (\U_{\oo}\oplus B(0,\alpha) )^c\Big),
\end{equation}
where $\U_{\oo}\oplus  B(0,\alpha)$ denotes the $\alpha$-parallel set of $\U_{\oo}$ defined by the Minkowski sum of $\U_{\oo}$ and $ B(0,\alpha)$. In this case \eqref{Delta} is also computable and the sum therein   reduces to the perimeter of the $\alpha$-parallel set: 
\begin{equation}\label{falphasum}
\sum_{x\in\oo_\L} f_\alpha((x,R), \oo\backslash (x,r))=\LL\left(\U_{\oo_\L}\oplus B(0,\alpha)\right).\end{equation}

\begin{figure}[htbp]
    \setlength{\tabcolsep}{0.1cm} \centerline{
  \begin{tabular}[]{ccc}
    \includegraphics[angle=0,scale=.2]{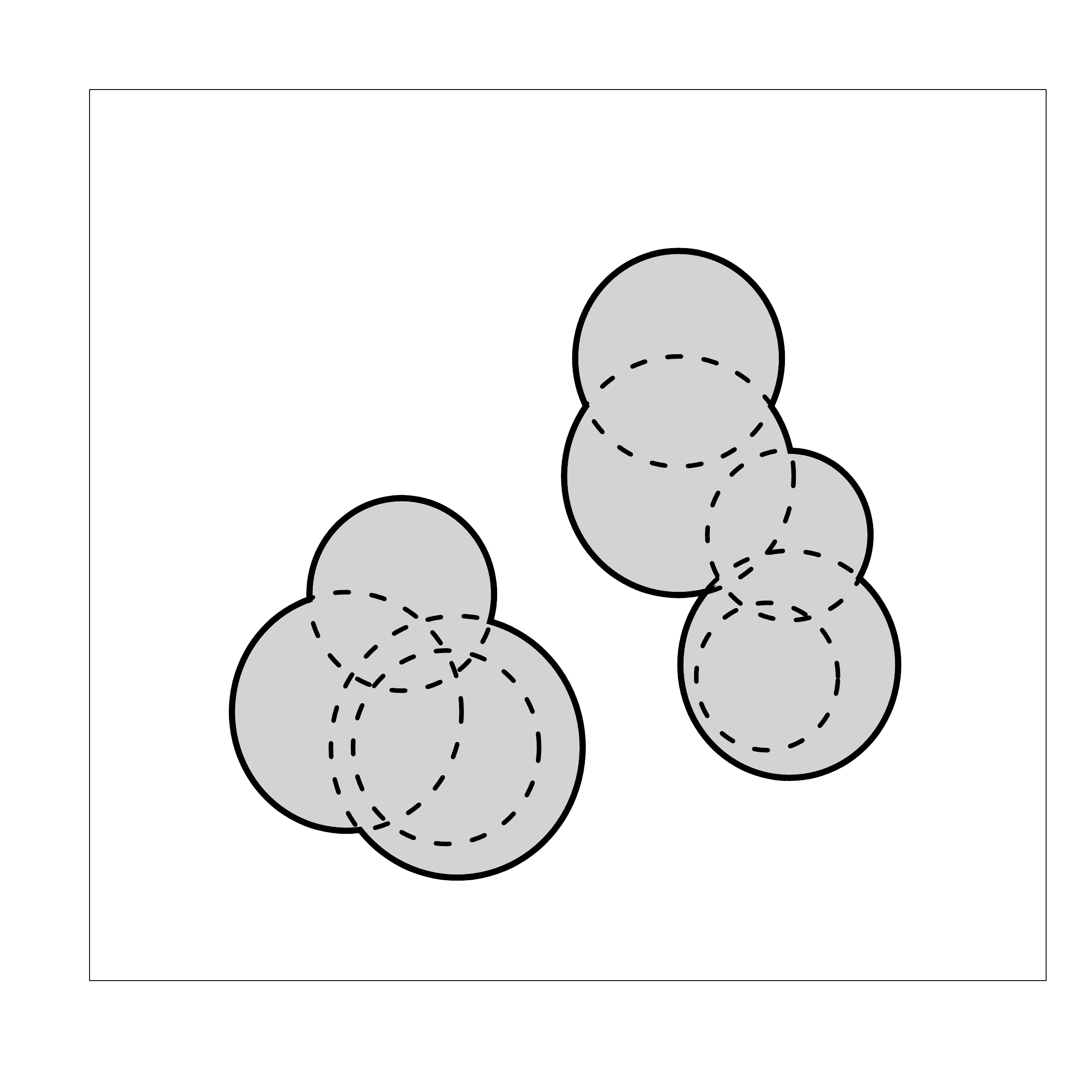}&
  \includegraphics[angle=0,scale=.2]{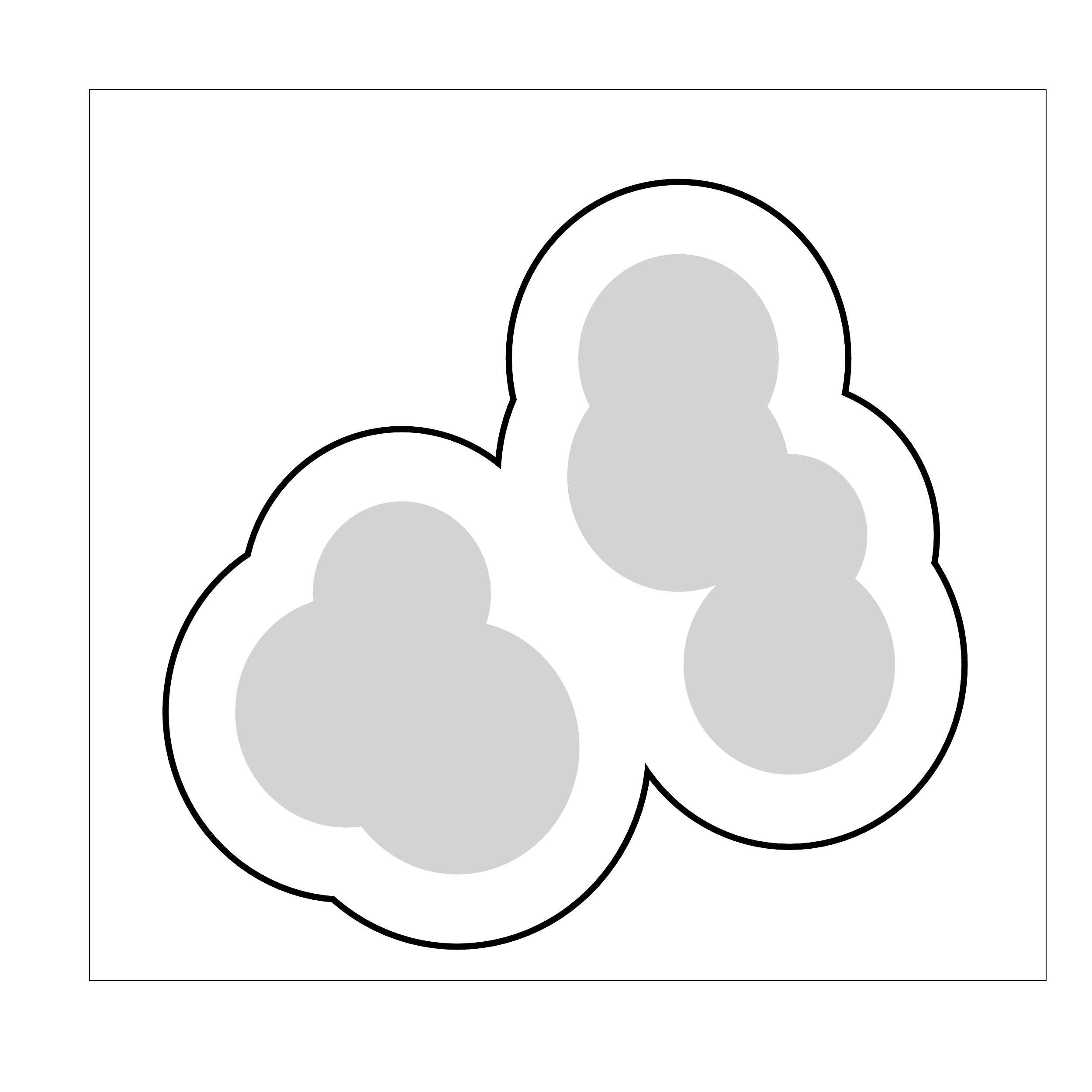}&
   \includegraphics[angle=0,scale=.2]{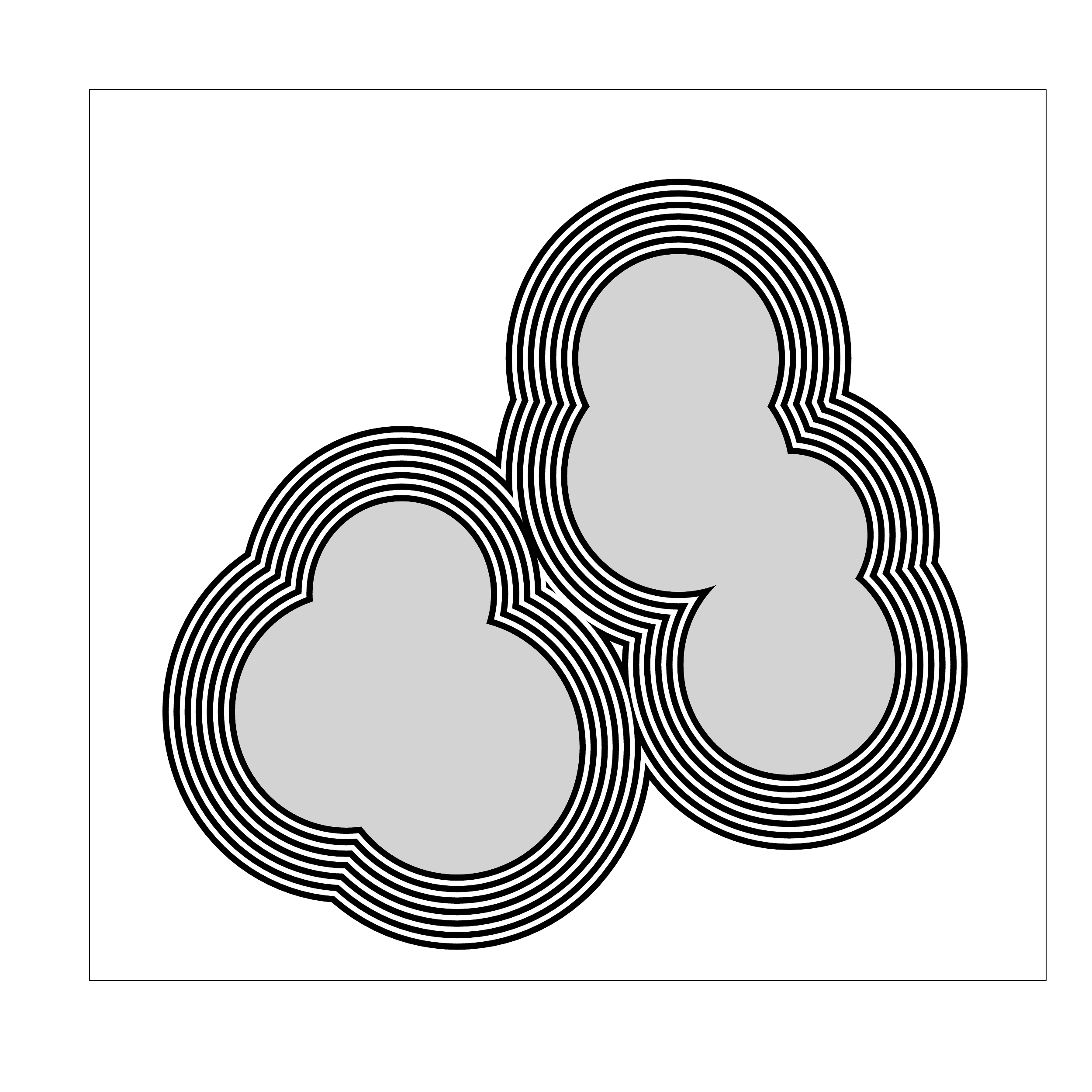}
  \end{tabular}
  }
      \caption{{\small Left: $\U_{\oo_\L}$ in grey, as the union of the balls bordered with dashed lines; the length of the solid black line is the quantity \eqref{f0sum}. Middle: the same $\U_{\oo_\L}$ in grey;  the length of the solid black line is the quantity \eqref{falphasum} for some $\alpha>0$. Right: the same $\U_{\oo_\L}$ in grey; the sum of lengths of the solid black lines is the quantity \eqref{fsumsum} for 7 $\alpha_i$'s.}}\label{fig:test-functions}
 \end{figure}

In the same spirit as  \eqref{f0sum}, it is natural to look for some test functions   $f_{\mathrm{area}}$ and $f_{\mathrm{ep}}$ such that the sum in \eqref{Delta}  reduces to $\A(\U_{\oo_\L})$ and $\x(\U_{\oo_\L})$ respectively, and  the integral term in \eqref{Delta} is calculable. However, we have not found any.
Nevertheless, considering the  test function
\begin{equation}\label{fsum} f_{\text{sum}}=\sum_i  f_{\alpha_i},\end{equation}
where the sum is over some suitable finite set of non negative $\alpha_i$'s, we obtain 
\begin{equation}\label{fsumsum}
\sum_{x\in\oo_\L} f_{\text{sum}}((x,R), \oo\backslash (x,r))=\sum_i \LL\left(\U_{\oo_\L}\oplus B(0,\alpha_i)\right),\end{equation}
which is  closely related to the area of the complementary set of $\U_{\oo_\L}$ in $\L$, in agreement with the so-called Cavalieri estimator in stereology, see Figure~\ref{fig:test-functions}. For this reason, $f_{\text{sum}}$ can be viewed as a substitute for $f_{\mathrm{area}}$.

Let us finally introduce the test function $f_{\text{iso}}$ which indicates,  for any $(x,R)\in\oo_\L$,  whether  $B(x,R)$ is an isolated ball in $\U_{\oo}$:
\begin{equation}\label{fiso}
f_{\text{iso}}((x,R), \oo)=\left\{ \begin{array}{ll}
1 & \text{ if } \mathcal{S}(x,R)\cap \U_{\oo}=\emptyset,\\
0 & \text{ otherwise. }
\end{array} \right.
\end{equation}

In this case, $\sum_{x\in\oo_\L} f_{\text{iso}}((x,R), \oo\backslash (x,R))$ is  the number of isolated balls in $\U_{\oo_\L}$. Note that an isolated ball can contain smaller balls completely included inside it.

\begin{remarque}
The test functions  $f_0$, $f_\alpha$, $f_{\text{sum}}$, $f_{\text{iso}}$ introduced above satisfy the regularity conditions assumed in \cite{CDDL}, which   imply   consistency and asymptotic normality of the associated TF estimators, provided identifiability holds (see Section~\ref{localmin}).

\end{remarque}

 \section{Practical aspects of Takacs-Fiksel estimation}\label{PA}

\subsection{Computation of the contrast function}\label{sec:computation}

The TF estimator in \eqref{Takacs} requires to compute $C_{\L}^{z,\theta}(\oo; f)$ defined in \eqref{Delta} for the different well-adapted  test functions $f$ introduced in Section~\ref{functions}. Note that  in Section~\ref{edge}, we will explain how to handle edge effects in the computation of \eqref{Delta}, so we do not address this question in this section. 

The integral term in  \eqref{Delta} can  be approximated by a Monte Carlo approach, based on $N$ independent points $x_1,\dots,x_N$ uniformly distributed on $\Lambda$ and $N$ independent realisations $R_1,\dots, R_N$ from $\mu$. This leads to the approximation 
\begin{equation}\label{Deltaapprox}
 C_{\L}^{z,\theta}(\oo; f)\approx \sum_{x\in\oo_{\L}} f((x,R),\oo \backslash (x,R))- \frac 1N \sum_{i=1}^N z\;\mathrm e^{-h^{\theta}((x_i,R_i),\oo)}f((x_i,R_i),\oo).
\end{equation}

Let us consider the computation of \eqref{Deltaapprox} for the test functions $f_0$, $f_\alpha$, $f_{\text{sum}}$ and some fixed values of $z$ and $\theta$. The first sum above reduces respectively to the perimeter of $\U_{\oo_\L}$, the perimeter of the   $\alpha$-parallel set of $\U_{\oo_\L}$, or the sum of such perimeters.   The construction of $\alpha$-parallel sets of $\U_{\oo_\L}$ and the calculation of  the above-mentioned perimeters may be achieved using some mathematical morphology tools, see e.g.  \cite{Ohser, Mrkvicka, SchmidtSpodarev}. The second sum in \eqref{Deltaapprox} involves the local energy of each new ball $(x_i,R_i)$, given by \eqref{hlocalized}, and $f((x_i,R_i),\oo)$, for $f=f_0, f_\alpha, f_{\text{sum}}$. While for each $i$ we can use again  some mathematical morphology tools to calculate these two terms, the computation for the whole sum, which requires to repeat $N$ times these calculations,  may become difficult to implement.

An alternative procedure to compute \eqref{Deltaapprox} is to  approximate  the observed set $\U_{\oo}\cap\L$ by any union of balls  (e.g. by the procedure described in \cite{Thiedmann}) and to consider the associated power tessellation (see \cite{MH1}). Then the construction of approximated $\alpha$-parallel sets of $\U_{\oo_\L}$ becomes straightforward since it suffices to increase the radius of each ball by $\alpha$. Moreover, the calculation of the above-mentioned perimeters,  the calculation of the local energy of any new ball $(x_i,R_i)$, and the calculation of  $f((x_i,R_i),\oo)$, for $f=f_0, f_\alpha, f_{\text{sum}}$,  are easily deduced from the associated power tessellation as described in \cite{MH1}. 

We emphasise that the above procedure does not depend on the union of balls used to approximate   $\U_{\oo}\cap\L$. In order to minimise the computational complexity, an approximation involving a low number of balls is therefore preferable.  Moreover, this union of balls is only used for computational reasons and is not related to  the germ-grain representation of  $\U_{\oo}$. In particular, the number of balls in this approximation is not relevant to estimate $z$.

Finally, let us turn to the computation of \eqref{Deltaapprox} for the test function $f_{\text{iso}}$.
Contrary to  $f_0$, $f_\alpha$ and $f_{\text{sum}}$, the test function $f_{\text{iso}}$ is strongly related to the ball structure of the germ-grain model, since the first sum in \eqref{Deltaapprox} corresponds to the number of isolated balls in $\U_{\oo_\L}$.  In this sense $f_{\text{iso}}$ seems less natural  for applications. However, as it will be demonstrated in Section~\ref{sec:simu}, the test function $f_{\text{iso}}$ appears to provide a relevant information for the estimation of the parameters. For this reason, it can be important to include it  in \eqref{Takacs}.  In practice, we have then to decide what can be considered as an isolated ball. A solution can be to view an isolated component of $\U_{\oo_\L}$ as an isolated ball, if its diameter is smaller than some constant chosen a priori. Another solution is to use  the approximation of $\U_{\oo}\cap\L$  by a union of balls, as mentioned above: An isolated component can then be considered as an isolated ball if it is approximated by an isolated ball. 

In this work, the practical implementation of the TF procedure has been conducted using an approximation of $\U_{\oo}\cap\L$ by a union of balls and the construction of the associated power tessellation.

 \subsection{Minimisation of the contrast function}\label{implementation}
 As described in Section~\ref{functions}, many computable test functions are available to implement the TF procedure for the Quermass model: $f_0$, $f_{\alpha}$  for any $\alpha>0$, $f_{\text{sum}}$, $f_{\text{iso}}$.

In order to estimate $p$ unknown parameters among $z$, $\t_1$, $\t_2$ and $\t_3$, it suffices to choose  $K\geq p$ test functions as above and to solve  \eqref{Takacs} where $\t=(\t_1,\t_2,\t_3)$. 

The TF optimisation \eqref{Takacs} in $z$ can be done  explicitly. We deduce that the solution $(\hat z,\hat\t)$ of  \eqref{Takacs} necessarily belongs to the implicit manifold $z=\tilde z(\theta)$ where 
\begin{equation}\label{ztheta}
 \tilde z(\theta): = \frac{ \sum_{k=1}^K  S_k I_k(\theta)}{   \sum_{k=1}^K  I_k^2(\theta)},
\end{equation} 
with
$$S_k=\sum_{x\in\oo_\L} f_k((x,R),\oo\backslash (x,R))$$ and $$I_k(\t)= \int_0^{R_0} \int_{\L}  e^{-h^{\theta}((x,R),\oo)}f_k((x,R),\oo) \mathrm d x\;\mu( \mathrm dR).$$
Therefore, in practice, $\theta$ is first estimated by the solution $\hat\t$ of  \eqref{Takacs} where $z$ is replaced by $\tilde z(\t)$, i.e.

\begin{equation*}
\hat \theta = \mathop{\arg\min}_{\theta}  \sum_{k=1}^K \big(S_k - \tilde z(\t) I_k(\t)\big)^2. 
\end{equation*}
This solution may be obtained by a grid search optimisation procedure. Then $z$ is estimated by $\hat z=\tilde z(\hat\theta)$. 

Recall that the practical computation of all terms involved in \eqref{ztheta} and \eqref{Takacs} can be conducted as explained in Section~\ref{sec:computation}.

This procedure allows to consider several estimators, depending on the number and the choice of test functions used in \eqref{Takacs}. It is not easy to find an optimal choice. The asymptotic variance of TF estimators is known (cf \cite{CDDL}),  but appears intractable to be optimised. Some simulations are thus mandatory to provide some recommendations, see Section~\ref{sec:simu}.

\subsection{Edge effects}\label{edge}

Based on the observation of $\U_\oo\cap\Lambda$, two types of edge effects may occur in the computation of \eqref{Delta} or \eqref{Deltaapprox}. 

First, edge effects may occur in the integral term in \eqref{Delta}, or the second sum in \eqref{Deltaapprox},  for the computation of $h^\theta((x,R),\oo)$ and $f((x,R),\oo)$, for $f=f_0,f_\alpha,f_{\text{sum}},$ or $f_{\text{iso}}$, when $x$ is close to the boundary of $\Lambda$. In view of \eqref{hlocalized}, \eqref{f0}, \eqref{falpha}, \eqref{fsum}, \eqref{fiso}, this type of edge effects does not occur for the marked points $(x,R)$ such that $B(x,R_{\text{max}})\subset\Lambda$, where $R_{\text{max}}=R_0+\max_i \alpha_i$. Therefore, in practice, this first type of edge effects can be avoided by considering minus sampling, i.e. the estimation on $\Lambda^-$, where $\Lambda^-$ denotes the eroded set of $\Lambda$ by $R_{\text{max}}$. The estimators are then defined by \eqref{Takacs} where $\L$ is replaced by $\L^-$. 

Second, edge effects can occur for the first sum term in \eqref{Delta} and \eqref{Deltaapprox}. Let us consider $f=f_0$. If $\oo=\oo_\Lambda$, as it was assumed for simplicity in Section~\ref{functions}, then from \eqref{f0sum} this sum term reduces to $\LL(\U_{\oo_\L})$. As we only observe $\U_\oo\cap\Lambda$, an approximation in this case could be $\LL(\U_{\oo}\cap\Lambda)$ and we have  $\LL(\U_{\oo}\cap\Lambda)\leq \LL(\U_{\oo_\L})$, where the equality occurs if and only if  $\U_\oo\subset\Lambda$. In the general case when $\oo\not=\oo_\L$,  meaning that the domain of definition of $\oo$ strictly contains $\L$, the equality \eqref{f0sum} does not hold any more and  some extra terms due to edge effects have to be added in the right hand side of \eqref{f0sum}. Furthermore, to avoid the first type of edge effects, the estimation procedure is implemented on the eroded domain $\Lambda^-$. In this general case we may also use  the approximation 
\begin{equation}\label{approxedge1}\sum_{x\in\oo_{\L^-}} f_0((x,R),\oo\backslash (x,R))  \approx \LL(\U_{\oo}\cap\Lambda^-),\end{equation} where the last perimeter can be computed by  some mathematical morphology tools. Alternatively, $\U_{\oo}\cap\Lambda$ can be approximated by a union of balls, as suggested in Section~\ref{sec:computation}, namely for some integer $n$, $\U_{\oo}\cap\Lambda\approx \bigcup_{i=1}^n B(y_i,r_i)$ where $y_i\in\L$ and $r_i>0$, and we can use the approximation 
\begin{equation}\label{approxedge2}\sum_{x\in\oo_{\L^-}} f_0((x,R),\oo\backslash (x,R))  \approx \sum_{y_i\in\L^-} \text{Length} \bigg(\mathcal S(y_i,r_i)\cap \Big( \bigcup_{\begin{subarray}{l} j=1\\ j\not=i \end{subarray}}^n B(y_j,r_j)\Big)^c\bigg).\end{equation}
This last sum can easily be implemented from the power tessellation based on  $\bigcup_{i=1}^n B(y_i,r_i)$.

Note that error of  approximations in \eqref{approxedge1} and \eqref{approxedge2} involve only edge effects that become negligible when $\Lambda$ tends to $\R^2$. Finally similar approximations can be used for $f=f_\alpha,f_{\text{sum}}$ or $f_{\text{iso}}$ instead of $f=f_0$ and we omit the details. 

In this work, we have decided to handle edge effects by working in \eqref{Takacs} with the eroded domain $\L^-$ instead of $\L$, and by using the approximation \eqref{approxedge2}, which is in agreement with our choice in Section~\ref{sec:computation}.

\subsection{Identifiability}\label{localmin}

The consistency of the TF procedure crucially depends on an identifiability assumption   which basically implies that the contrast function in \eqref{Takacs} has a unique minimum when $\Lambda$ tends to $\R^2$ (see  \cite{CDDL}).
This  assumption is in general satisfied if $K>p$, where $p$ is the number of parameters to estimate. When $K=p$,  identifiability is not easy to check, but has been proved to  hold in some cases (see Example 2 in \cite{CDDL}).   

Nevertheless, these theoretical considerations  only ensure  {\it asymptotic}  identifiability.
In practice, where $\Lambda\not=\R^2$, several local minima of the contrast function in \eqref{Takacs} may arise and the global minimum might lead to an improper estimation. To avoid this problem, a solution is to consider several TF contrast functions, coming from various choices of test functions. They should all share a local minimum in the same region, which allows us to restrict the domain of optimisation.

 \section{Simulations}\label{sec:simu}

 \subsection{Estimation of the intensity parameter $z$}\label{sec:z}
 Let us first assume that $\t_1^*$, $\t_2^*$, $\t_3^*$ are known and let us consider the estimation of $z$, which appears as the main challenge for the Quermass-interaction model when the germs are not observed. In this setting, the equation \eqref{ztheta} with $\theta=\theta^*$ provides an explicit estimator of $z$.  The computation of \eqref{ztheta} is conducted as explained in Section~\ref{PA}:  we use minus sampling and the Monte Carlo approximation \eqref{Deltaapprox} is applied, where we choose $N=2500$.
 
More specifically, we consider the estimators $\hat z_0$, $\hat z_{\alpha_1}$, $\cdots$, $\hat z_{\alpha_{10}}$,  $\hat z_{\text{iso}}$ and $\hat z_{\text{sum}}$,  defined by \eqref{ztheta} where $K=1$ and where  $f_k$ is respectively $f_0$, $f_{\alpha_1}$, $\cdots$, $f_{\alpha_{10}}$,  $f_{\text{iso}}$ and $f_{\text{sum}}$. The latest test function $f_{\text{sum}}$ is defined by the sum  \eqref{fsum} over the ten previous $\alpha_i$, $i=1,\dots,10$.
 
\begin{figure}[htbp]
  \centerline{
  \begin{tabular}[]{ccc}
  \includegraphics[angle=0,scale=.4]{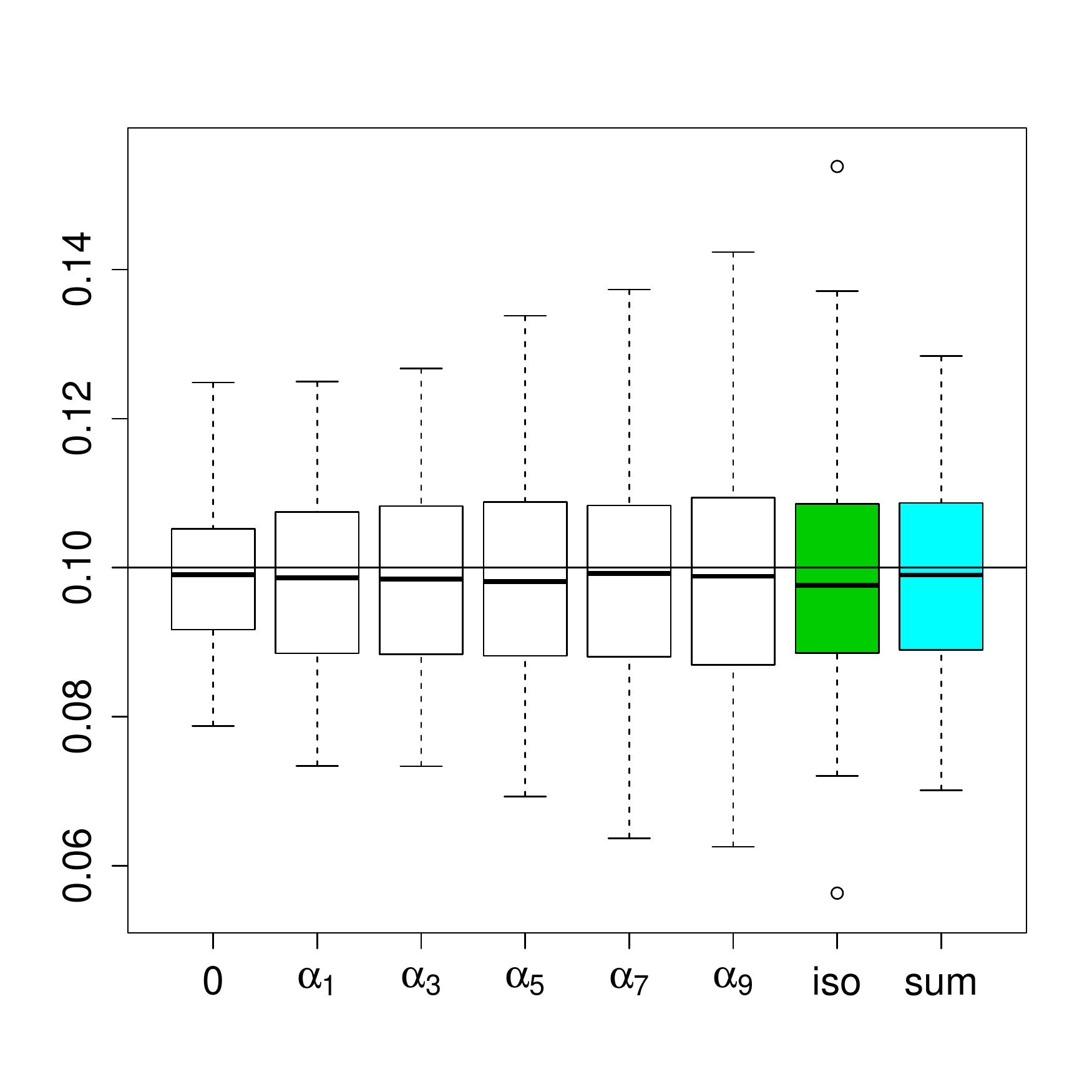}&$\hspace{0.1cm}$&
  \includegraphics[angle=0,scale=.4]{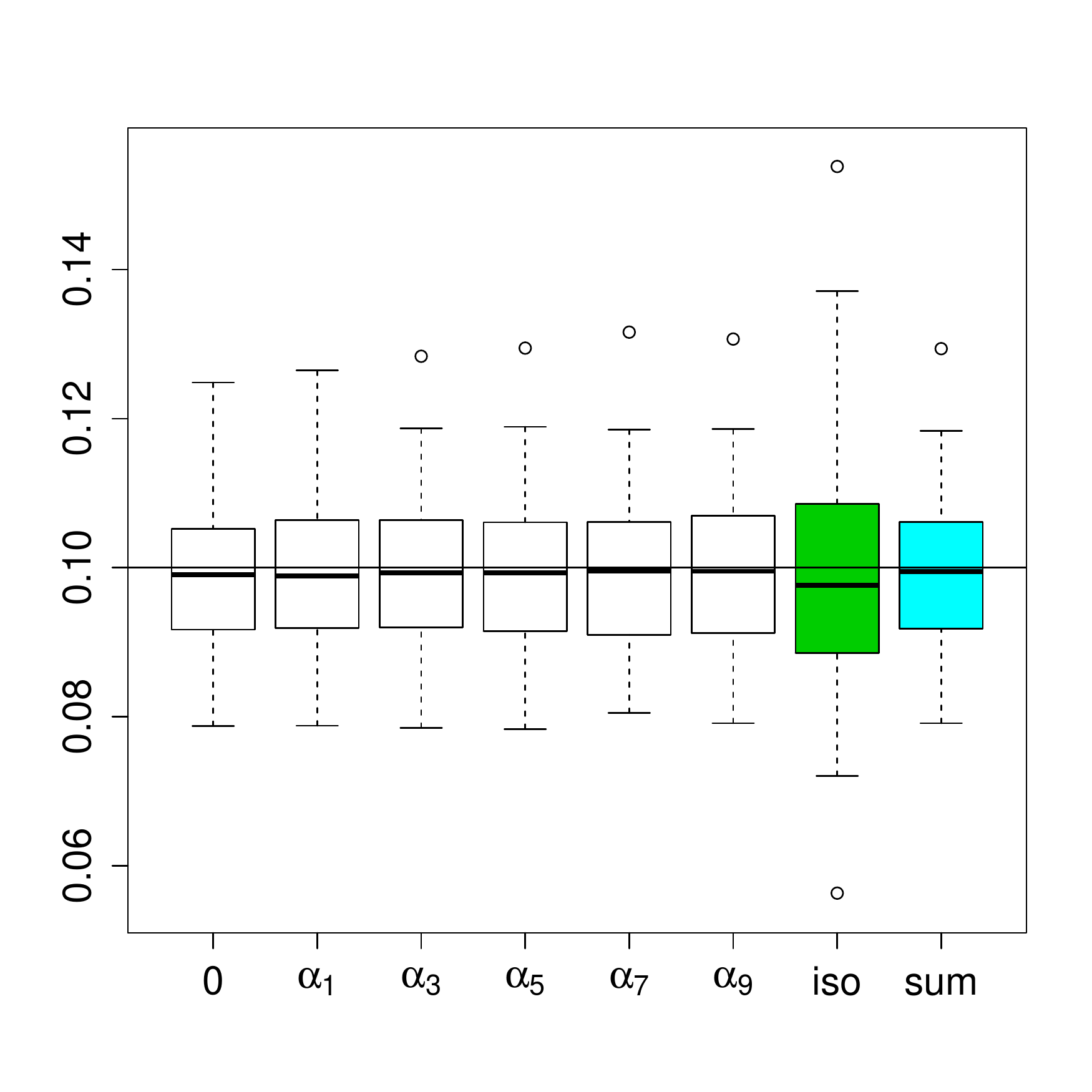}
 \end{tabular}
 }
 \caption{Estimations of $z^*=0.1$ from 100 replicates of an area-interaction process ($\t_1^*=0.2$). Left: Boxplots of  $\hat z_0$, $\hat z_{\alpha_i}$ ($\alpha_i=i/5$), $\hat z_{\text{iso}}$ and $\hat z_{\text{sum}}$. Right: The same boxplots but with $\alpha_i=i/50$.}\label{fig-z}
 \end{figure}

Some results for the estimation of $z$ are displayed in Figure~\ref{fig-z}, based on 100 replicates of the area-interaction process ($\t_2=0$, $\t_3=0$) on $[0,50]^2$, when $z^*=0.1$, $\t_1^*=0.2$ and $\mu$ is the uniform law on $[0.5,2]$. In the left plot, large parallel sets have been used: $\alpha_i=i/5$, $i=1,\dots,10$, while the right plot shows the results with smaller parallel sets: $\alpha_i=i/50$, $i=1,\dots,10$. 

Figure~\ref{fig-z}  confirms that our procedure allows to estimate $z$, even if the germs are not observed. Moreover, it appears that small parallel sets seem to provide better estimates. The same conclusion holds from intensive estimations of $z$ for other  values of $\t_1^*$, $\t_2^*$ and $\t_3^*$ (not presented in this article).

\subsection{Estimation for the $\A$, $\LL$ and $\x$ process}\label{one-interaction}
In this section, we consider the  estimation of the area-interaction process, for short $\A$-process ($\t_2=\t_3=0$), of the perimeter-interaction process, for short $\LL$-process ($\t_1=\t_3=0$) and of the Euler-Poincar\'e-interaction process, for short $\x$-process ($\t_1=\t_2=0$). In these cases, two parameters have to be estimated: the interaction parameter ($\t_1$, $\t_2$ or $\t_3$) and the intensity parameter $z$. 

Many TF estimators are conceivable. For instance, for the area-interaction process, the estimation of $(z,\t_1)$ can be done by the following TF estimators with $K=2$:  ($\hat z_{\text{iso}}, \hat \theta_{\text{iso}}$) based on $(f_{0},f_{\text{iso}})$; ($\hat z_{\alpha}, \hat \theta_{\alpha}$) based on $(f_0,f_{\alpha})$ for some $\alpha>0$; ($\hat z_{\text{sum}}, \hat \theta_{\text{sum}}$) based on $(f_0,f_{\text{sum}})$ where the sum  \eqref{fsum} is done over ten  $\alpha_i$'s.  It is also possible to consider more test functions, as for example with $K=11$: ($\hat z_{\text{all}}, \hat \theta_{\text{all}}$) based on $f_0,f_{\alpha_1},\dots,f_{\alpha_{10}}$ for ten different $\alpha_i$'s. 

Some simulations (omitted here) show that  given an observation $\U_{\oo_\L}$, the  estimations ($\hat z_{\alpha}, \hat \theta_{\alpha}$) may be very variable with $\alpha$, even when using small parallel sets as suggested in Section~\ref{sec:z}. As far as we do not know which estimator is the most efficient, we prefer to use estimators that combine several small parallel sets. For this reason, we will concentrate on ($\hat z_{\text{sum}}, \hat \theta_{\text{sum}}$) and ($\hat z_{\text{all}}, \hat \theta_{\text{all}}$). 

Another natural way to combine information coming from various parallel sets is to consider an aggregated estimator. How to aggregate several estimators is an interesting but difficult question, particularly in presence of strong dependent observations as in the Quermass model.  We simply consider here the median of several estimators, that provides a natural robust combination:  
$$(\hat z_{\text{med}}, \hat \theta_{\text{med}}) = \textrm{median} \big((\hat z_{\alpha_1}, \hat \theta_{\alpha_1}),\dots,(\hat z_{\alpha_{10}}, \hat \theta_{\alpha_{10}})\big).$$

Therefore, the following simulations concern the three estimators: ($\hat z_{\text{sum}}, \hat \theta_{\text{sum}}$), $(\hat z_{\text{med}}, \hat \theta_{\text{med}})$ and ($\hat z_{\text{all}}, \hat \theta_{\text{all}}$), where $\alpha_i=i/50$, $i=1,\dots,10$. Moreover, since the test function $f_{\text{iso}}$ is peculiar as it is strongly related to the ball structure of the germ-grain model, we also assess the performance of ($\hat z_{\text{iso}}, \hat \theta_{\text{iso}}$), based on $(f_{0},f_{\text{iso}})$.

\begin{figure}[htbp]
  \centerline{
  \begin{tabular}[]{ccc}
 \includegraphics[angle=0,scale=0.45]{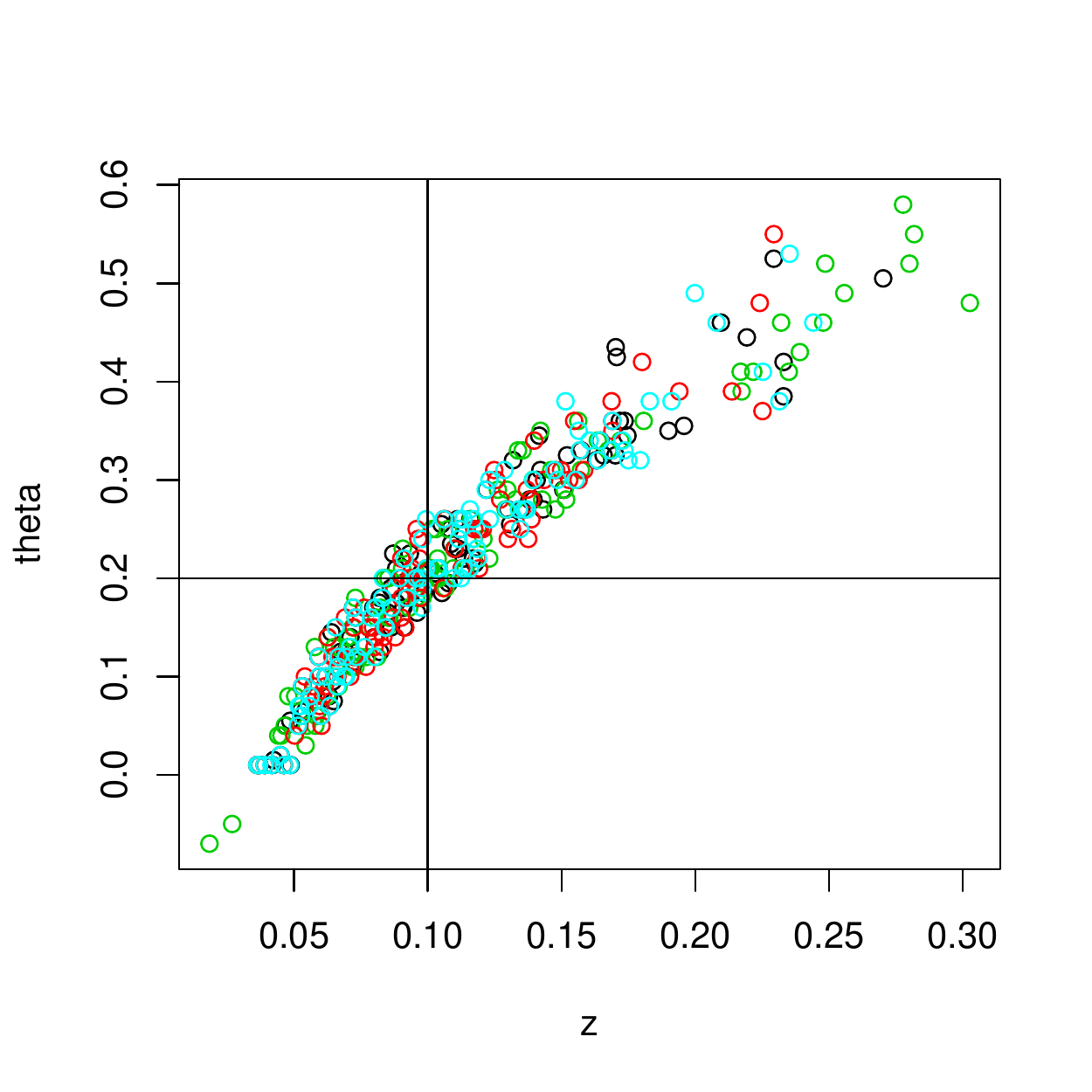} & \hspace{-1cm}
 \includegraphics[angle=0,scale=0.45]{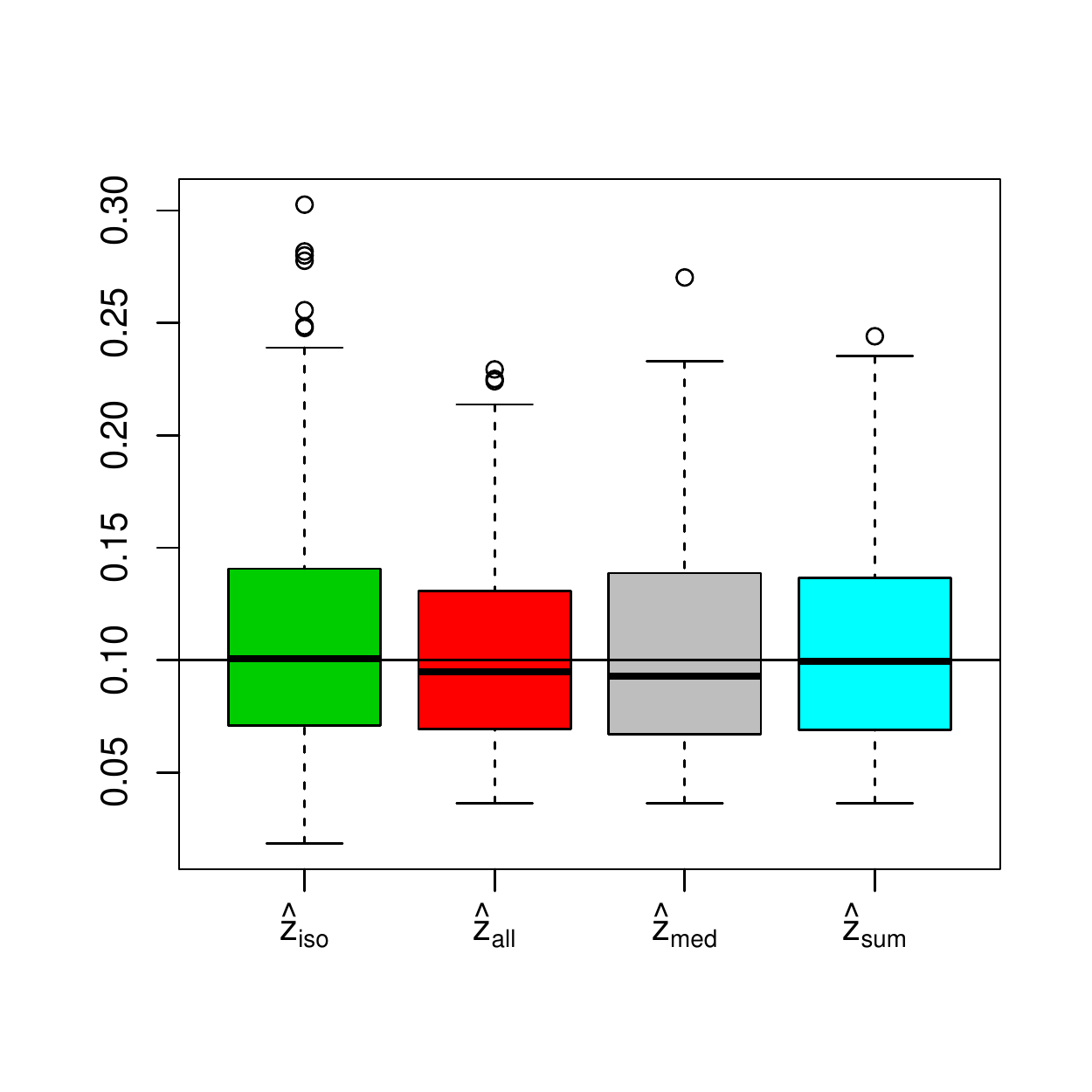} & \hspace{-1cm}
  \includegraphics[angle=0,scale=0.45]{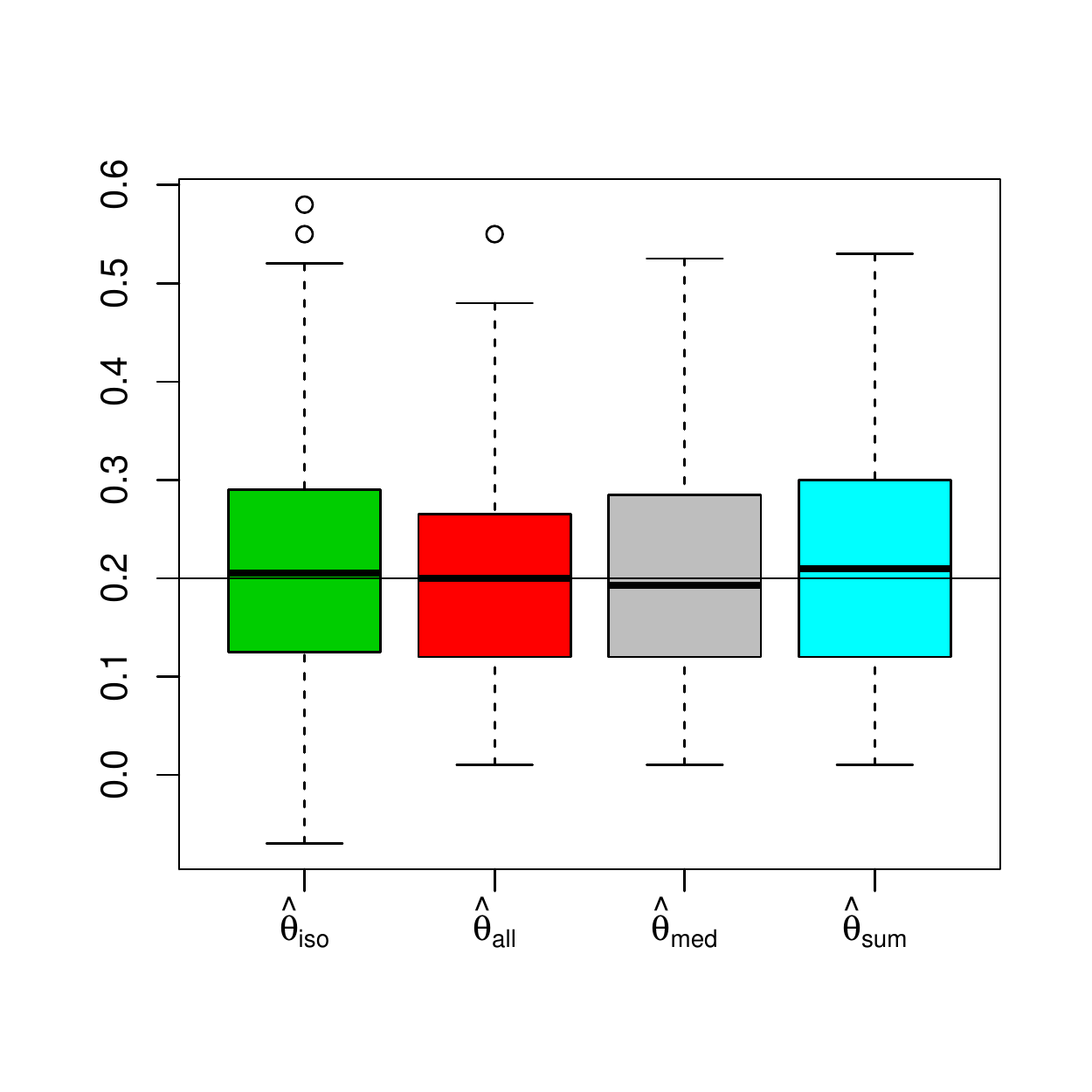} \\
 \includegraphics[angle=0,scale=0.45]{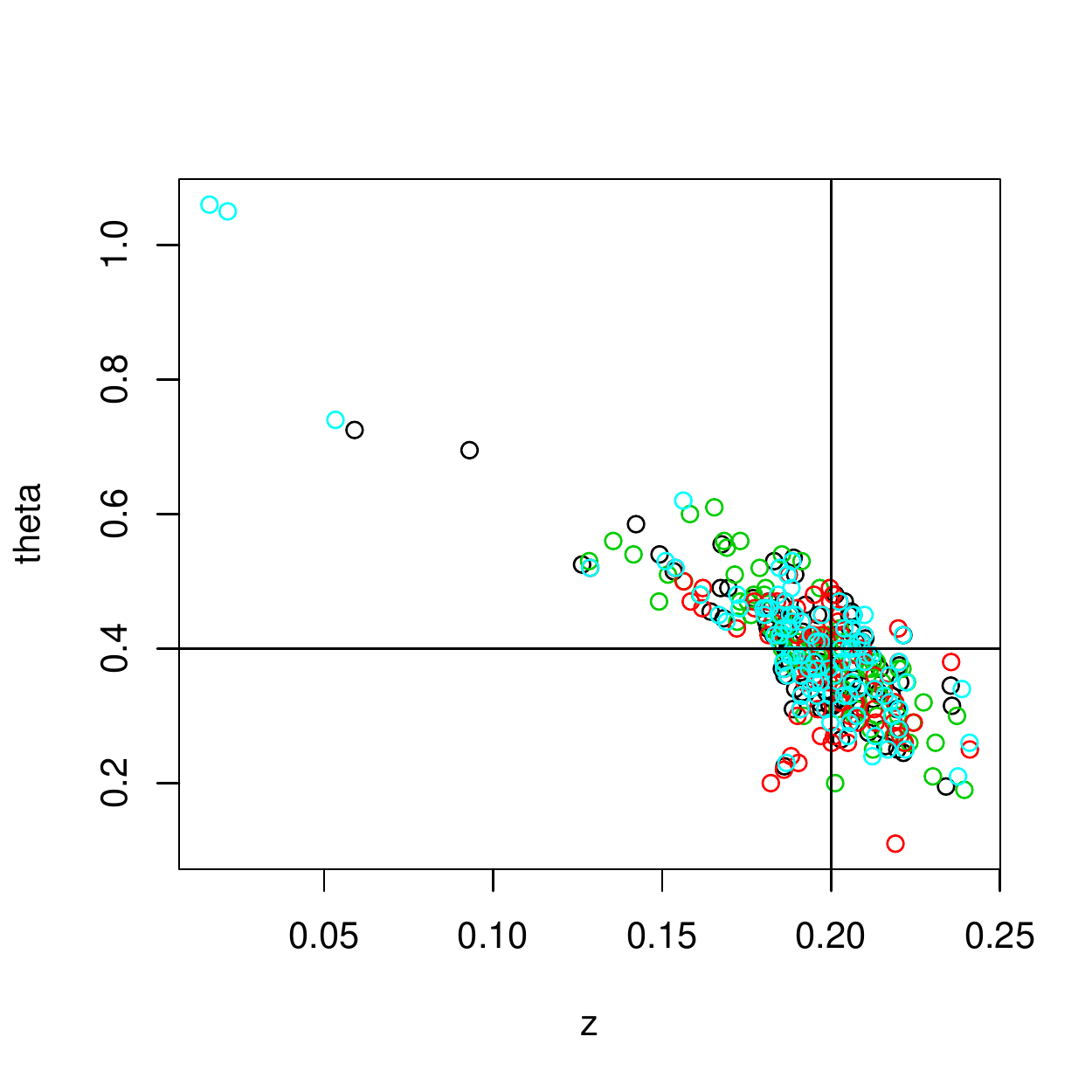} & \hspace{-1cm}
 \includegraphics[angle=0,scale=0.45]{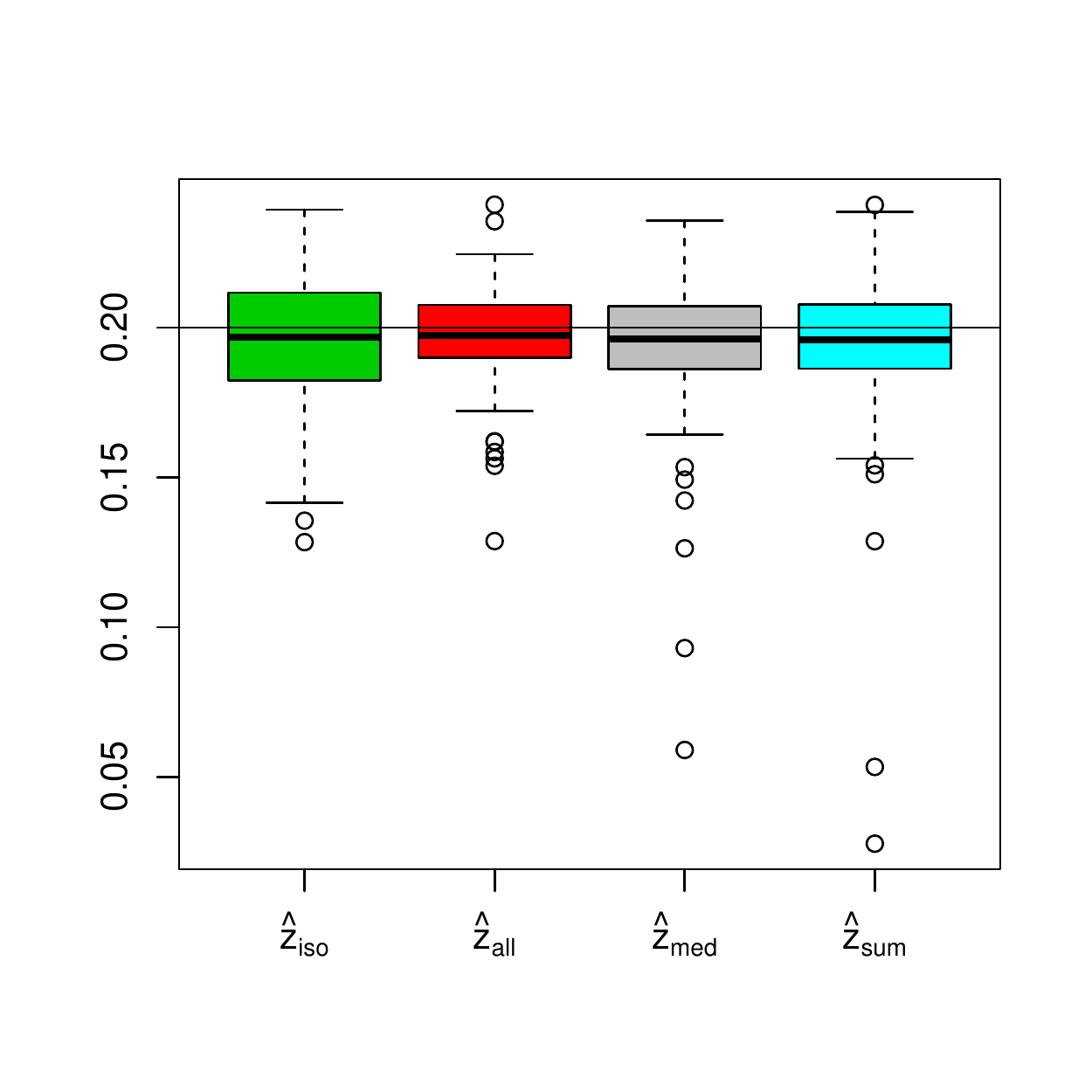} & \hspace{-1cm}
  \includegraphics[angle=0,scale=0.45]{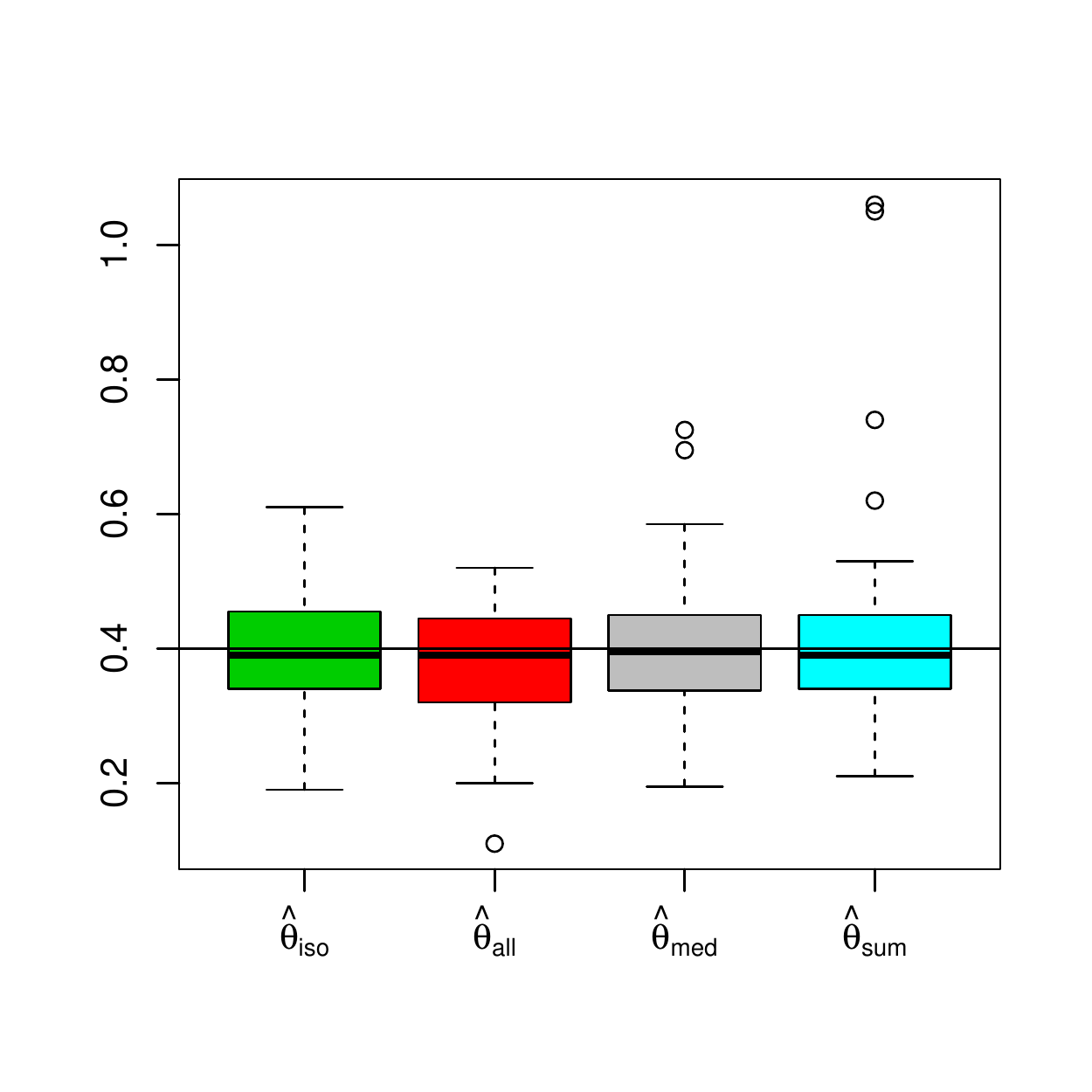} \\
 \includegraphics[angle=0,scale=0.45]{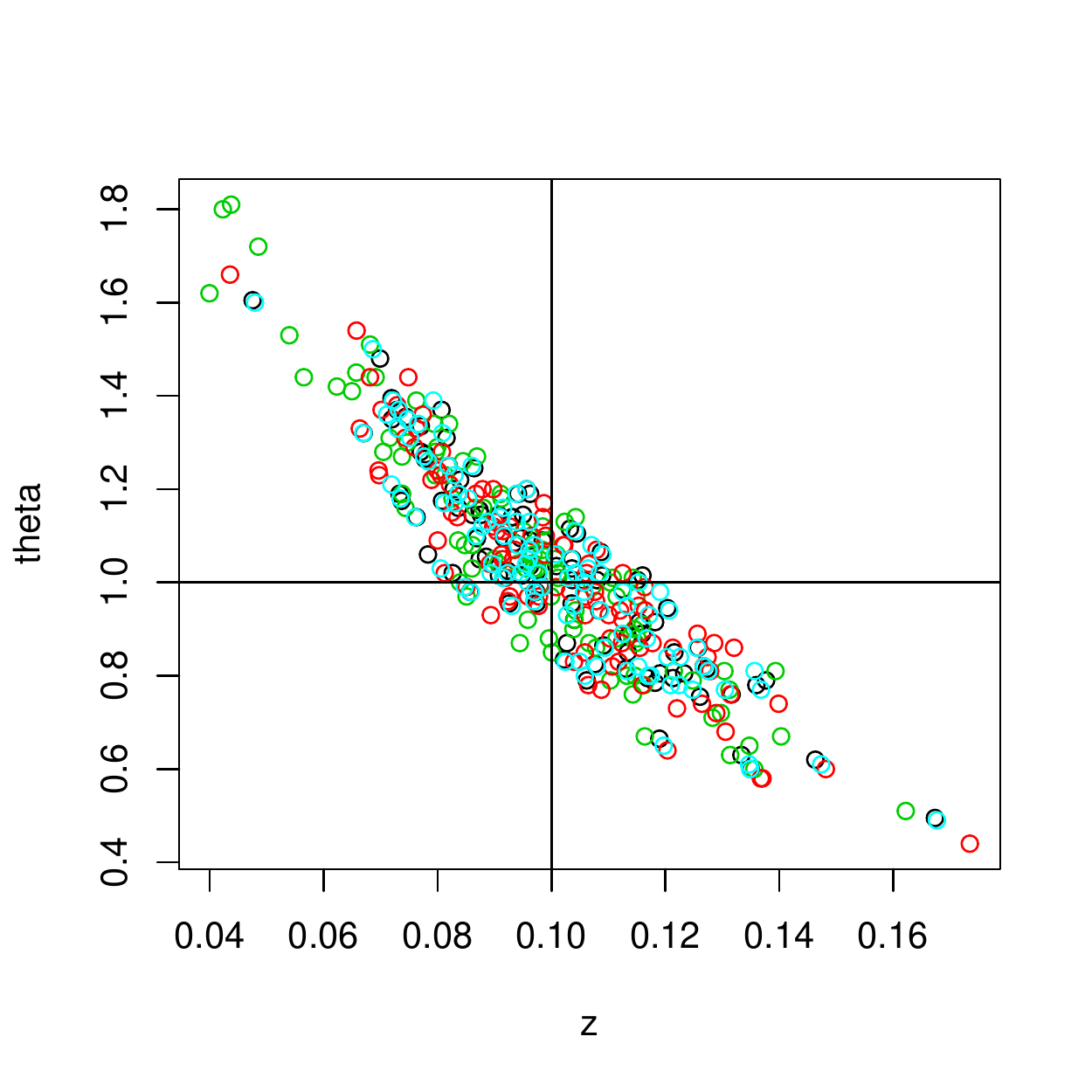} & \hspace{-1cm}
 \includegraphics[angle=0,scale=0.45]{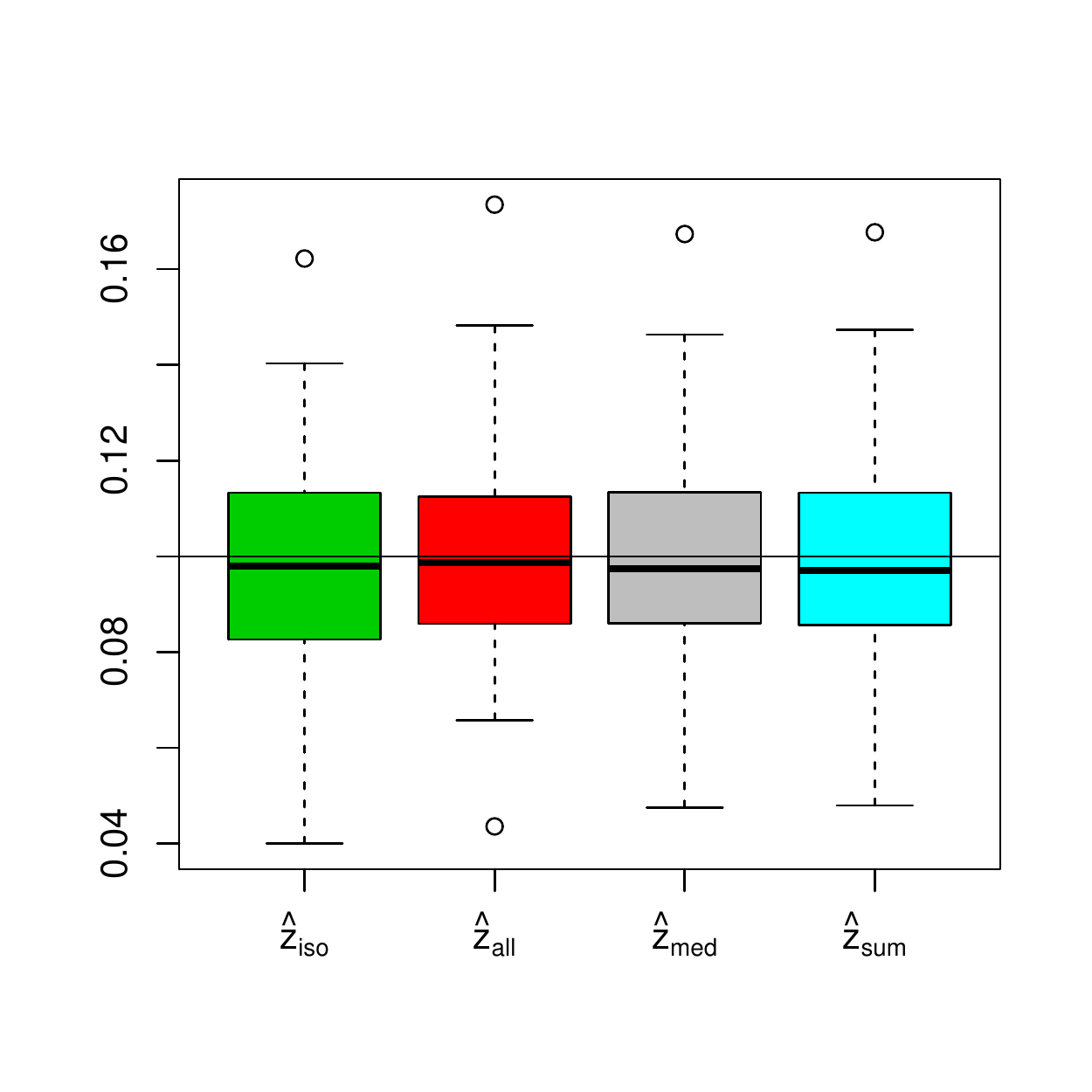} & \hspace{-1cm}
  \includegraphics[angle=0,scale=0.45]{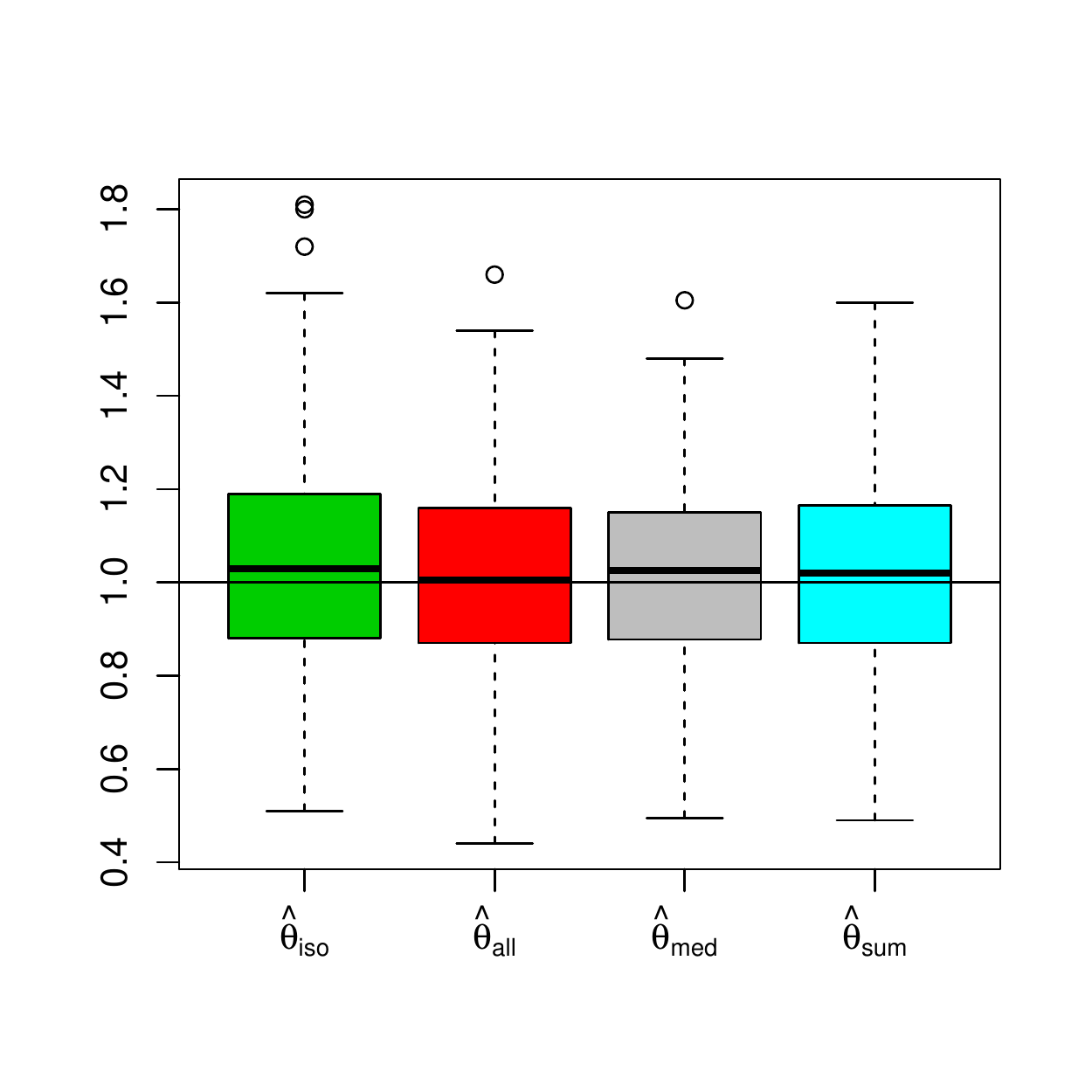} 
  \end{tabular}
 }
 \caption{Estimations from 100 replicates of, from top to bottom: An $\A$-process ($z^*=0.1$, $\t_1^*=0.2$);   An $\LL$-process ($z^*=0.2$, $\t_2^*=0.4$); A $\x$-process ($z^*=0.1$, $\t_3^*=1$). The estimators are  ($\hat z_{\text{iso}}, \hat \theta_{\text{iso}}$) (green),  ($\hat z_{\text{all}}, \hat \theta_{\text{all}}$) (red), $(\hat z_{\text{med}}, \hat \theta_{\text{med}})$ (black) and ($\hat z_{\text{sum}}, \hat \theta_{\text{sum}}$) (blue). From left to right: scatterplot of $(\hat z,\hat\t)$;  boxplots of $\hat z$; boxplots of $\hat \t$.}    \label{fig-area}
 \end{figure}

All simulations are done on $[0,50]^2$. The reference law of radii $\mu$ is uniform on $[0.5,2]$. The practical implementation is conducted as explained in Section~\ref{PA}, that includes minus sampling and the Monte Carlo approximation \eqref{Deltaapprox}.
Figure~\ref{fig-area} displays  the estimation results on 100 replicates of: The area-interaction process with $z^*=0.1$, $\t_1^*=0.2$, the perimeter-interaction process with $z^*=0.2$, $\t_2^*=0.4$, and the Euler-Poincaré-interaction process with $z^*=0.1$, $\t_3^*=1$,  respectively.

As a first conclusion, we see from Figure~\ref{fig-area} that all estimators  allow to estimate the two parameters and perform more or less equivalently. However,  ($\hat z_{\text{all}}, \hat \theta_{\text{all}}$) seems to slightly surpass the other estimators, in terms of dispersion and outliers. 

\begin{figure}[htbp]
  \centerline{
  \begin{tabular}[]{ccc}
 \includegraphics[angle=0,scale=0.45]{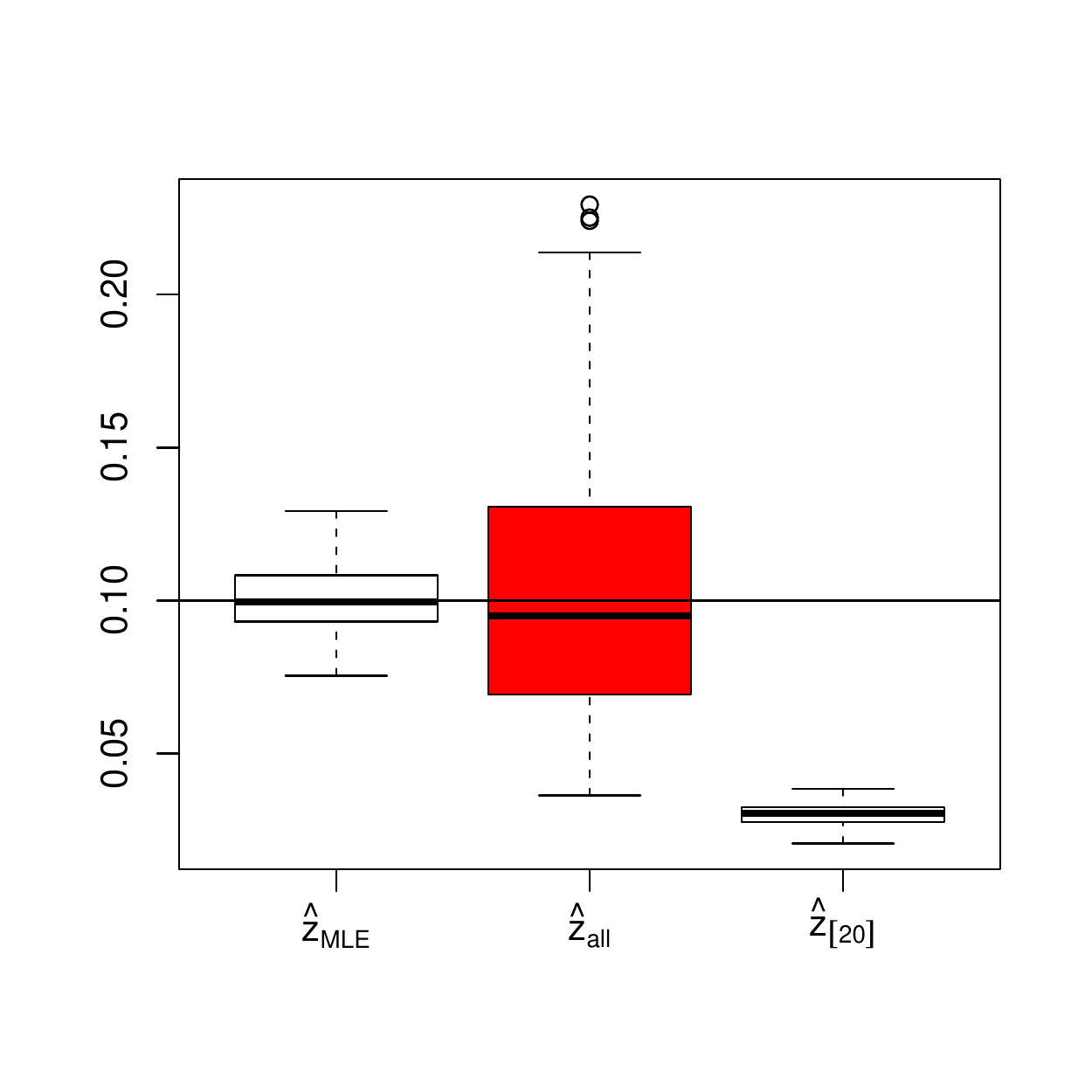} & \hspace{-1cm}
 \includegraphics[angle=0,scale=0.45]{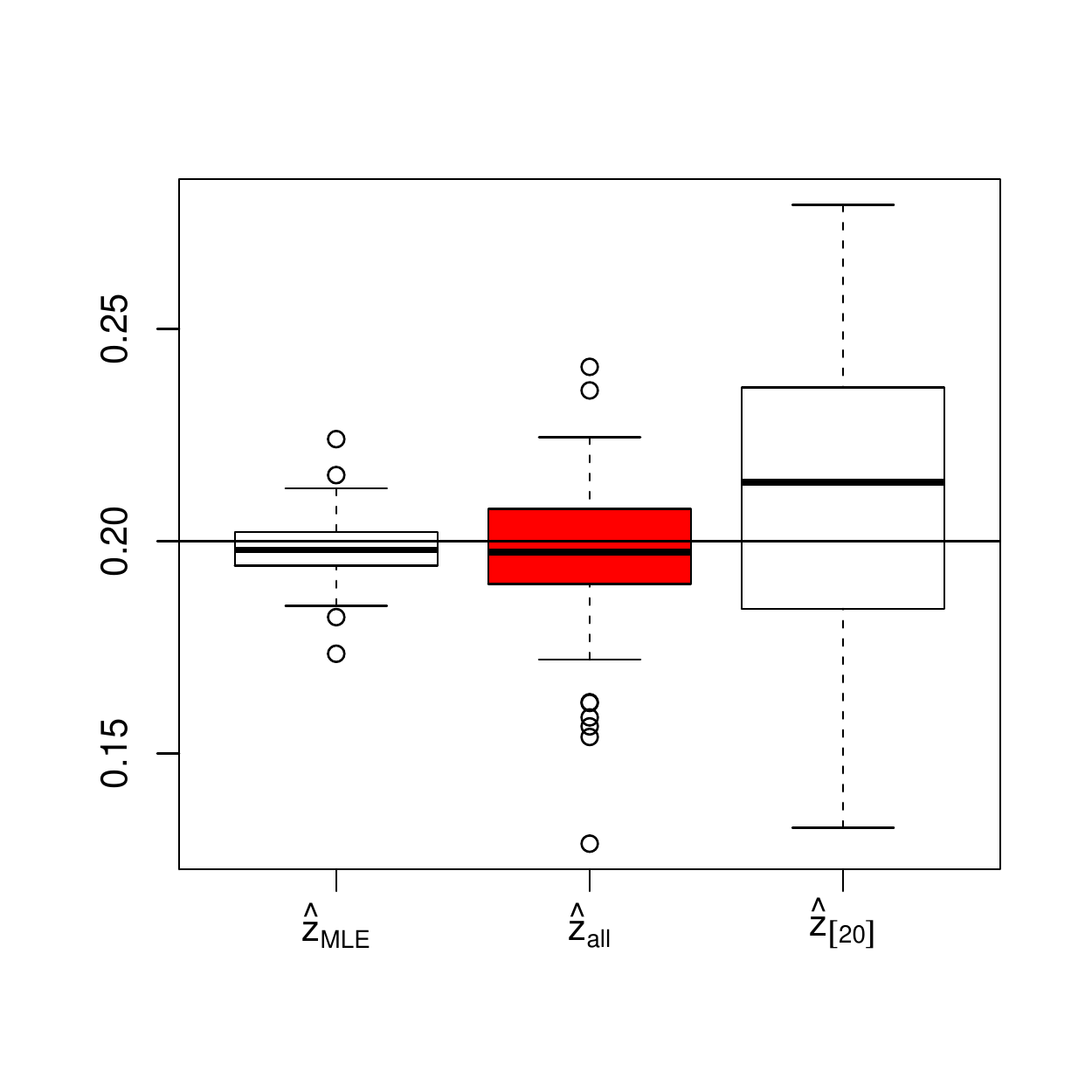} & \hspace{-1cm}
  \includegraphics[angle=0,scale=0.45]{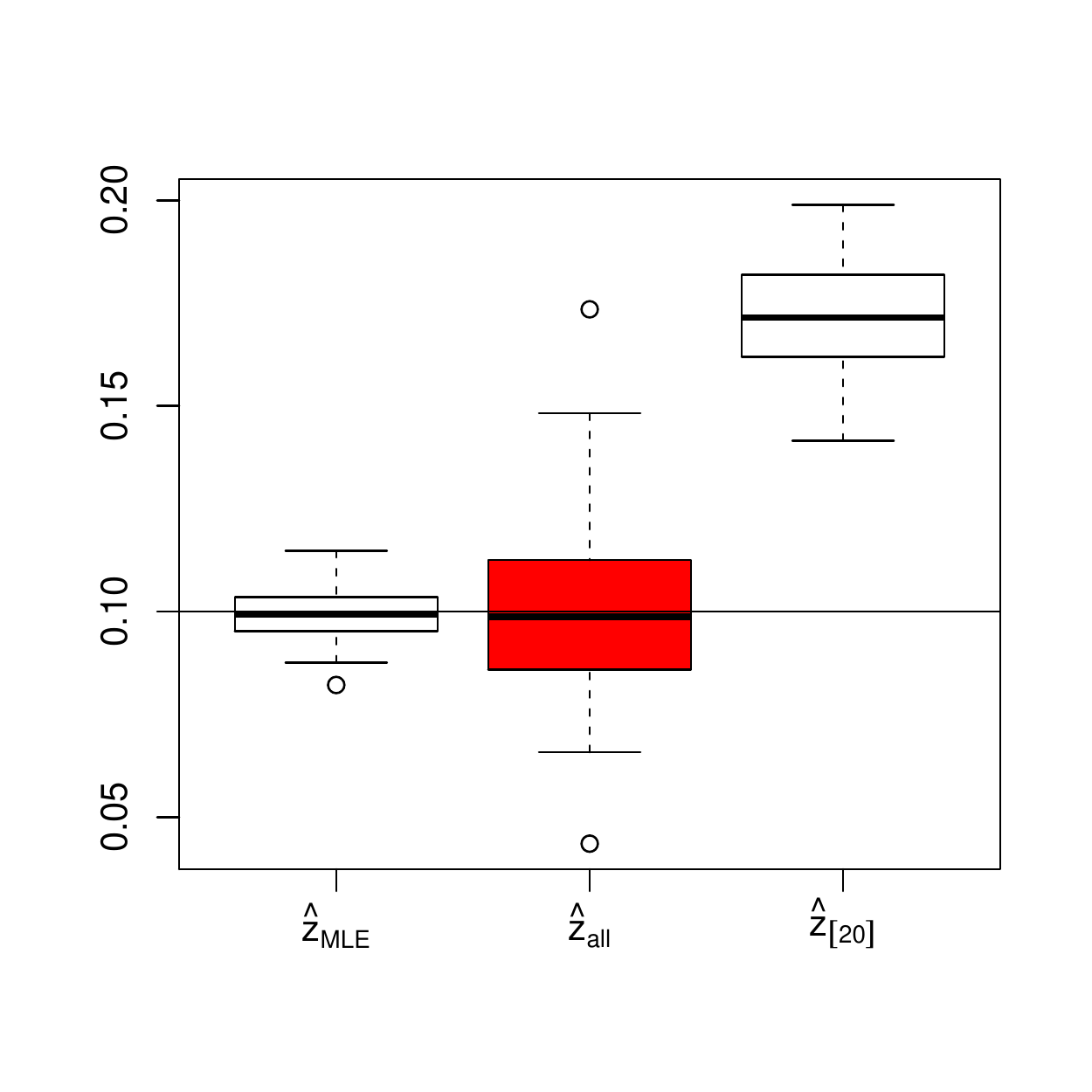} \\
  $\A$-process & $\LL$-process & $\x$-process \\
 \includegraphics[angle=0,scale=0.45]{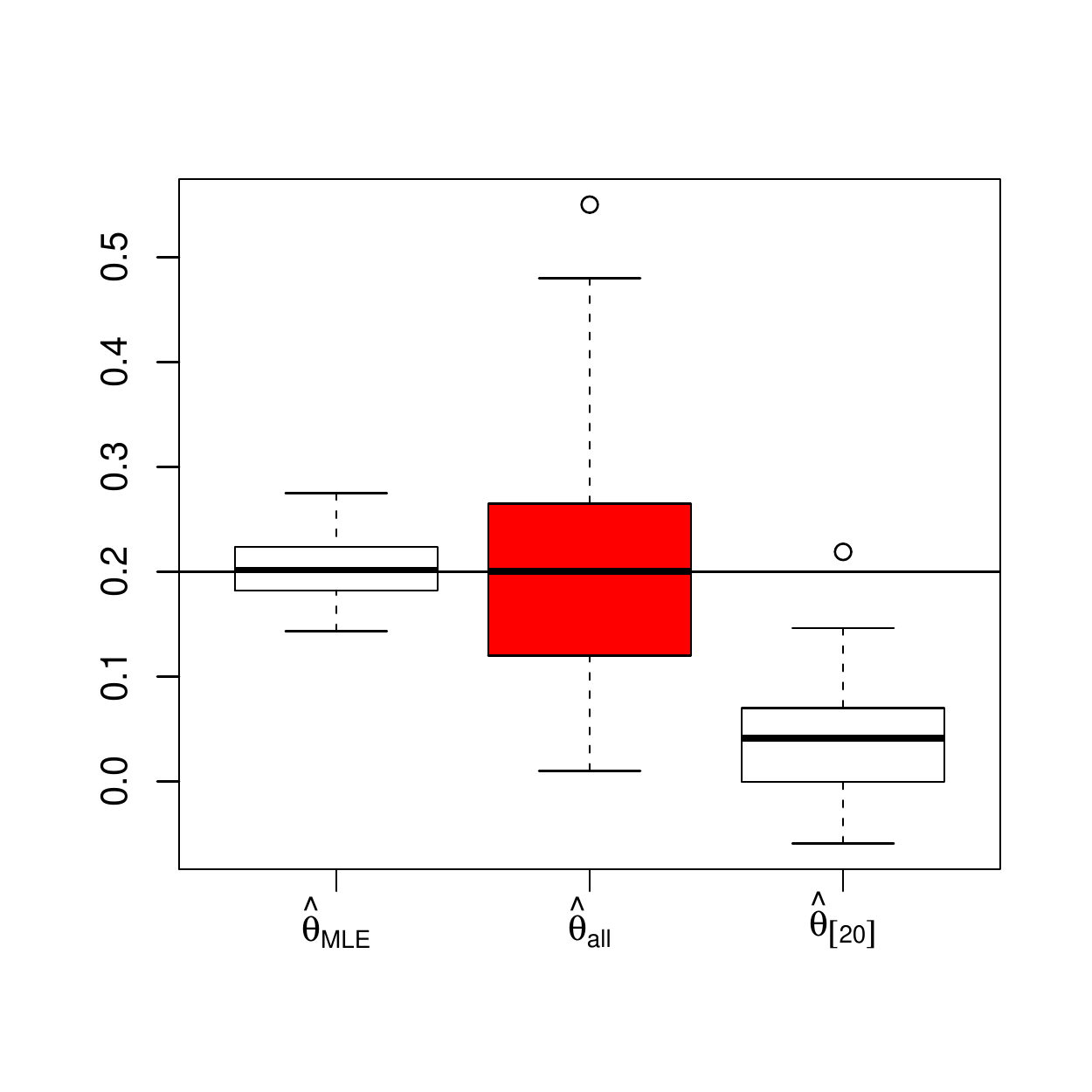} & \hspace{-1cm}
 \includegraphics[angle=0,scale=0.45]{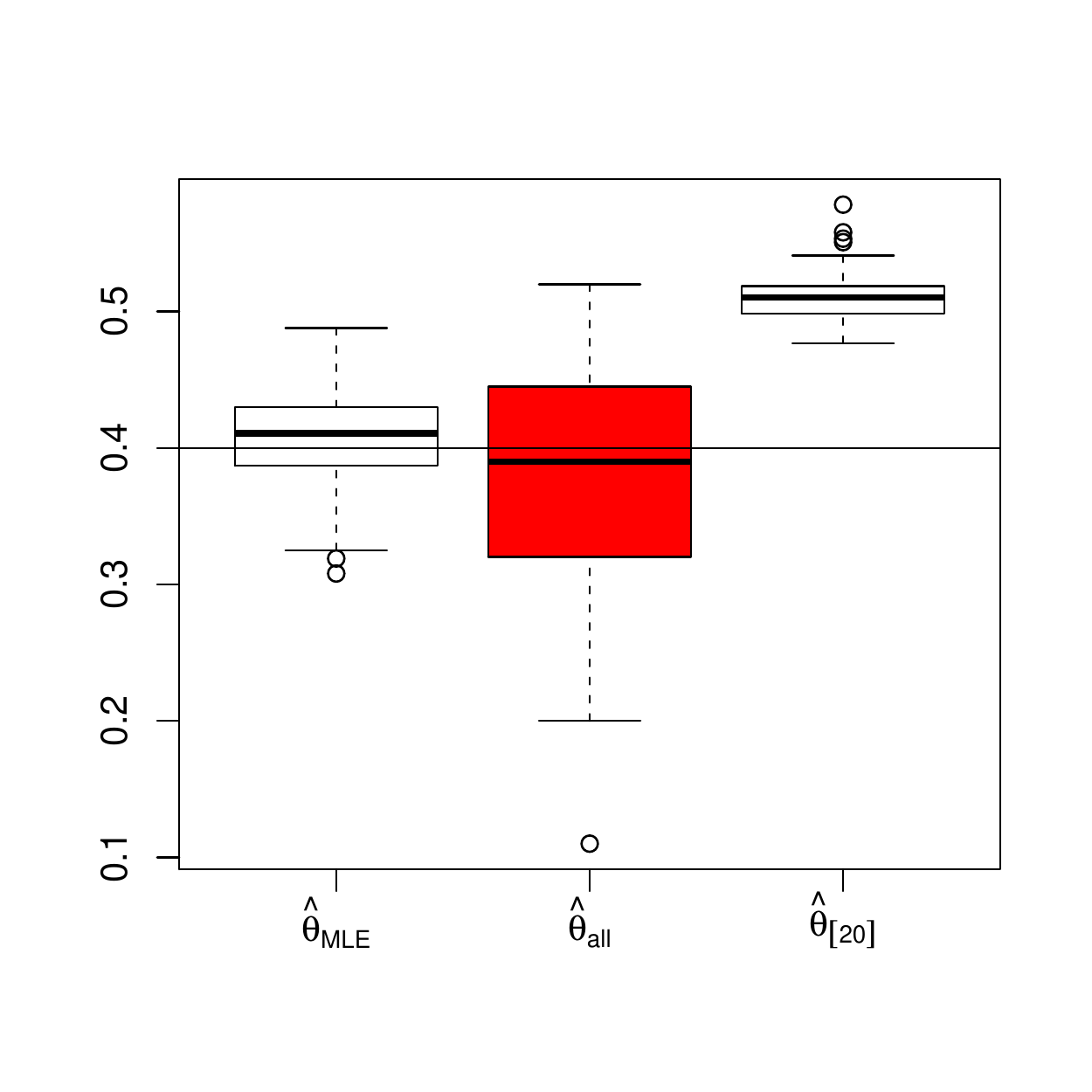} & \hspace{-1cm}
  \includegraphics[angle=0,scale=0.45]{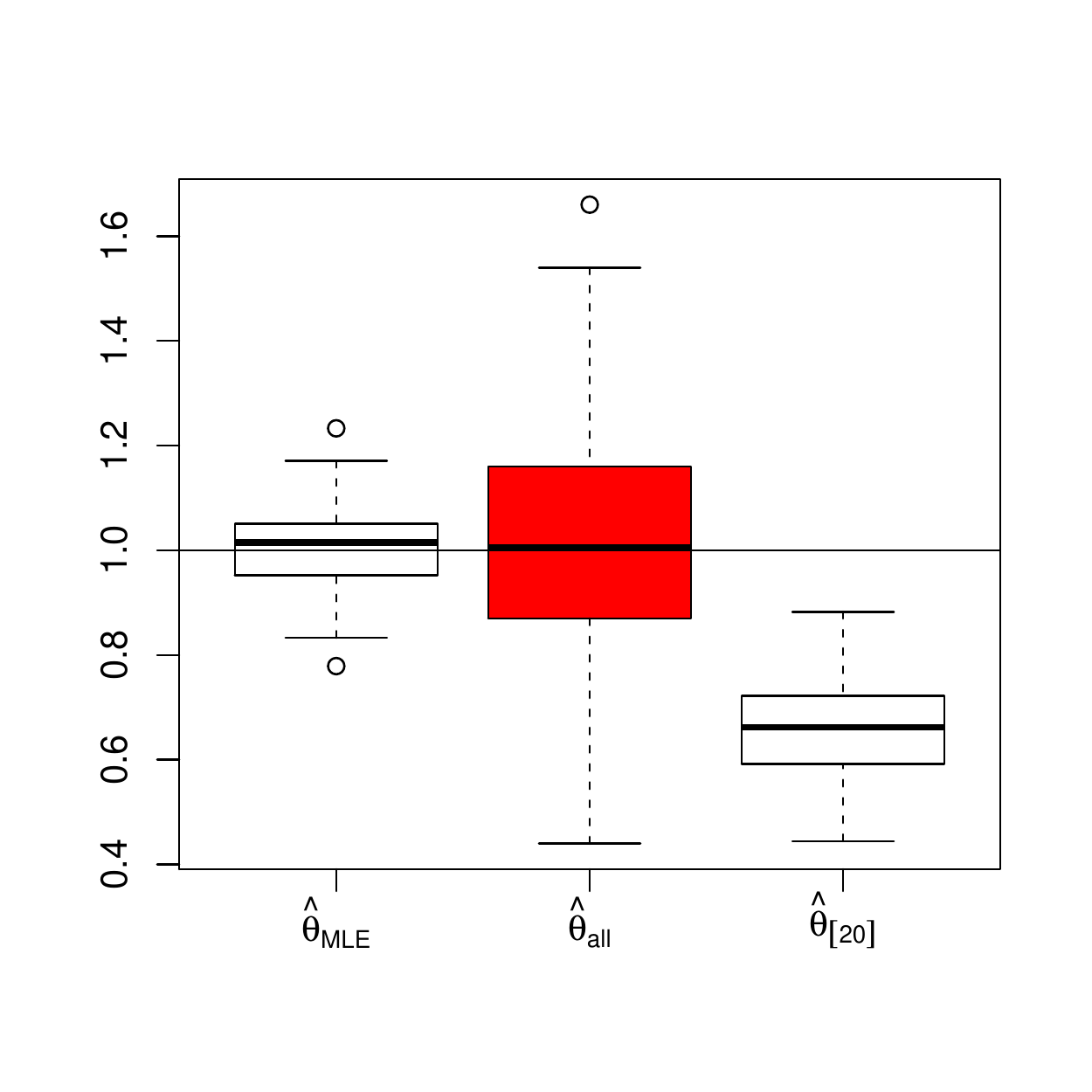} 
  \end{tabular}
 }
 \caption{Comparison of three estimators of $z$ (top) and $\theta$ (bottom) for the examples considered in Figure~\ref{fig-area}, i.e. 100 replicates of an $\A$-process, an $\LL$-process and a $\x$-process. 
 For each plot: boxplot of the MLE assuming the centres of the balls are known (left); $\hat z_{\text{all}}$ or  $\hat \theta_{\text{all}}$ (middle); the two-step estimator considered in \cite{MH2} (right).}    \label{fig-comparison}
 \end{figure}

Finally, in Figure~\ref{fig-comparison}, the fluctuations of  ($\hat z_{\text{all}}, \hat \theta_{\text{all}}$) are compared to the behaviour of the MLE ($\hat z_{\text{MLE}}, \hat \theta_{\text{MLE}}$),  for the three examples explained above. The MLE is computed in assuming that the centres of the balls are known (recall that this is not possible  in practice for a germ-grain set). So Figure~\ref{fig-comparison} shows what we lose when we do not observe the germs of the Quermass-interaction process and we use the TF procedure. On the right hand side of each plot of  Figure~\ref{fig-comparison}, the boxplot of the two-step estimator ($\hat z_{[20]}, \hat \theta_{[20]}$) considered in \cite{MH2}   is also represented (see Section~\ref{MLE} for a discussion). The strong bias of  this procedure is clearly illustrated.

\subsection{Estimation of several interaction parameters}\label{estimation-all}

In the previous section, we have considered models with two parameters : the intensity and one interaction parameter.
To fit  models with more interaction parameters, we follow the conclusion of the previous section: Given ten small $\alpha_i$'s,  we consider the TF estimator based on all available test functions $f_0$, $f_{\alpha_1}$, $\cdots$, $f_{\alpha_{10}}$, $f_{\text{sum}}$ and  $f_{\text{iso}}$. Since $f_{\text{iso}}$ is strongly related to the ball structure of the germ-grain model, which might be a restriction for practical computation, we also consider the TF estimator based on all previous test functions except $f_{\text{iso}}$.

\begin{figure}[H]
    \setlength{\tabcolsep}{0.1cm} \centerline{
  \begin{tabular}[]{cc}
    \includegraphics[angle=0,scale=.2]{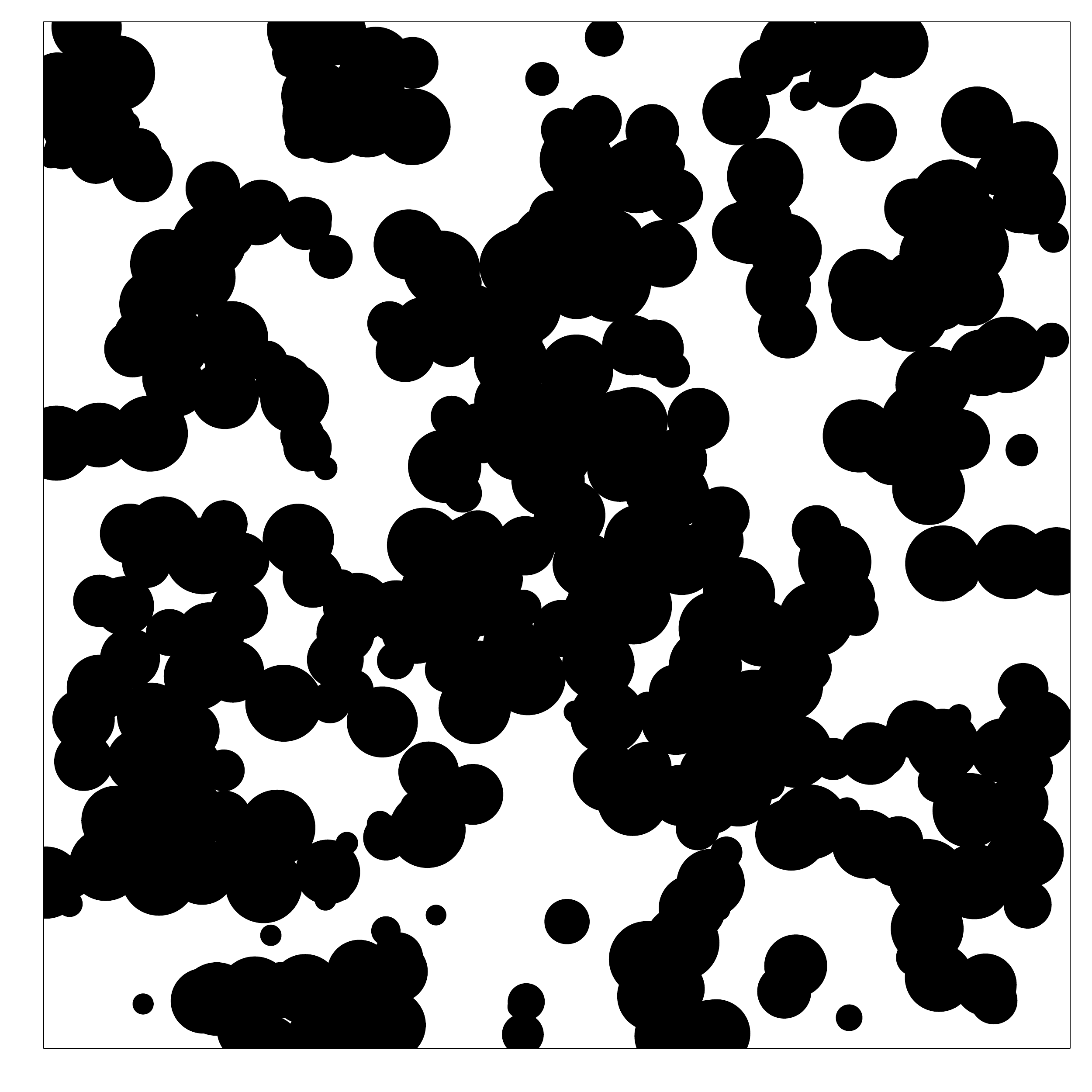}&\hspace{1cm}
  \includegraphics[angle=0,scale=.2]{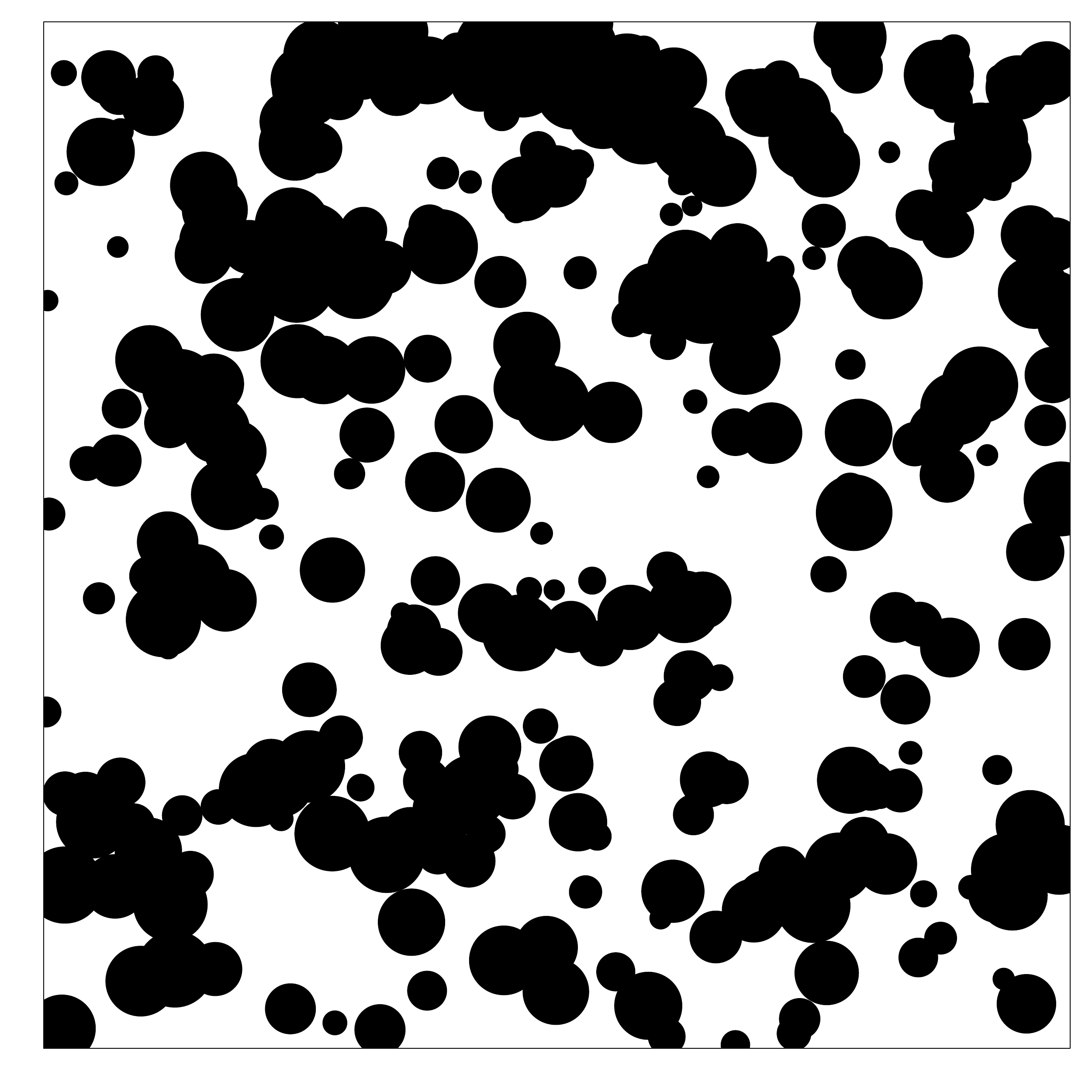}
  \end{tabular}
  }
      \caption{{\small Samples of the $(\A,\LL)$-process with  $z^*=0.1$,  $\t_1^*=-0.2$ and $\t_2^*=0.3$ (left) and of  the Quermass process with $z^*=0.1$, $\t_1^*=-0.2$, $\t_2^*=0.3$ and $\t_3^*=-1$ (right). The window is $[0,50]^2$ and $\mu$ is the uniform law on  $[0.5,2]$.}}\label{samples-all}
 \end{figure}
\begin{figure}[H]
  \setlength{\tabcolsep}{0.1cm} \centerline{
  \begin{tabular}[]{ccc}
 \includegraphics[angle=0,scale=.45]{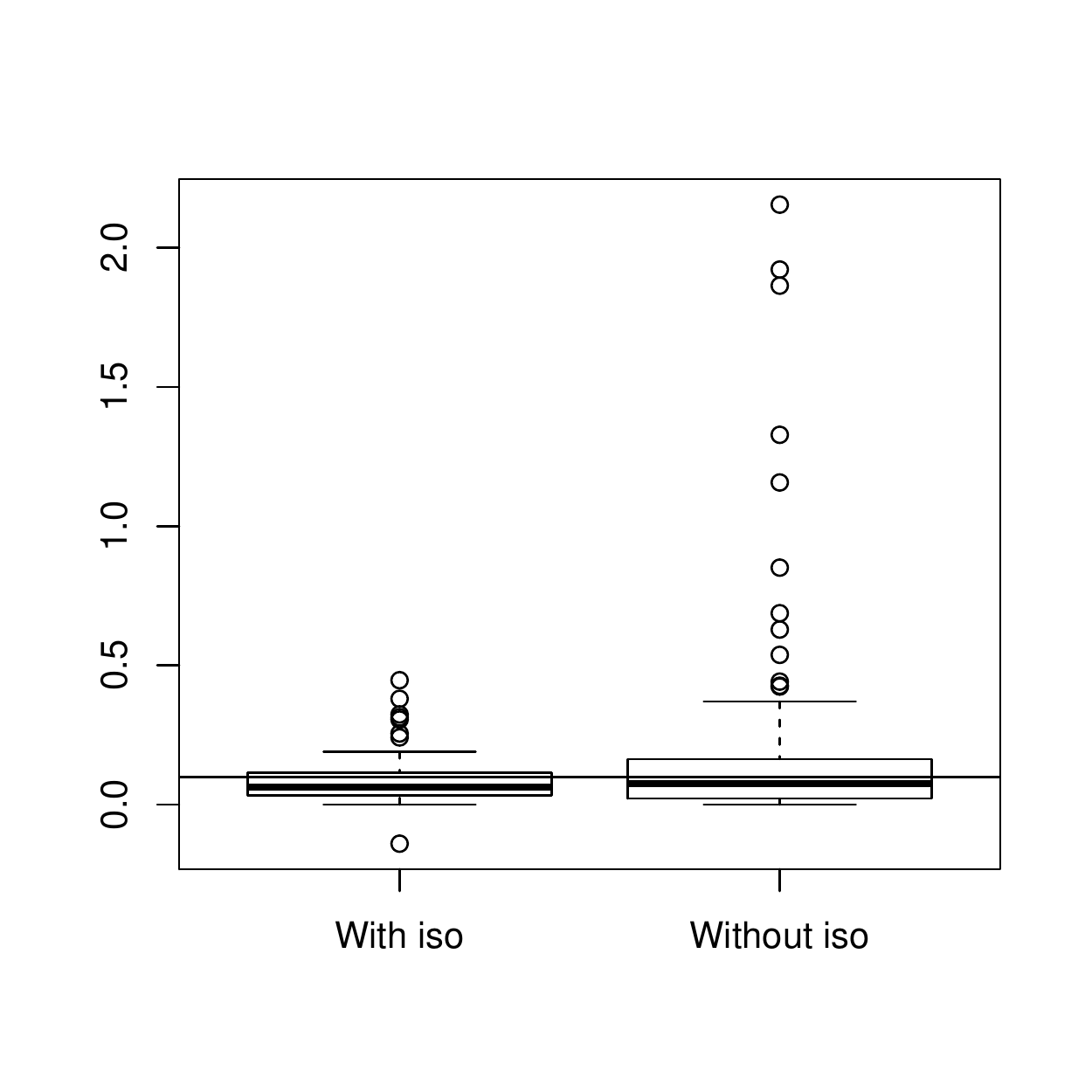} &\hspace{-1cm}
 \includegraphics[angle=0,scale=.45]{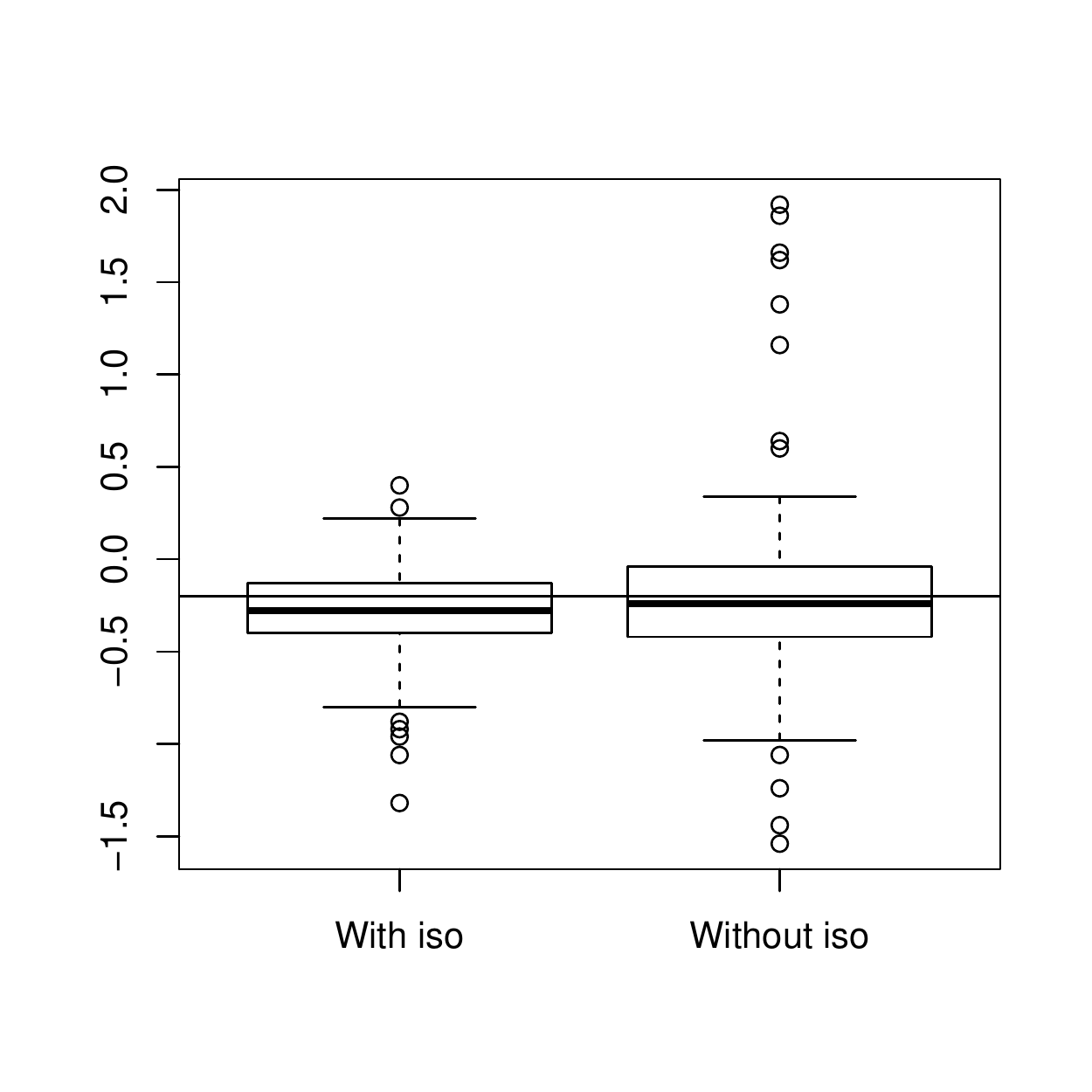}  &\hspace{-1cm}
 \includegraphics[angle=0,scale=.45]{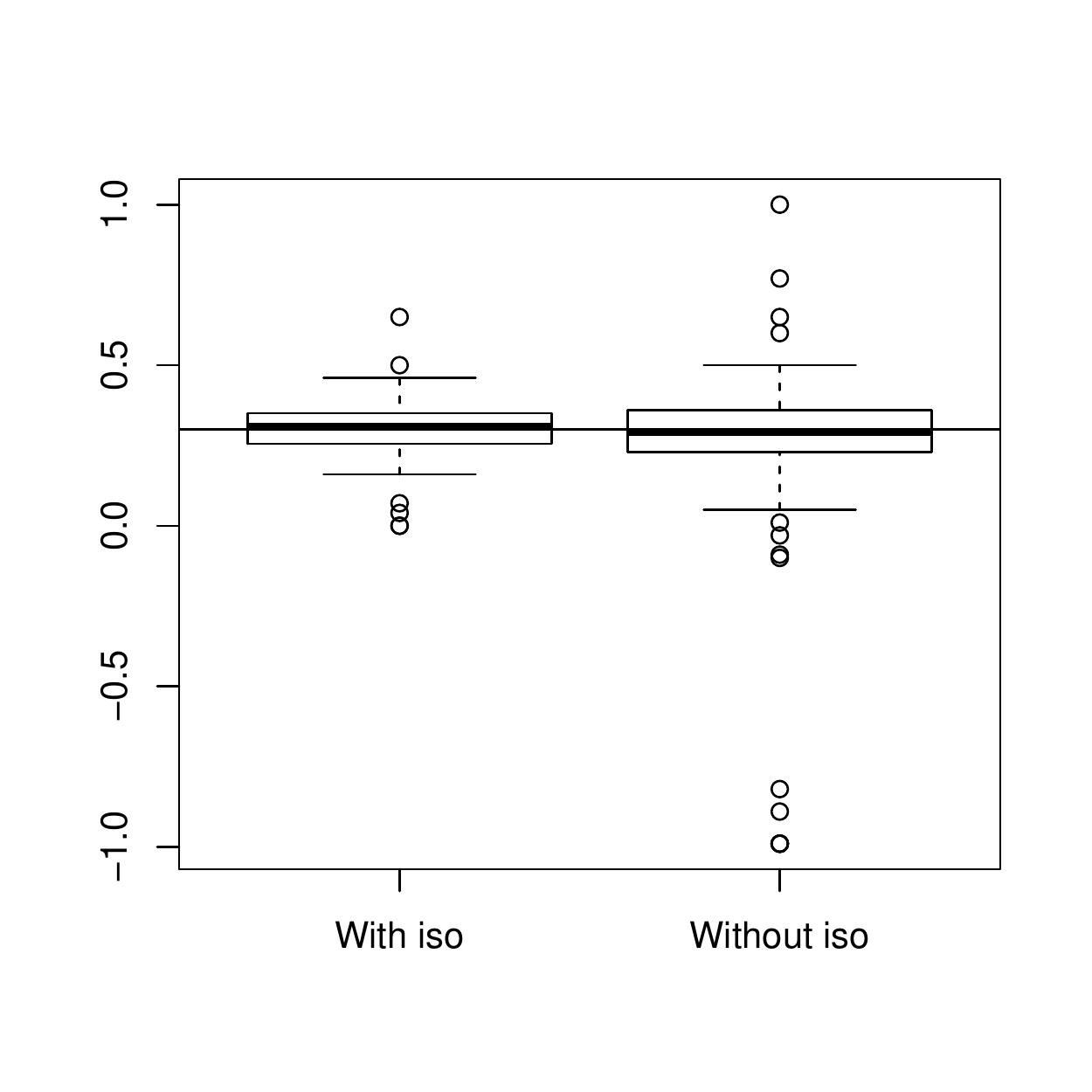} \end{tabular}
  }
  \caption{Estimation of $z$ (left), $\t_1$ (middle) and $\t_2$ (right) from 100 replicates of the $(\A,\LL)$-process with $z^*=0.1$, $\t_1^*=-0.2$ and $\t_2^*=0.3$. In each plot, the estimator involves  $f_{\text{iso}}$ (on the left) or not (on the right).}\label{fig-area-perim} 
 \end{figure}
\begin{figure}[H]
  \setlength{\tabcolsep}{0.1cm} \centerline{
  \begin{tabular}[]{cccc}
 \includegraphics[angle=0,scale=.3]{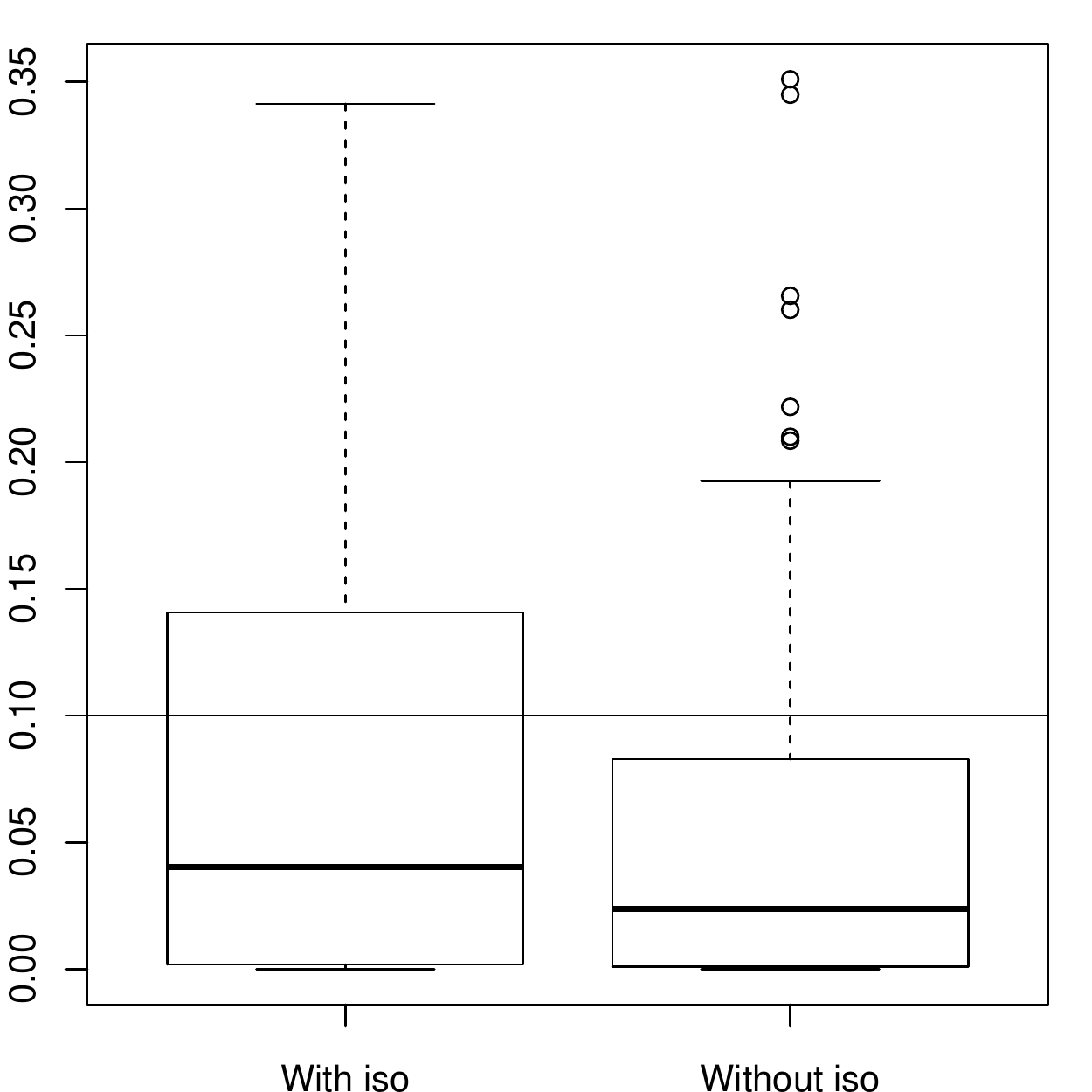} &
 \includegraphics[angle=0,scale=.3]{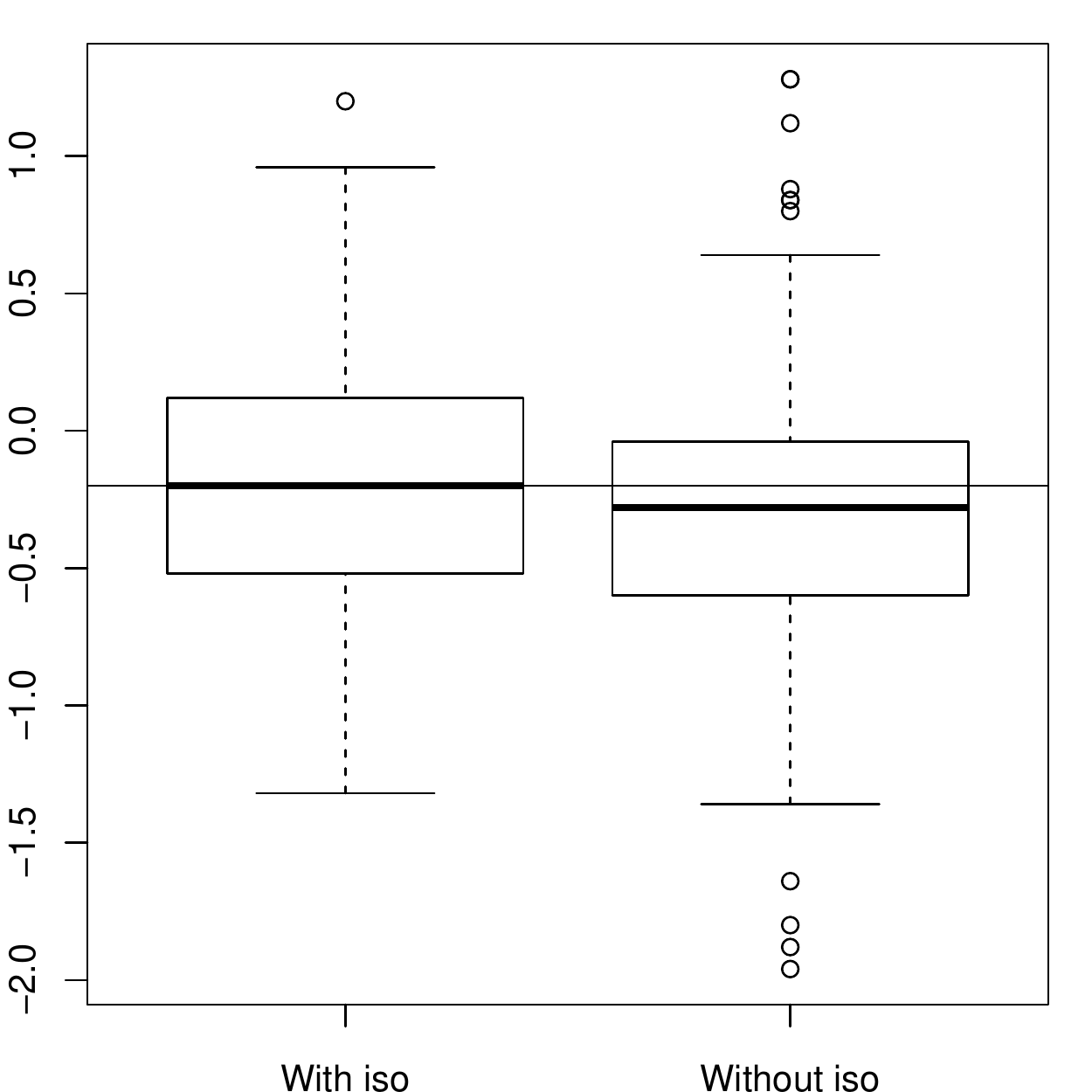}&
 \includegraphics[angle=0,scale=.3]{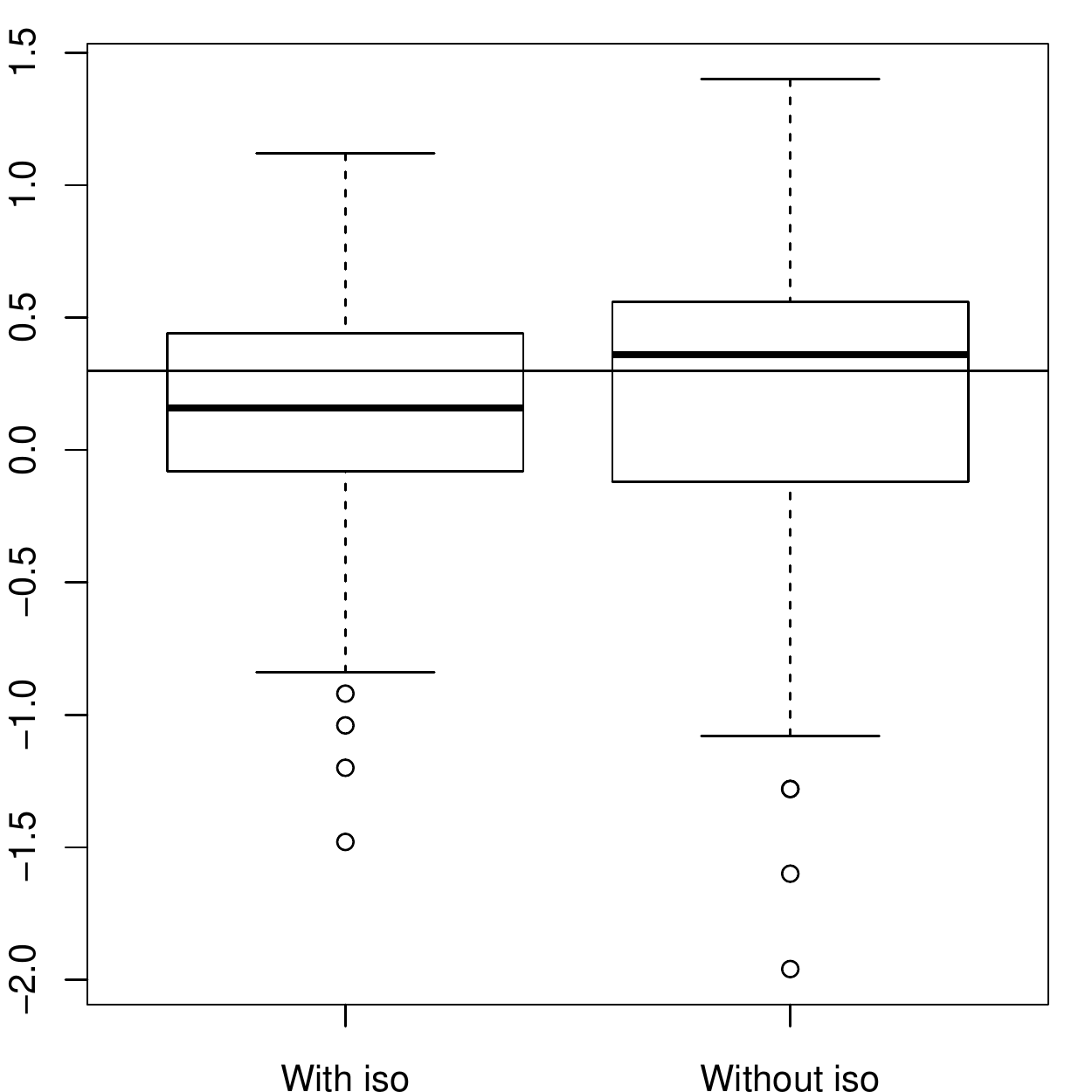} &
 \includegraphics[angle=0,scale=.3]{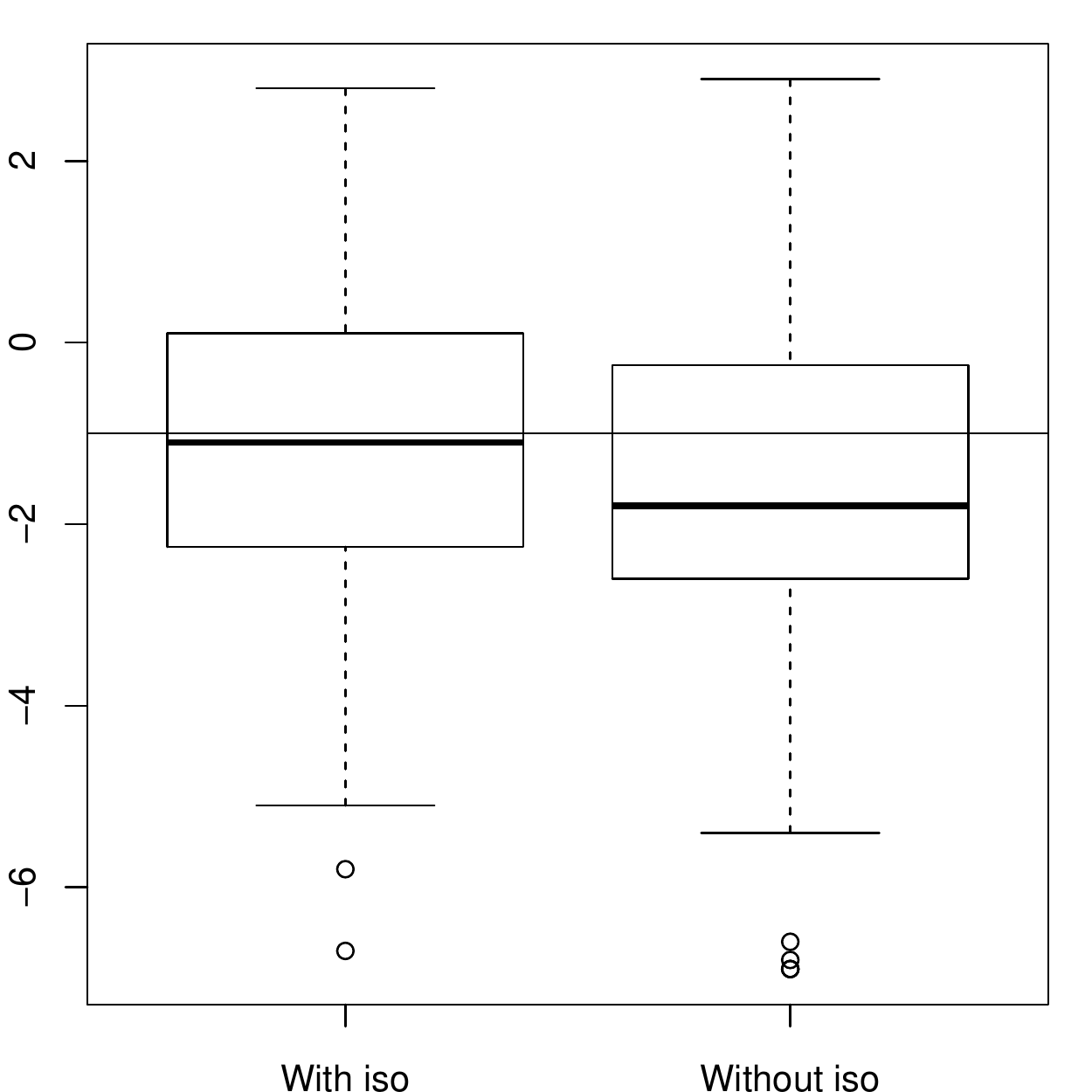}  \end{tabular}
  }
  \caption{Estimation (from left to right) of $z$, $\t_1$, $\t_2$ and $\t_3$ from 100 replicates of the Quermass-interaction process with $z^*=0.1$, $\t_1^*=-0.2$, $\t_2^*=0.3$ and $\t_3^*=-1$. In each plot, the estimator involves  $f_{\text{iso}}$ (on the left) or not (on the right).}\label{fig-quermass} 
 \end{figure}

Some results of estimation are presented in Figure~\ref{fig-area-perim} for 100 replicates of the $(\A,\LL)$-process ($\t_3=0$) on $[0,50]^2$, with $z^*=0.1$, $\t_1^*=-0.2$, $\t_2^*=0.3$ and the reference law of radii is uniform on $[0.5,2]$.  The same kind of results are displayed in Figure~\ref{fig-quermass} for the Quermass-interaction process with $z^*=0.1$, $\t_1^*=-0.2$, $\t_2^*=0.3$ and $\t_3^*=-1$. See some examples of samples in Figure~\ref{samples-all}.

While the quality of estimation turns out to be satisfactory in the presence of one interaction parameter (i.e. for the $\A$-process, the $\LL$-process or the $\x$-process) as seen in Section~\ref{one-interaction}, it naturally decreases in the presence of more interaction parameters (Figure~\ref{fig-area-perim}), to become poor for the full Quermass-interaction model (Figure~\ref{fig-quermass}).  Moreover,  Figures~\ref{fig-area-perim} and \ref{fig-quermass} show that using the test function $f_{\text{iso}}$ seems to improve the estimation, but  its computation in practice requires to decide which components can be considered as isolated balls.

\section{Application to heather data}\label{sec:heather}

\begin{figure}[htbp]
  \setlength{\tabcolsep}{0.1cm} \centerline{
  \begin{tabular}[]{cc}
 \includegraphics[angle=0,scale=.3]{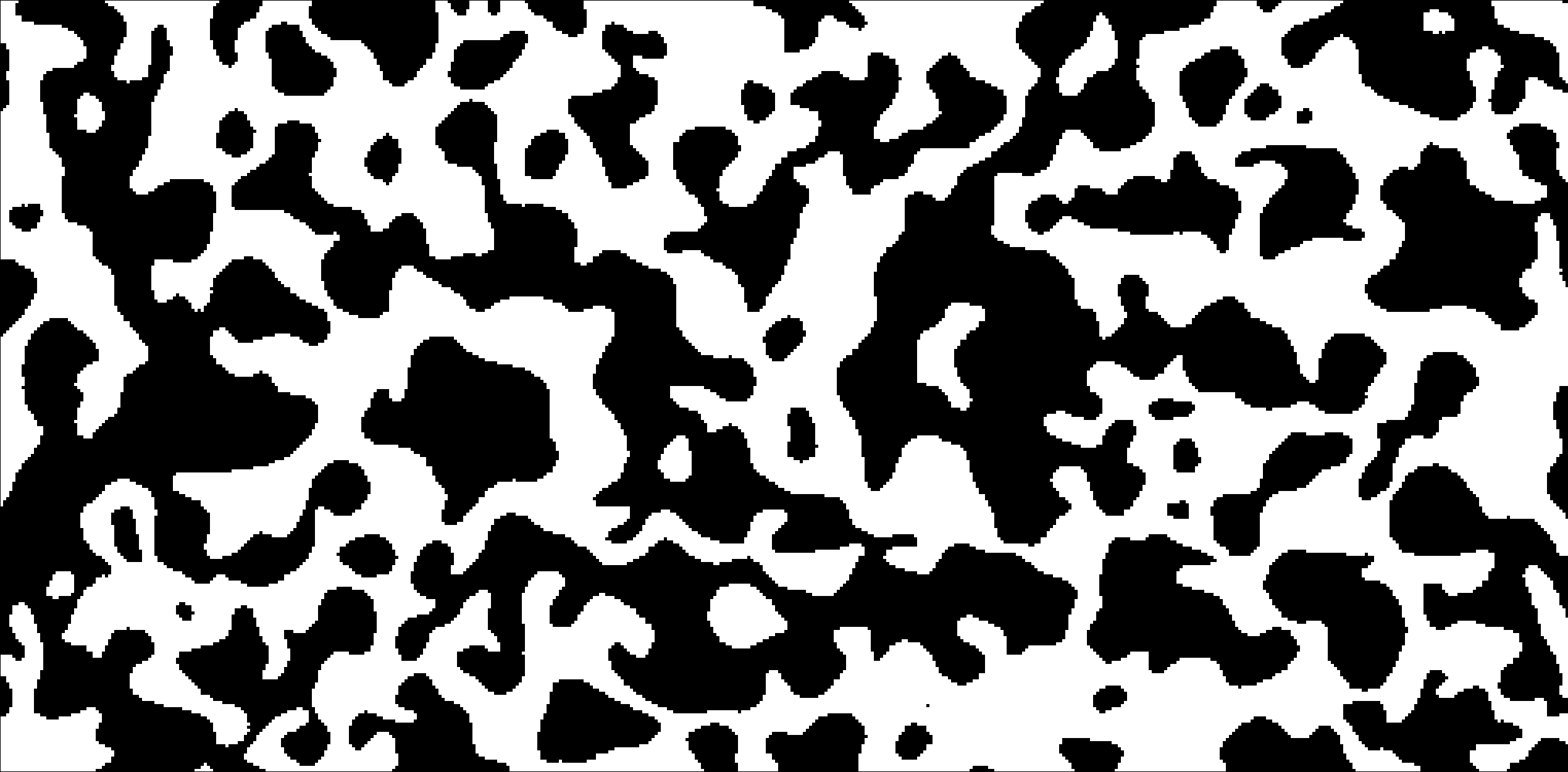} & 
  \includegraphics[angle=0,scale=.3]{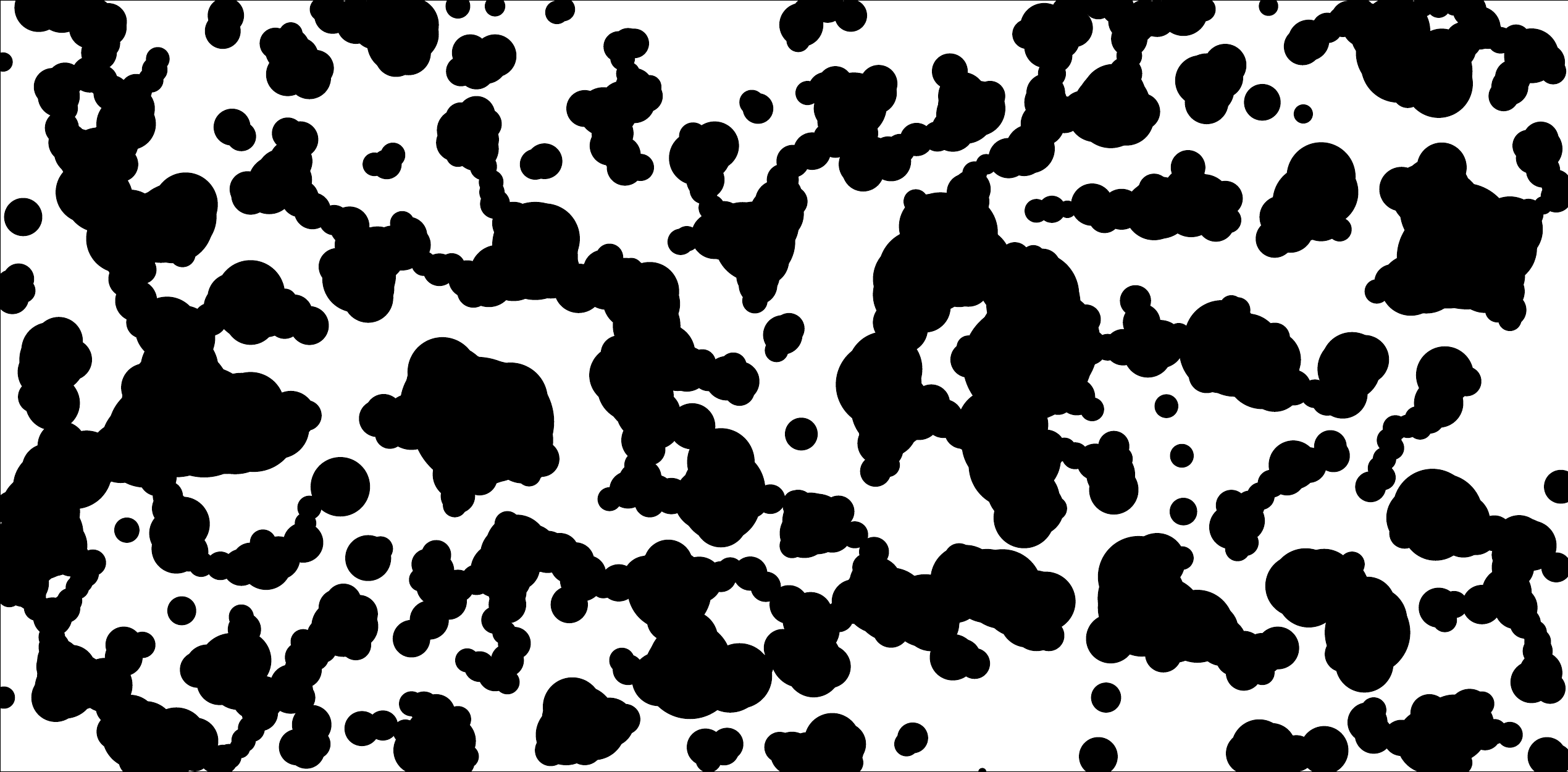} 
   \end{tabular}
  }
  \caption{Heather data (left) and its approximation by a union of balls (right)}\label{fig:heather}
 \end{figure}

The plot on the left hand side of Figure~\ref{fig:heather} shows a binary image of the presence of heather in a $10\times 20$ m rectangular region at J\"adra\r{a}s, Sweden. This heather dataset has been widely studied. It was first presented by P. Diggle in \cite{Diggle81}, where it was modelled by a stationary spherical Boolean model. This point of view has been considered further in \cite{Hall85, Hall88} and in \cite{Cressie}, where alternative estimation methods for the Boolean model have been implemented.

\begin{figure}[htbp]
  \setlength{\tabcolsep}{0.1cm} \centerline{
  \begin{tabular}[]{cc}
 \includegraphics[angle=0,scale=.35]{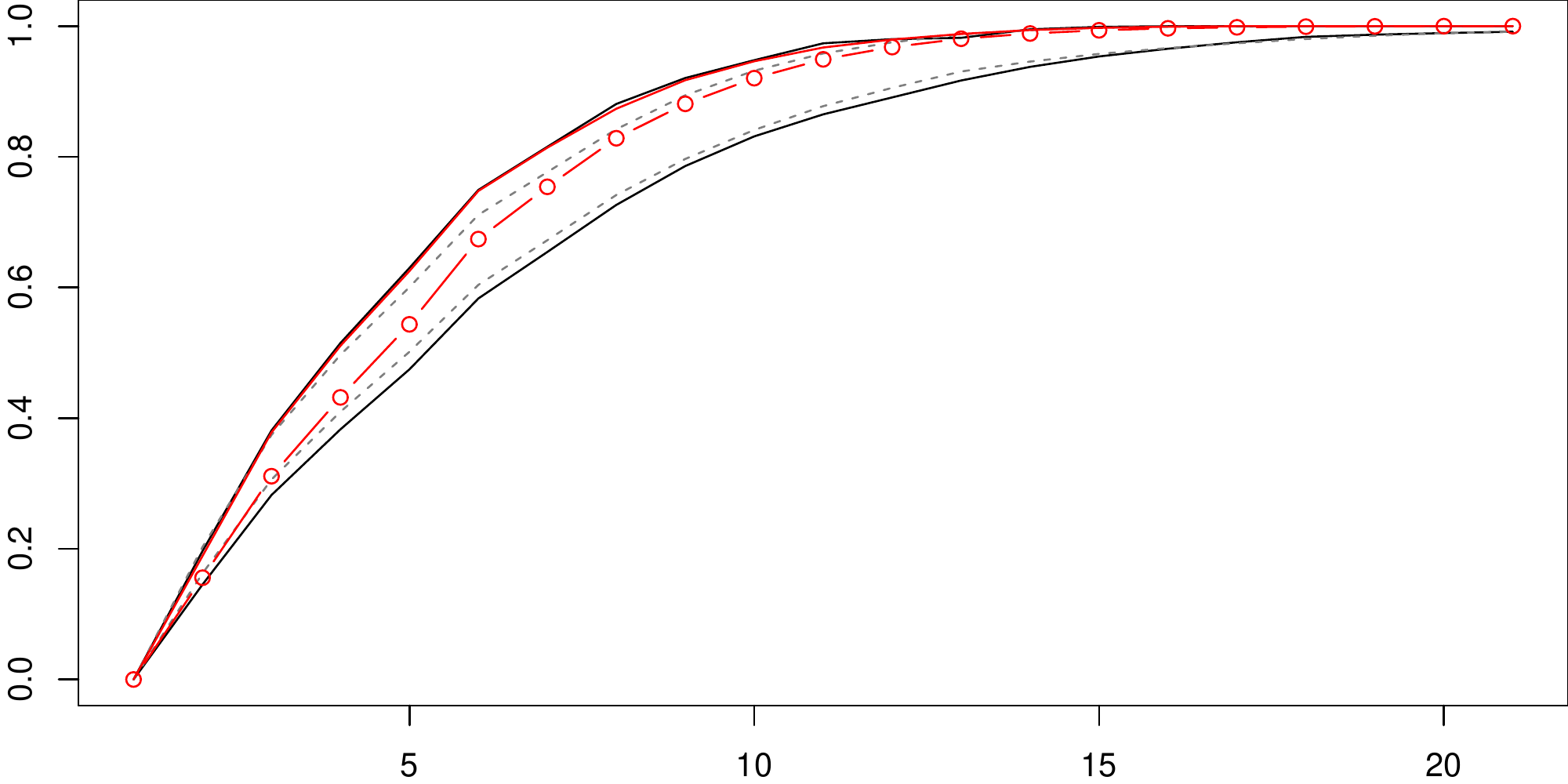} & \includegraphics[angle=0,scale=.35]{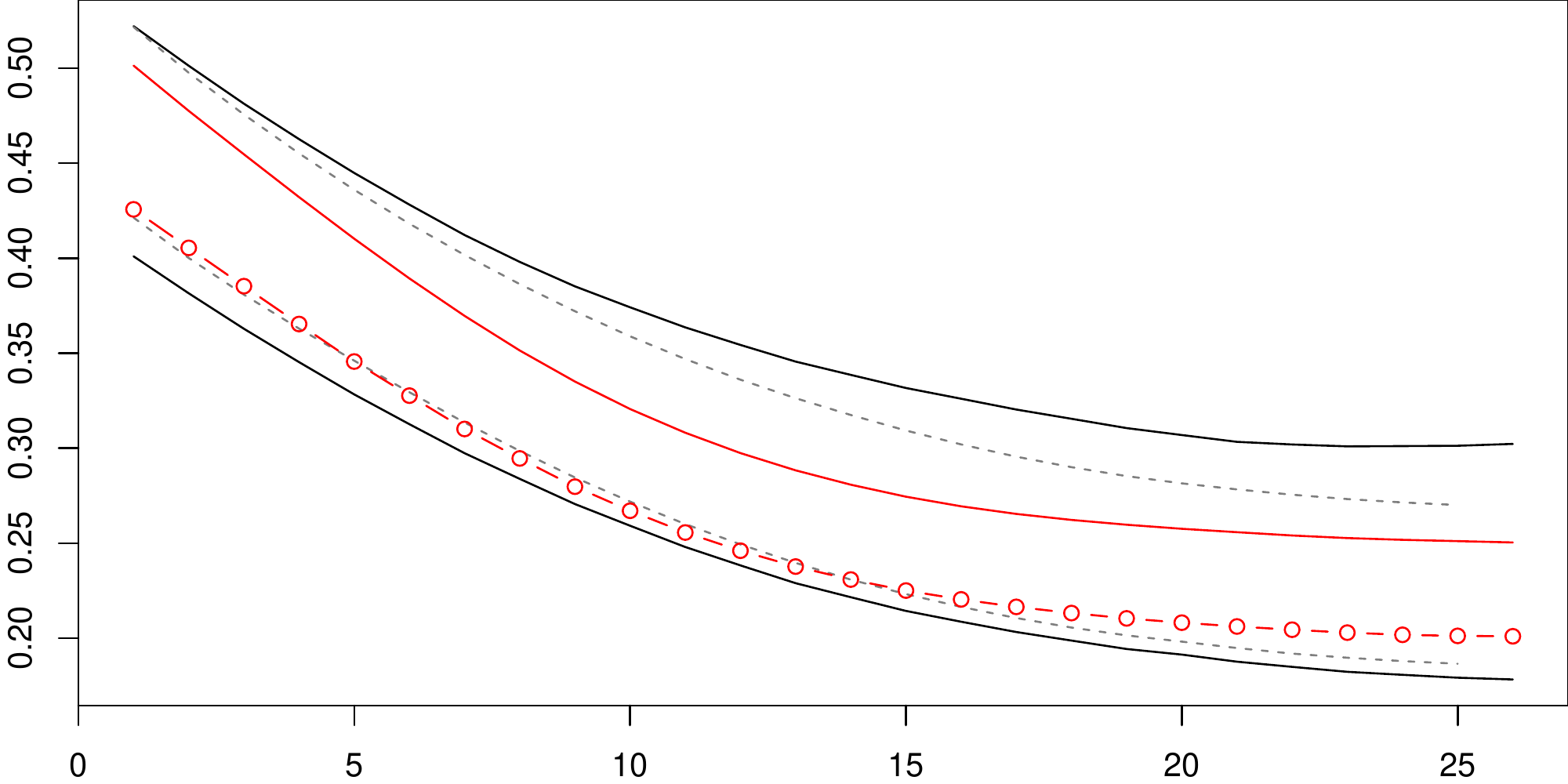} \\ \includegraphics[angle=0,scale=.35]{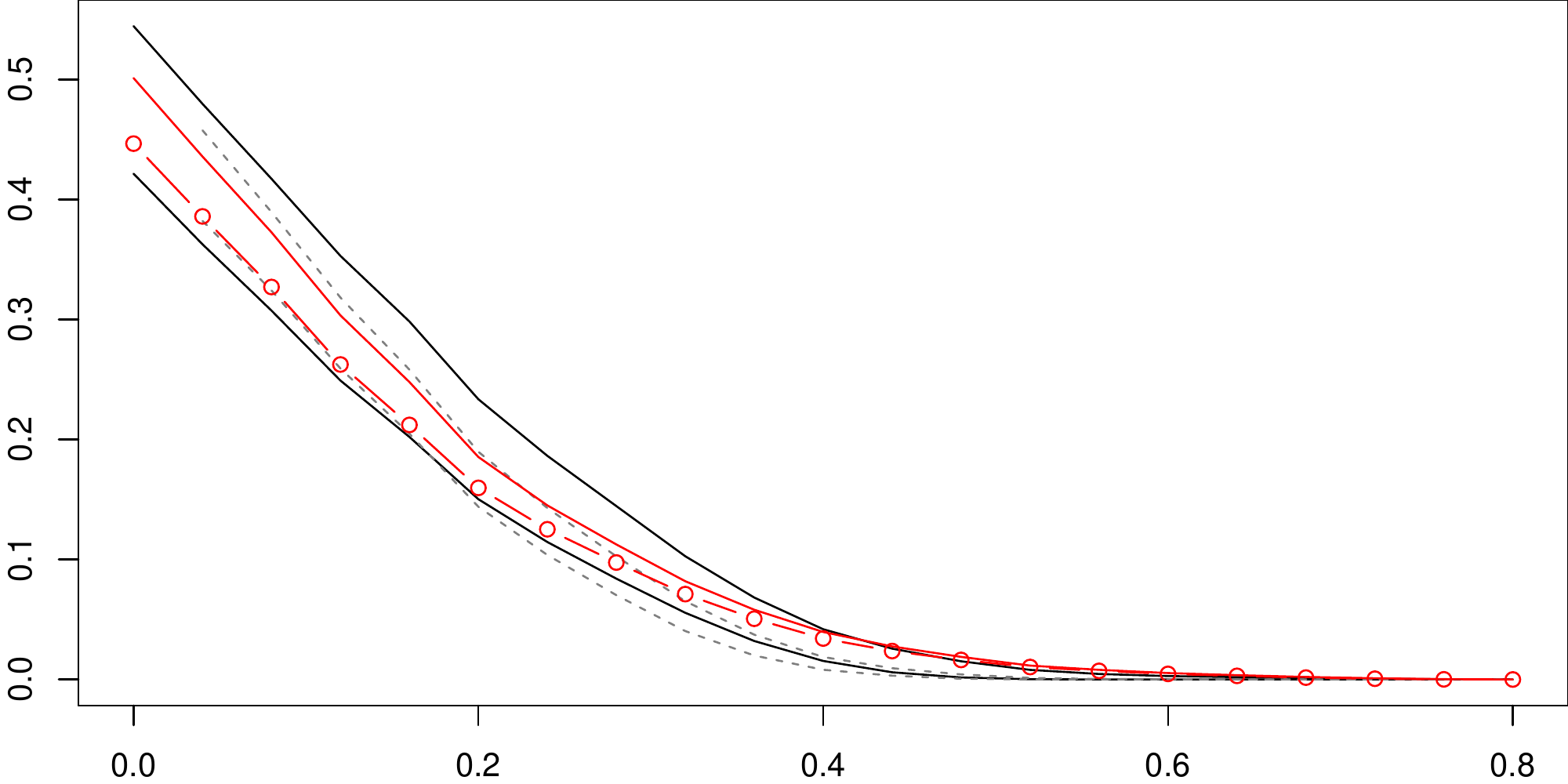} &
  \includegraphics[angle=0,scale=.35]{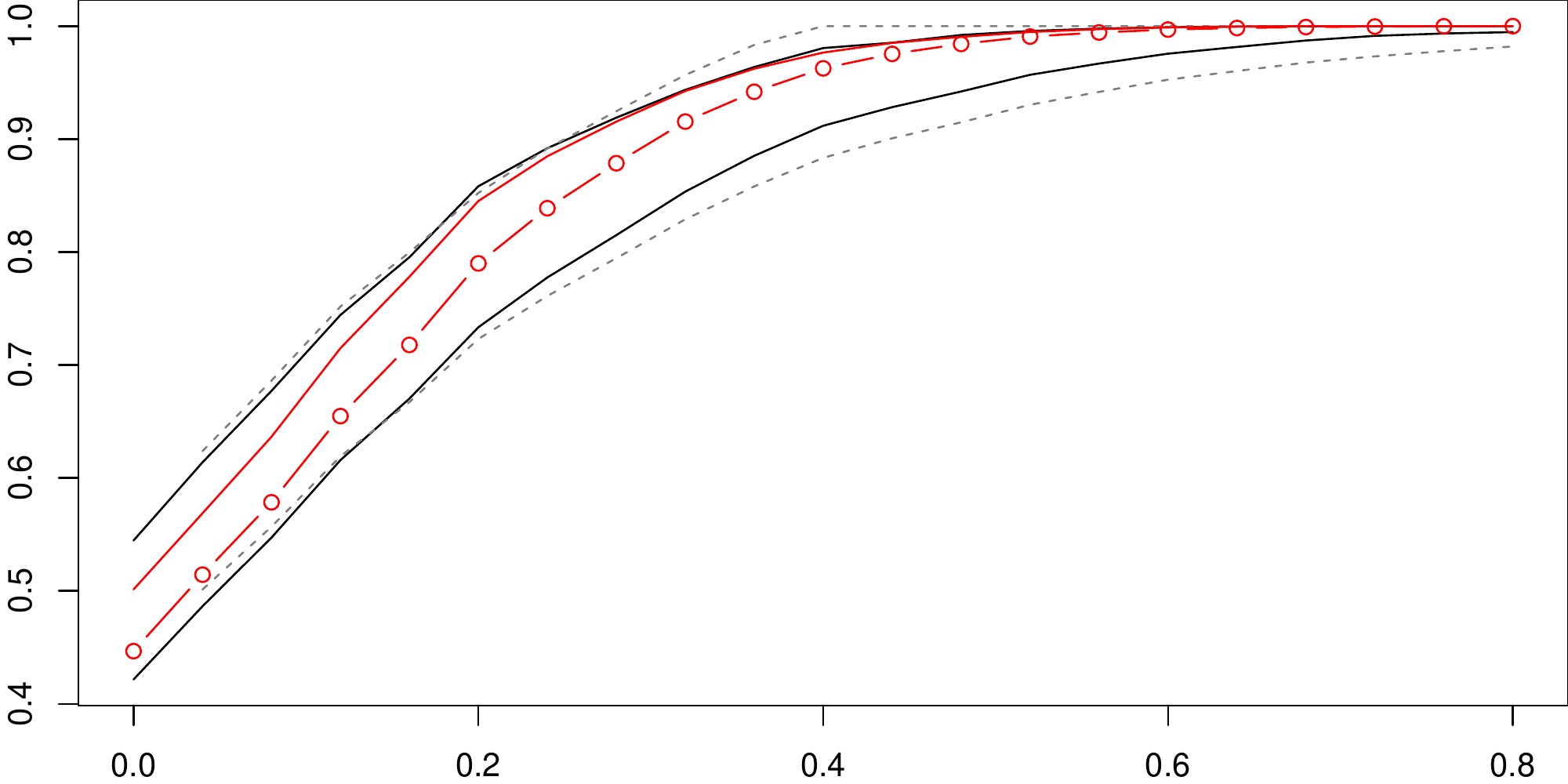} \\ \includegraphics[angle=0,scale=.35]{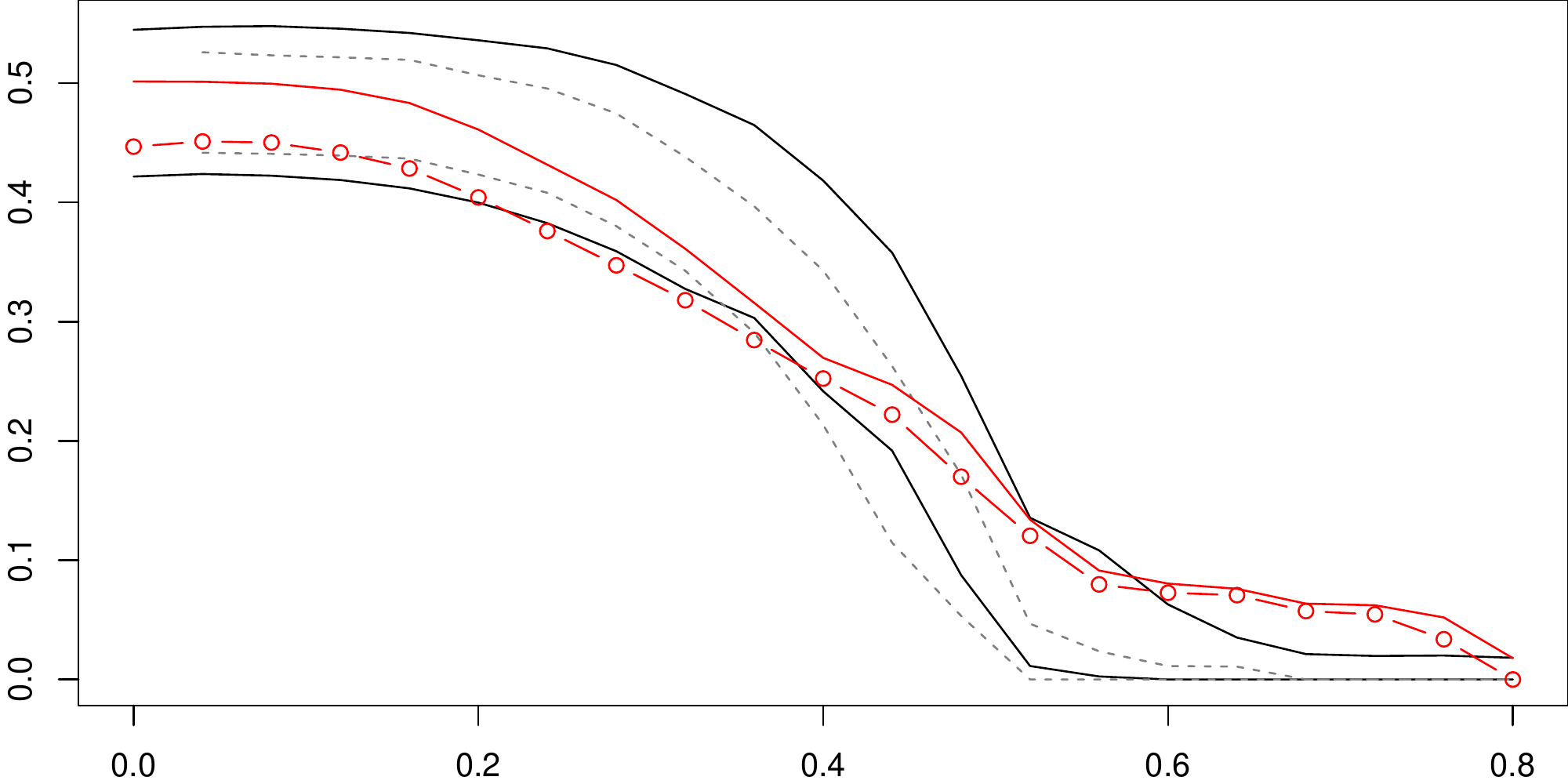} & \includegraphics[angle=0,scale=.35]{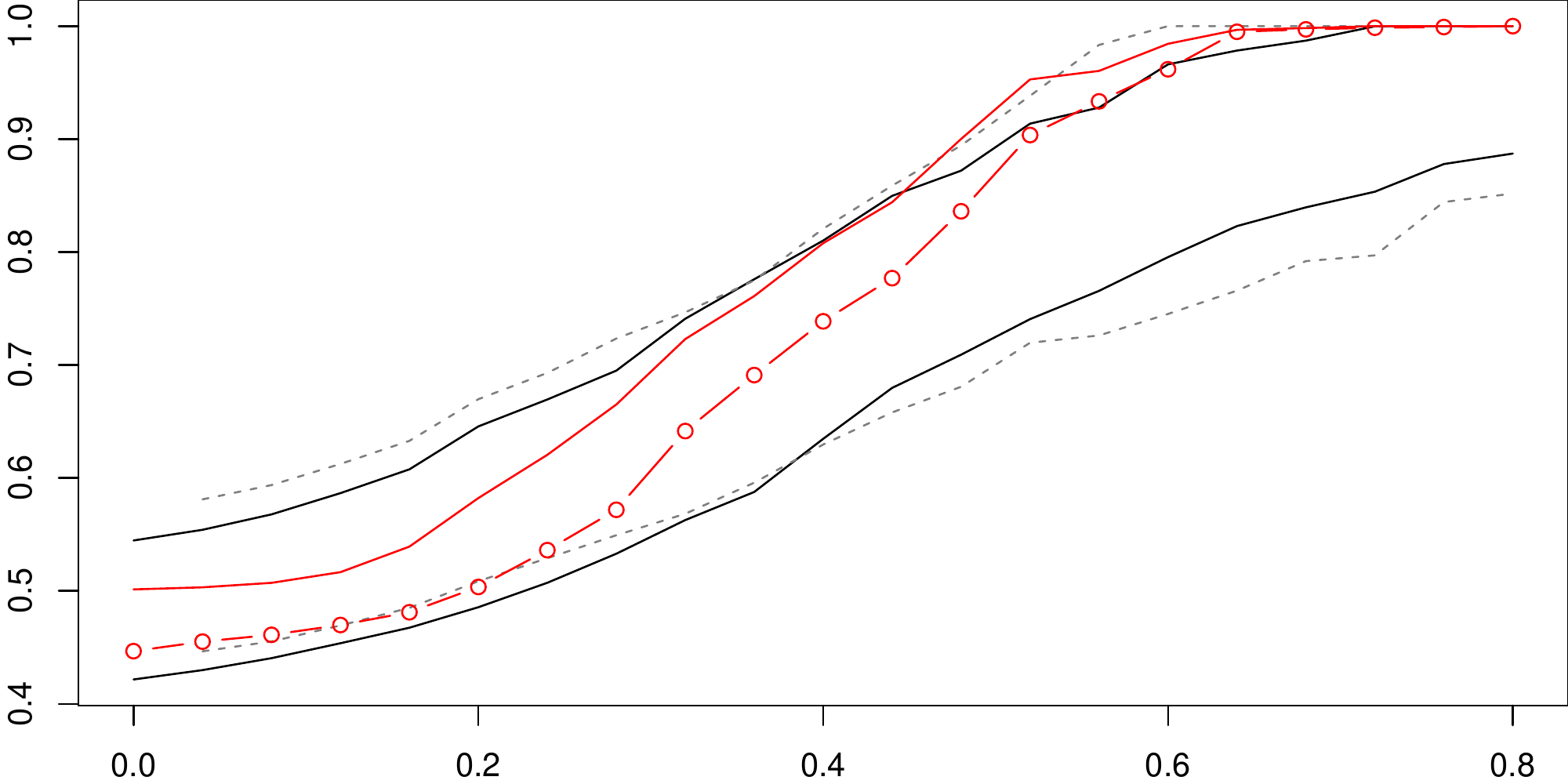}
   \end{tabular}
  }
  \caption{From top left to bottom right: contact distribution function, covariance function, erosion, dilatation, opening, closing, for the estimation from the heather dataset (solid red line), the estimation for the balls approximation of the heather dataset (red line with circles),  $95\%$ envelopes under  \eqref{ourestimation} (solid black line), $95\%$ envelopes  under  \eqref{energyJesper} (dashed  black line)}\label{fig:control}
 \end{figure}

Considering heather bushes as discs, the spherical germ-grain representation of the data seems to be a natural approximation. However, the independence between the location of the grains and their radii, implied by the Boolean model, appears more questionable. As a matter of fact, simulations of the fitted Boolean models from the aforementioned studies do not look visually similar to the heather dataset, as initially observed by P. Diggle in \cite{Diggle81}.  This lack of fit was confirmed in \cite{MH2} by various diagnostic plots (as in Figure~\ref{fig:control}). The same conclusion is drawn in \cite{Mrkvicka09}, where a Mat\'ern's cluster model is proposed as an alternative to the Boolean model. However no estimation procedure for this Mat\'ern's cluster model is given in  \cite{Mrkvicka09}. Instead, a brute-force approach is used to choose the parameters in order to pass some diagnostic plots. While the latter model appears to be well adapted  to the heather dataset, a comparison between our statistical method and the approach in  \cite{Mrkvicka09} is  difficult and therefore will not be conducted further.

Moreover,  a model of interacting discs has been fitted to the heather dataset  in \cite{MH2}. 
The model is very similar to the Quermass-interaction model, except that in the Hamiltonian \eqref{energy}, the Euler characteristic is replaced by the number of connected components. This model has been fitted in \cite{MH2} by a maximum likelihood approach as described in Section~\ref{MLE}, where three reference measures $\mu$ have been tested. As a result, the best fitted model was: $\mu$ is the uniform law on $[0,0.53]$, $\hat z=2.45$ and the estimated Hamiltonian, denoted $\widetilde H^{\hat\theta}$, is 
\begin{equation}\label{energyJesper}
\widetilde H^{\hat\theta}(\oo_\L)= 4.91\;  \A(\U_{\oo_\L}) - 1.18\; \LL(\U_{\oo_\L}) + 2.25\;   \mathcal N_{cc}(\U_{\oo_\L}),
\end{equation}   
where $\mathcal N_{cc}(\U_{\oo_\L})$ denotes the number of connected components of $\U_{\oo_\L}$.

As explained in Section~\ref{MLE}, the estimation procedure used in \cite{MH2} implies some rather strong restriction on the parameter space. For comparison, we fit hereafter a Quermass-interaction model  using the TF approach.

 Following the practical recommendations in Section~\ref{PA}, the binary image is first approximated by a union of balls, using the procedure described in \cite{Thiedmann}. The result of this approximation is shown on the right plot of Figure~\ref{fig:heather}. Second, we choose for the reference measure $\mu$ the uniform law on $[0.05,0.55]$, which is in agreement with the final choice made in \cite{MH2}, except that a minimal value for the radii is fixed to 0.05, corresponding to the smallest radius of the spots observed in the heather dataset.

The estimation of the full Quermass model (i.e. all parameters $z$, $\t_1$, $\t_2$, $\t_3$ are unknown) did not give satisfactory results. This  is not very surprising in view of the strong variability of the estimates shown in Figure~\ref{fig-quermass}. The contrast function in \eqref{Takacs} computed from the heather dataset actually exhibits in this general case many local minima, leading to an identifiability issue. Actually, it appears that any value of the intensity parameter $z$ can somehow be compensated by some value of the area interaction parameter $\t_1$. This identifiability issue between the area interaction and the intensity was already visible in the top-left plot of Figure~\ref{fig-area}, where we observe a strong correlation between the two estimates in the $\A$-process. Therefore, we have decided to choose a simpler model. Since the heather dataset seems not rich enough to distinguish the role played by the intensity parameter $z$ and the area interaction parameter $\t_1$, a natural choice is to fit an $(\LL,\x)$-process. In fact, fitting an $(\A,\LL)$-process leads to the same kind of identifiability issue as  before.  This basically means that only one parameter becomes relevant  in the $(\A,\LL)$-process for the heather dataset, namely the perimeter interaction parameter, which is not sufficient to obtain a good fit. 

So, the estimation of the $(\LL,\x)$-process has been implemented on the balls approximations of the heather dataset (right plot of Figure~\ref{fig:heather}), using the TF procedure associated to the test functions $f_0$, $f_{\alpha_1}$, $\cdots$, $f_{\alpha_{10}}$, $f_{\text{sum}}$ and  $f_{\text{iso}}$, where $\alpha_i=0.005\,i$, $i=1,\dots,10$, and  where $\mu$ is the uniform law on $[0.05,0.55]$. As a result, we obtained $\hat z=2.12$ and 

\begin{equation}\label{ourestimation}
H^{\hat\theta}(\oo_\L)=  0.14\; \LL(\U_{\oo_\L}) + 0.22\;   \x(\U_{\oo_\L}).
\end{equation} 

In order to check the quality of fit of \eqref{ourestimation} to the heather dataset, we use some diagnostic plots in Figure~\ref{fig:control}, where a comparison is also provided with the fitted model \eqref{energyJesper} of \cite{MH2}. These plots are the same as in \cite{MH2}, they correspond to  an estimation of the contact distribution function (top left): $H(r)=P(D\leq r|D>0)$ where $D=\inf\{\rho\geq 0: \U\cap B(0,\rho)\not=\emptyset\}$; the covariance function (top right)  $C(r)=P(x \in\U, y\in\U, \|x-y\|=r)$; and some shape-characteristics (see \cite{Ripley88} or \cite{MH2} for a definition): erosion $e_r$ (middle left), dilatation $d_r$ (middle right), opening $o_r$ (bottom left) and closing $c_r$ (bottom right). 

\begin{figure}[htbp]
  \setlength{\tabcolsep}{0.1cm} \centerline{
  \begin{tabular}[]{cc}
  \includegraphics[angle=0,scale=.3]{heather_balls.pdf} &
 \includegraphics[angle=0,scale=.3]{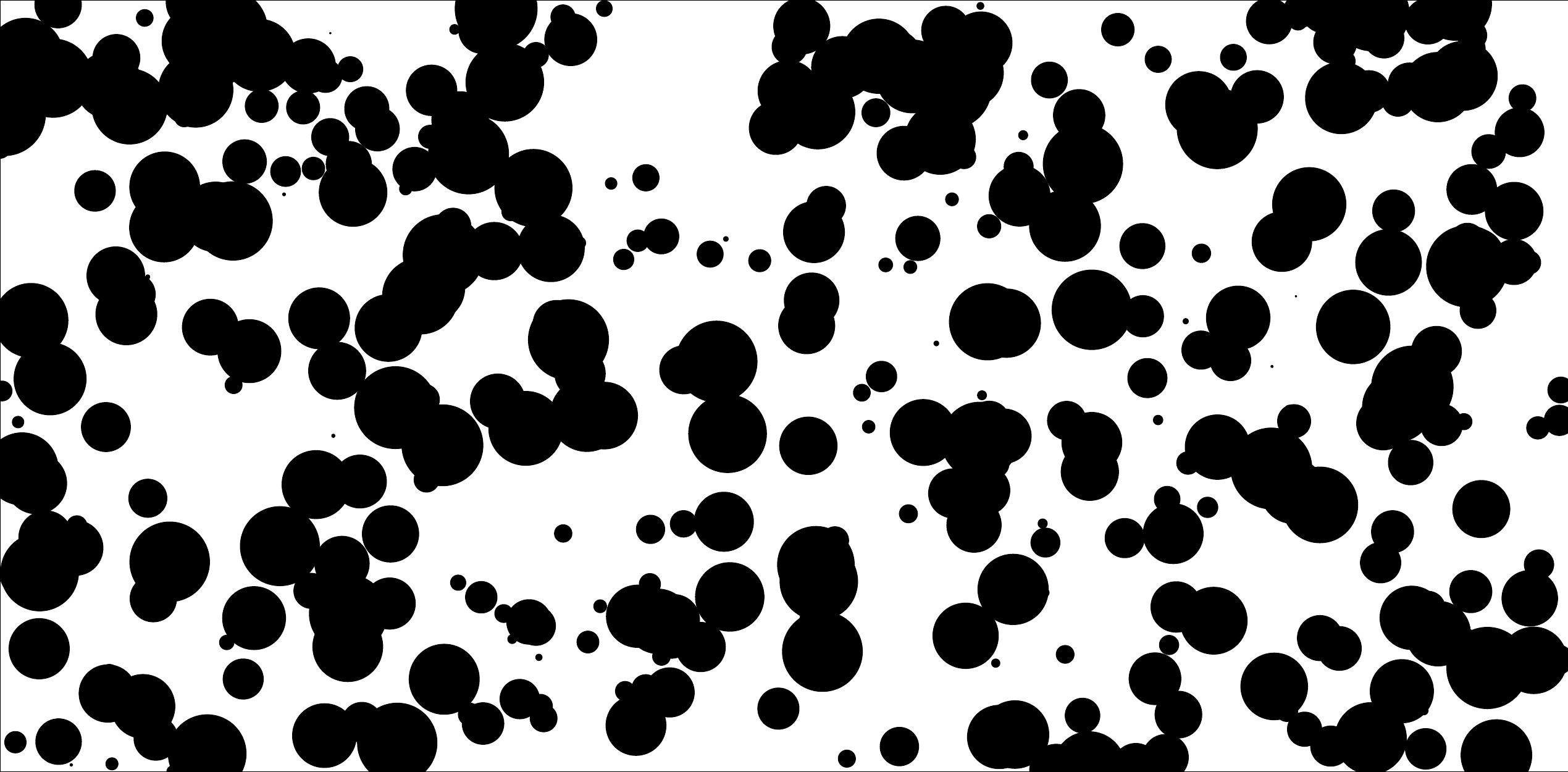} \\
  \includegraphics[angle=0,scale=.3]{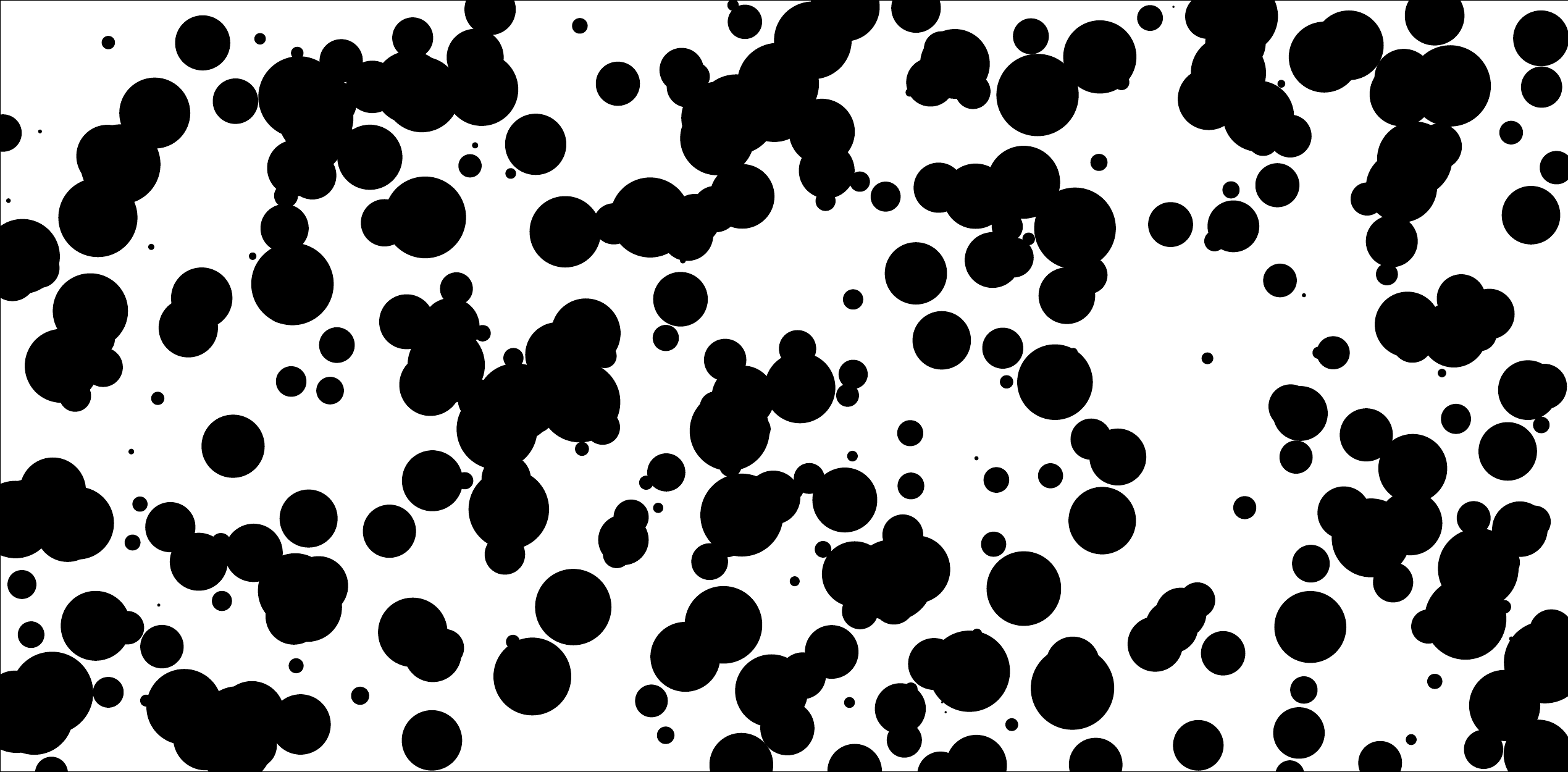} &
  \includegraphics[angle=0,scale=.3]{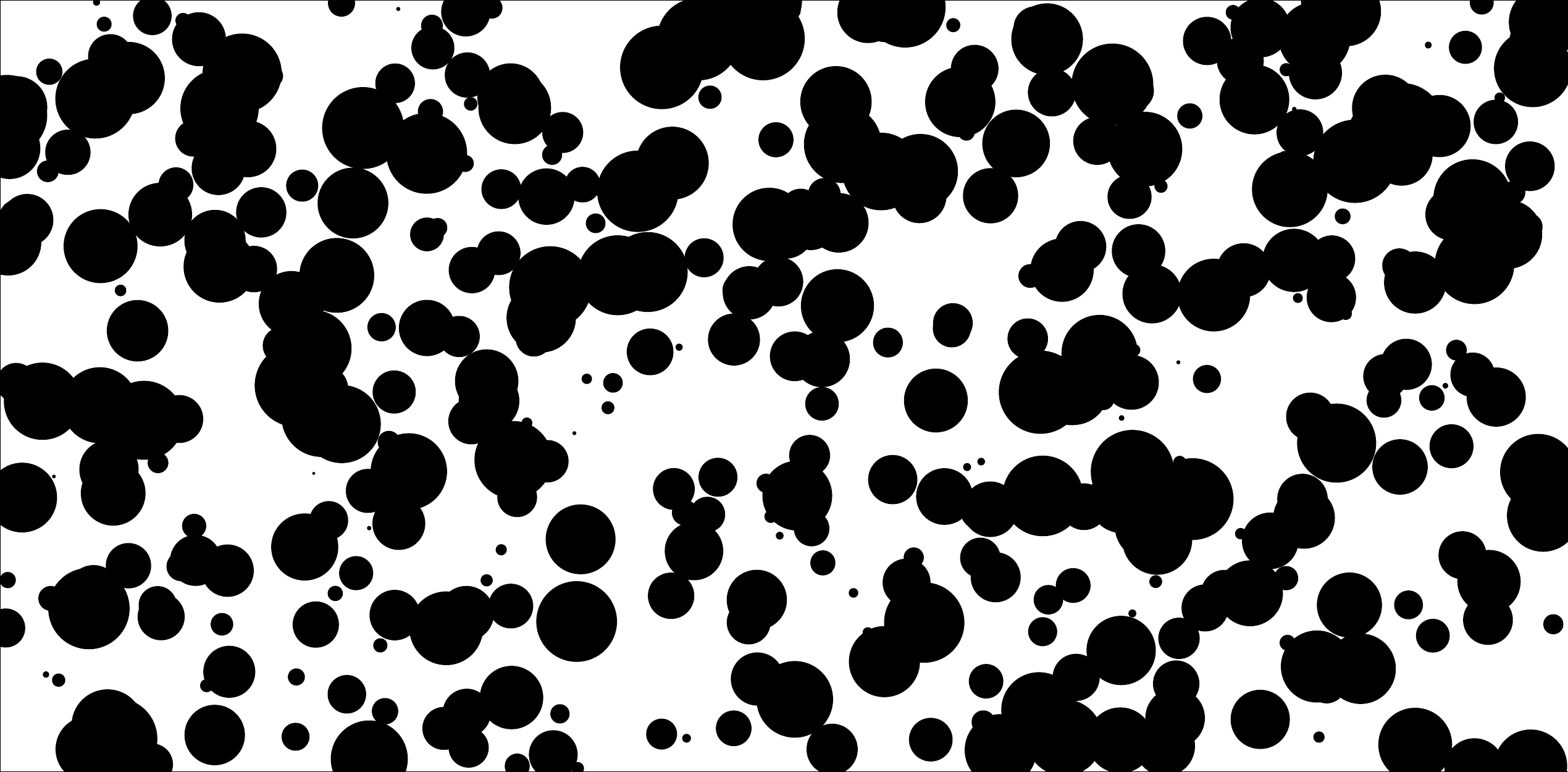} \\
   \end{tabular}
  }
  \caption{Heather data approximation (top left) and three samples from the fitted model \eqref{energyJesper}}\label{fig:samples_Jesper}
 \end{figure}

\begin{figure}[htbp]
  \setlength{\tabcolsep}{0.1cm} \centerline{
  \begin{tabular}[]{cc}
  \includegraphics[angle=0,scale=.3]{heather_balls.pdf} &
 \includegraphics[angle=0,scale=.3]{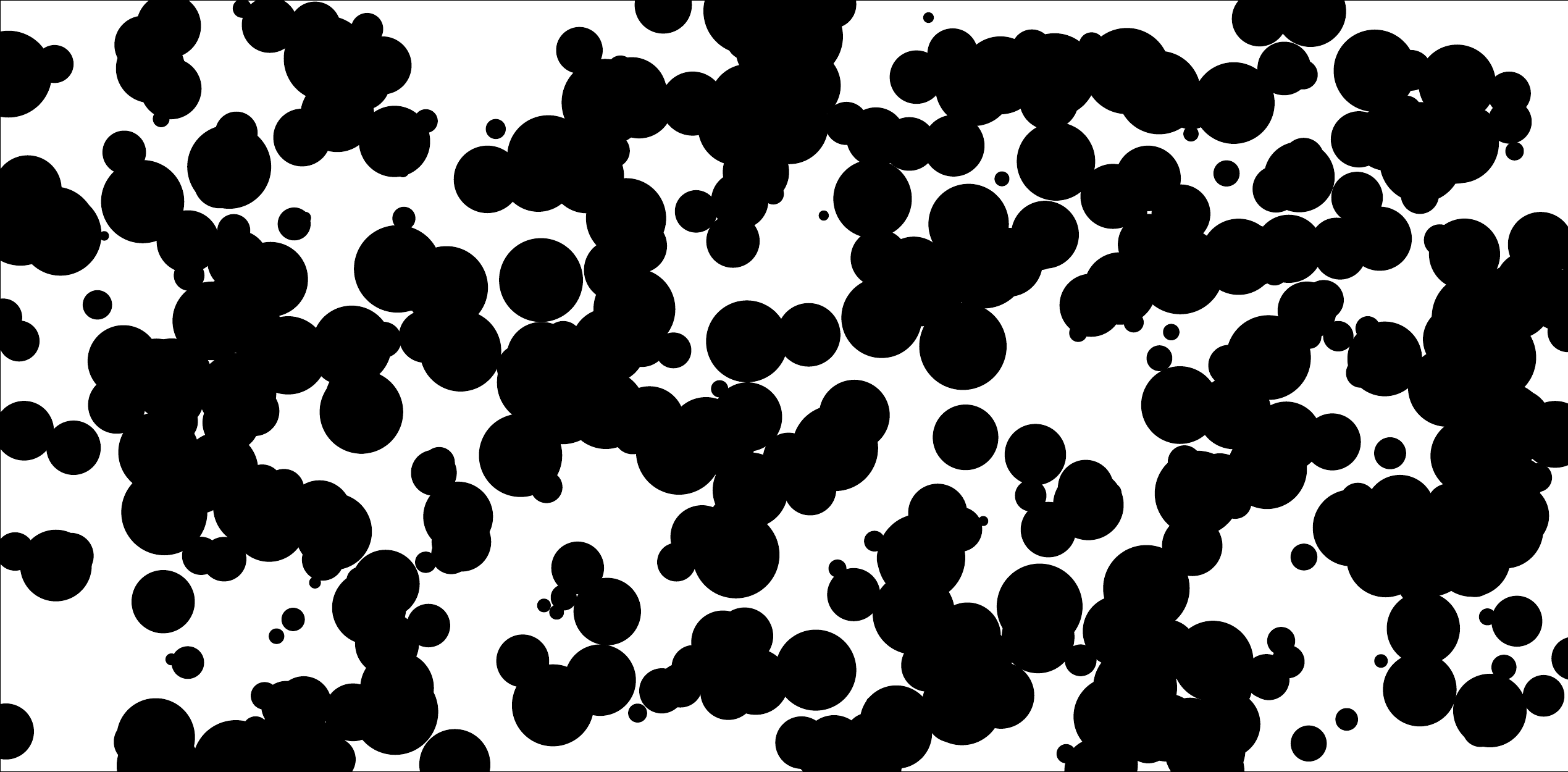} \\
  \includegraphics[angle=0,scale=.3]{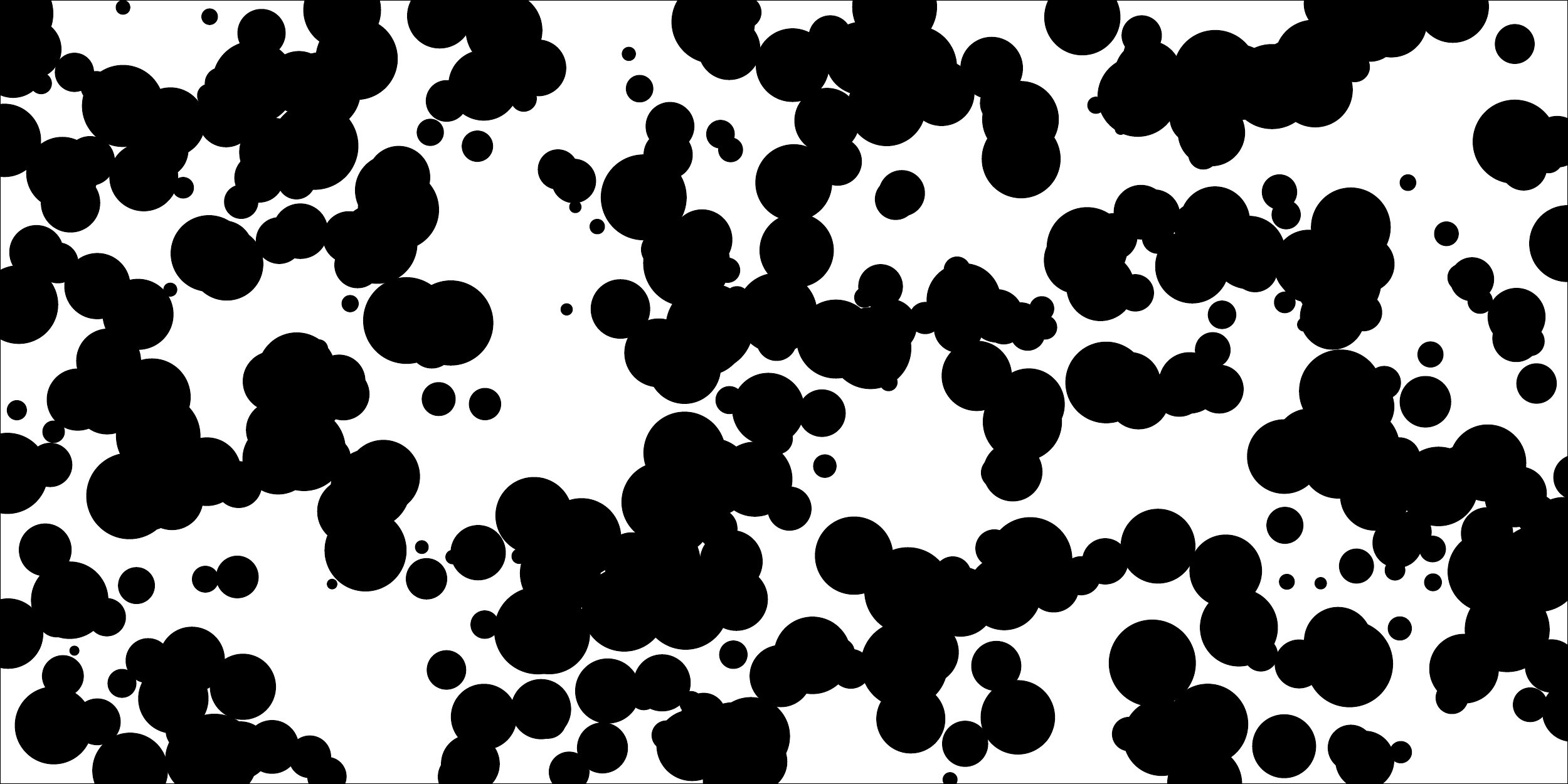} &
  \includegraphics[angle=0,scale=.3]{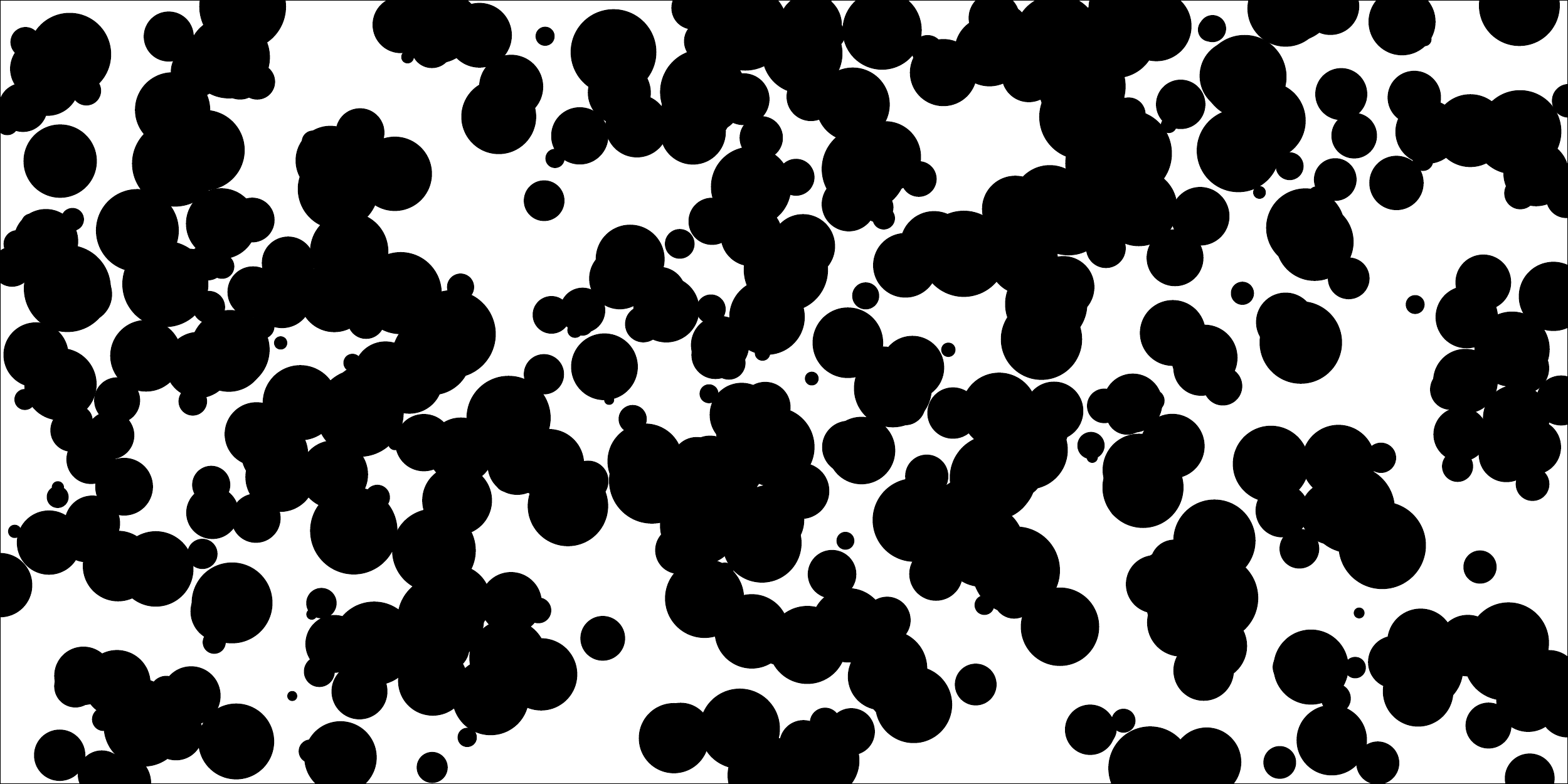} \\
   \end{tabular}
  }
  \caption{Heather data approximation (top left) and three samples from the fitted model \eqref{ourestimation}}\label{fig:heather_samples}
 \end{figure}

A slightly better fit is observed for \eqref{ourestimation}, in particular from the contact distribution function, the erosion and the closing, even if the closing is not fitted very well. This conclusion seems confirmed by the visual impression from samples of \eqref{energyJesper}  and  \eqref{ourestimation},  shown respectively in Figure~\ref{fig:samples_Jesper} and \ref{fig:heather_samples}, next to the balls approximation of the heather dataset.

\subsection*{Acknowledgments}
This work was partially supported by the  research project MATPYL of the F\'ed\'eration de Math\'ematiques des Pays de la Loire. Kate\v{r}ina {\sc Sta\v{n}kov\'{a} Helisov\'{a}}'s research was funded by Grant Agency of Czech Republic,  projects No. P201/10/0472 and 13-05466P.

\end{document}